\documentclass[11pt,a4]{article}
\usepackage[T2A]{fontenc}
\usepackage[cp866]{inputenc}
\usepackage[mathscr]{eucal}
\usepackage{amssymb}
\usepackage{amsmath}
\usepackage{fancyhea}
\renewcommand{\leq}{\leqslant}
\renewcommand{\geq}{\geqslant}

\renewcommand{\C}{{\mathbb C}}

\newcommand{\R}{{\mathbb R}}

\renewcommand{\k}{\rule{0.7em}{0.7em}}
\begin{document}
\sloppy
\headrulewidth = 2pt
\pagestyle{fancy}
\lhead[\scriptsize V.N.Gorbuzov,  A.F.Pranevich]
{\scriptsize V.N.Gorbuzov,  A.F.Pranevich}
\rhead[\it \scriptsize ${\mathbb R}\!$-holomorphic solutions and $\!\!{\mathbb R}\!$-differentiable integrals 
of multidimensional differential systems]{\it \scriptsize ${\mathbb R}\!$-holomorphic solutions and 
$\!\!{\mathbb R}\!$-differentiable integrals ...}
\headrulewidth=0.25pt

{\normalsize 

\thispagestyle{empty}

\mbox{}
\\[-2.5ex]
\centerline{
{
\bf 
${\mathbb R}\!$-HOLOMORPHIC SOLUTIONS AND 
$\!{\mathbb R}\!$-DIFFERENTIABLE INTEGRALS 
}
}
\\[0.75ex]
\mbox{}\hfill
{
\bf 
OF MULTIDIMENSIONAL DIFFERENTIAL SYSTEMS}
\hfill\mbox{}
\\[2.75ex]
\centerline{
\bf 
V.N. Gorbuzov$\!{}^{*},$   A.F. Pranevich$\!{}^{**}$ 
}
\\[2.75ex]
\centerline{
\it 
$\!{}^{*}\!\!$Department of Mathematics and Computer Science, 
Yanka Kupala Grodno State University,
}
\\[1ex]
\centerline{
\it 
Ozeshko 22, Grodno, 230023, Belarus
}
\\[1ex]
\centerline{
E-mail: gorbuzov@grsu.by
}
\\[2.5ex]
\centerline{
\it 
$\!{}^{**}\!\!$Department of Economics and Management, 
Yanka Kupala Grodno State University,
}
\\[1ex]
\centerline{
\it 
 Ozeshko 22, Grodno, 230023, Belarus
}
\\[1ex]
\centerline{
E-mail: pronevich@tut.by
}
\\[4.75ex]
\centerline{{\large\bf Abstract}}
\\[1ex]
\indent
We consider multidimensional differential systems (total differential systems  
and partial differential systems)
with $\!{\mathbb R}\!$-differentiable coefficients. 
We investigate the problem of the exis\-tence of $\!{\mathbb R}\!$-holomorphic solutions, 
$\!{\mathbb R}\!$-dif\-fe\-ren\-ti\-ab\-le in\-teg\-rals, and last multipliers.
The theo\-rem of existence and uniqueness of ${\mathbb R}\!$-holomorphic solution is proved.
The necessary con\-di\-ti\-ons and criteria for the existence of  ${\mathbb R}\!$-differentiable 
first integrals, partial integrals, and last multipliers are given. 
For a completely solvable total differential equation   
with ${\mathbb R}\!$-ho\-lo\-mor\-p\-hic right hand side are con\-s\-t\-ruc\-ted 
the classification of ${\mathbb R}\!$-singular points of solutions and 
proved sufficient conditions that equation have no movable nonalgebraical 
${\mathbb R}\!$-sin\-gu\-lar points.
The spectral method for building $\!{\mathbb R}\!$-differentiable first integrals 
for linear homogeneous first-order partial differential systems 
with ${\mathbb R}\!$-linear coefficients is developed.
\\[2ex]
\\[0.75ex]
\\[5ex]
{\large\bf Contents} 
\\[1ex]
{\bf  
1. System of total differential equations}\footnote[1]{
This Section  has been published in {\it Vestsi Nats. Akad. Navuk Belarusi}
(Ser. fiz.-matem. Navuk), 1999, No.\! 3, 124-126;\!
{\it Differential\! Equations and Control Processes} (http://www.neva.ru/journal), 2008, No.\! 1, 35-49;
{\it Differential Equations} 35 (1999), No. 4,  447-452.
}
                                                                                                                    \hfill\ 2
\\[0.5ex]
\mbox{}\hspace{1em}
1.1. Introduction                                                                                             \dotfill\ 2
\\[0.5ex]
\mbox{}\hspace{1em}
1.2. ${\mathbb R}\!$-holomorphic functions                                                     \dotfill\ 3
\\[0.5ex]
\mbox{}\hspace{1em}
1.3. The Cauchy existence and uniqueness theorem 
for an ${\mathbb R}\!$-holomorphic solution                                                     \dotfill\ 4
\\[0.5ex]
\mbox{}\hspace{1em}
1.4. ${\mathbb R}\!$-differentiable integrals and last multipliers                            \dotfill\ 5
\\[0.5ex]
\mbox{}\hspace{2.75em}
1.4.1. ${\mathbb R}\!$-differentiable partial integrals                                              \dotfill\ 6
\\[0.5ex]
\mbox{}\hspace{2.75em}
1.4.2. ${\mathbb R}\!$-differentiable first integrals                                                 \dotfill\ 10
\\[0.5ex]
\mbox{}\hspace{2.75em}
1.4.3. ${\mathbb R}\!$-differentiable last multipliers                                               \dotfill\ 12
\\[0.5ex]
\mbox{}\hspace{1em}
1.5. $\!{\mathbb R}\!$-regular solutions of an algebraic equation have no movable 
\\
\mbox{}\hspace{2.75em}
nonalgebraic ${\mathbb R}\!$-singular point                                              \dotfill\ 13
\\[0.75ex]
\noindent
{\bf  
2. 
System of first-order partial differential equations}\footnote[2]{
Section 2 has been published in  {\it Vestnik of the Yanka Kupala Grodno State Univ.}, 
2005, Ser.\! 2, No.\! 1(31), 45-52;
Proceedings of Scientific Conference  "\textit{Herzen Readings}-2004" (Ed. Professor V.F. Zaitsev), 
Russian State Uni\-ver\-si\-ty of Pedagogic, Saint-Petersburg, Russia, 2004, 65-70. 
}                                                                                                                  \hfill\ 16
\\[0.5ex]
\mbox{}\hspace{1em}
2.1. $\!{\mathbb R}\!$-differentiable integrals and last multipliers
                                                                                                                                   \dotfill\ 16
\\[0.5ex]
\mbox{}\hspace{2.75em}
2.1.1. $\!(n-k_1,n-k_2)\!$-cy\-lin\-dricality  partial integrals                                           \dotfill\ 17
\\[0.5ex]
\mbox{}\hspace{2.75em}
2.1.2. $\!(n-k_1,n-k_2)\!$-cy\-lin\-dricality first integrals                                               \dotfill\ 20
\\[0.5ex]
\mbox{}\hspace{2.75em}
2.1.3. $\!(n-k_1,n-k_2)\!$-cy\-lin\-dricality last multipliers                                             \dotfill\ 21
\\[0.5ex]
\mbox{}\hspace{1em}
2.2. First integrals of  linear homogeneous system with $\!{\mathbb R}\!$-linear coefficients 
                                                                                                                              \dotfill\ 23
\\[0.5ex]
\mbox{}\hspace{2.75em}
2.2.1. $\!{\mathbb R}\!$-linear partial integral                                                                           \dotfill\ 23
\\[0.5ex]
\mbox{}\hspace{2.75em}
2.2.2. $\!{\mathbb R}\!$-differentiable  first integrals                                                   \dotfill\ 24
\\[0.35ex]
{\bf  References}                                                                                                  \dotfill\ 27
\\[5ex]
\centerline{
\large\bf  
1. System of total differential equations}
\\[1.25ex]
\indent
{\bf  1.1. Introduction}
\\[0.5ex]
\indent
Let us consider a system of total differential equations 
\\[1.75ex]
\mbox{}\hfill                               
$
dw = X_1(z,w)\,dz + X_2(z,w)\,d\,\overline{z}\,,
$
\hfill (1.1)
\\[2ex]
where  
\vspace{0.35ex}
$w=(w_1,\ldots,w_n)\in {\mathbb C}^n,\ z=(z_1,\ldots,z_m)\in {\mathbb C}^m;$ 
the entries of the $n\times m$ matrices  
$X_1(z,w)=\|X_{\xi  j}(z,w)\|$ and  
\vspace{0.5ex}
$X_2(z,w)=\|X_{\xi,m+j}(z,w)\|$ are ${\mathbb R}\!$-differentiable 
[1, pp. 33 -- 35; 2, p. 22] in a domain 
\vspace{0.35ex}
$G\subset {\mathbb C}^{m+n}$ scalar functions 
\vspace{0.35ex}
$X_{\xi l}\colon G\to \C,\ \xi=1,\ldots, n,\ l=1,\ldots, 2m;$
$dw=\mbox{colon}(dw_1,\ldots,dw_n),$ $dz=\mbox{colon}(dz_1,\ldots,dz_m),$ and
$d\,\overline{z}=\mbox{colon}(d\,\overline{z}_1,\ldots,d\,\overline{z}_m)$
\vspace{0.25ex}
are vector columns; 
the $\overline{z}_j$  are the complex conjugates of   $z_j, \ j=1,\ldots,m.$ 
\vspace{0.5ex}

The notion of an $\!\!{\mathbb R}\!$-differentiable function is consistent with the 
approach of  I.N.Vekua~[3] and G. N. Polo\v{z}ii [4]
in the case of one complex variable and \'{E}. I. Grudo [5] in 
the case of two com\-p\-lex variables. 
Let $u\colon V\to {\mathbb C}$ be a one variable ${\mathbb R}\!$-differentiable function on the domain $V\subset {\mathbb C}.$  
The function $u$ is {\it holomorphic} if  $\partial_{{}_{\scriptstyle \overline{z}}}\,u(z)=0$ for all  $z\in V.$ 
\vspace{0.35ex}
The function  $u$ is called {\it an\-ti\-ho\-lo\-mor\-p\-hic} if  
$\partial_{z}u(z)=0$ for all  $z\in V$  [1, p. 42]. 
\vspace{0.35ex}
The function  $u$  is said to be $(p,q)\!$-{\it ana\-ly\-ti\-cal} if   
$(p(x,y)-iq(x,y))
\partial_{{}_{\scriptstyle \overline{z}}}\,{\rm Re}\,u(z) +
i\partial_{{}_{\scriptstyle \overline{z}}}\,{\rm Im}\,u(z) =0
$
\vspace{0.5ex}
for all $(x,y)\in V$  and $p(x,y)>0$  for all $(x,y)\in V$  [4]. 
\vspace{0.35ex}
If  $\partial_{{}_{\scriptstyle \overline{z}}}\,u(z) +A(z)u(z)+ B(z)\,\overline{u}(z) =C(z)$
for all  $z\in V,$ then we say that $u$  is a {\it generalized analytic} function [3]. 
In the case of several complex 
va\-ri\-ab\-les, 
the theorem of existence and uniqueness 
of ${\mathbb R}\!$-holomorphic solution for first-order partial differential system was proved in [6].
The spectral method for building first integrals of  
completely solvable multidimensional 
$\!\!{\mathbb R}\!$-linear differential systems was elaborated\! [7 -- 9].

In this paper we study the problem of the existence of 
$\!{\mathbb R}\!$-holomorphic solutions, 
$\!{\mathbb R}\!$-dif\-fe\-ren\-ti\-ab\-le first integrals, partial integrals,
and last multipliers for total differential systems (Section 1)
and partial differential systems (Section 2). 
The article is organized as follows. 

In Subsection 1.2 we define the basic notions of 
$\!{\mathbb R}\!$-holomor\-p\-hic functions and $\!{\mathbb R}\!$-singular points.
\!There we also investigate some relations between them.

In Subsection 1.3 we consider the completely solvable total differential system  (1.1) 
with ${\mathbb R}\!$-holomorphic right hand side. The theorem of existence and uniqueness 
of ${\mathbb R}\!$-holomorphic solution (analogous to the Cauchy theorem) is proved.

In Subsection 1.4 we investigate the problem of the existence of 
${\mathbb R}\!$-differentiable integrals and 
last multipliers for the system of total differential equations (1.1). The necessary 
con\-di\-tions and criteria for the existence of  ${\mathbb R}\!$-differentiable 
first integrals, 
${\mathbb R}\!$-differentiable partial integrals, and 
${\mathbb R}\!$-differentiable last multipliers are given.

In Subsection 1.5, for a completely solvable total differential equation  
with ${\mathbb R}\!$-ho\-lo\-mor\-p\-hic right hand side are constructed 
the classification of ${\mathbb R}\!$-singular points of solutions and 
proved sufficient conditions that equation have no movable nonalgebraical 
${\mathbb R}\!$-sin\-gu\-lar points (analogous to the Painleve theorem and Fuchsian's theorem). 

In Subsection 2.1 the necessary con\-di\-tions and criteria for the existence of  ${\mathbb R}\!$-differentiable 
integrals and last multipliers of linear homogeneous partial differential systems are given.

In Subsection 2.2 the spectral method for building $\!{\mathbb R}\!$-differentiable first integrals 
for linear homogeneous first-order partial differential systems 
with ${\mathbb R}\!$-linear coefficients is developed.
\\[2.75ex]
\indent
{\bf 1.2. ${\mathbb R}\!$-holomorphic functions}
\\[1.25ex]
\indent
{\bf Definition 1.1.}\!
{\it 
A function $g\colon {\mathscr D}\to {\mathbb C},\, {\mathscr D}\subset {\mathbb C}^m,$
is said to be ${\mathbb R}\!$-holomorphic at a point 
$z_0\!=\!(z_1^0,\ldots, z_m^0)\!\in {\mathscr D}\!$ if there exists a neighborhood 
$U(z_0)\!\subset\! {\mathscr D}\!$ of the point $\!z_0\!$ such that in this neighborhood the 
function $\!g\!$ can be represented by the absolutely convergent function se\-ries}
\\[1.75ex]
\mbox{}\hfill
$
\displaystyle 
g(z)=\sum_{k_1+l_1+\ldots+k_m+l_m=0}^{{}+\infty}
c_{{}_{\scriptstyle k_1l_l\ldots k_m l_m}}
\prod_{j=1}^m 
\bigl(z_j-z_j^0\bigr)^{k_j}\bigl(\overline{z}_j-\overline{z}_j^{\,0}\bigr)^{l_j}$ 
\ \ for all 
$ z\in U(z_0),
$
\hfill (1.2)
\\[2ex]
{\it where $c_{{}_{\scriptstyle k_1l_l\ldots k_m l_m}}\in {\mathbb C}$ 
and the exponents $k_j$ and $l_j$ are nonnegative integers.
}
\vspace{1.25ex}

The term {\it ${\mathbb R}\!$-holomorphic} is introduced by analogy with the term 
{\it ${\mathbb R}\!$-differentiable}. Indeed, it follows from the absolute convergence 
of the series (1.2) that the real and imaginary parts in the representation $g=u+iv$ of an 
${\mathbb R}\!$-holomorphic function are real holomorphic functions in a neighborhood 
of the point $(x_0,y_0),$ where $x_0={\rm Re}\,z_0$ and $y_0={\rm Im}\, z_0.$ 
\vspace{1ex}

{\bf Definition 1.2.}\!
{\it 
A function $\overline{g} \colon {\mathscr D}\to {\mathbb C},\ 
{\mathscr D}\subset {\mathbb C}^m,$
is said to be conjugate to the ${\mathbb R}\!$-ho\-lo\-mor\-p\-hic 
function {\rm (1.2)} at a point $z_0\in {\mathscr D}$ if in some neighborhood 
$U(z_0)\subset {\mathscr D}$ of the point $z_0$ the 
function $\overline{g}$ can be represented by the function se\-ries}
\\[1.75ex]
\mbox{}\hfill
$
\displaystyle 
\overline{g}(z)=\sum_{k_1+l_1+\ldots+k_m+l_m=0}^{{}+\infty}
\overline{c}_{{}_{\scriptstyle k_1l_l\ldots k_m l_m}}
\prod_{j=1}^m 
\bigl(\overline{z}_j-\overline{z}_j^{\,0}\bigr)^{k_j}\bigl(z_j-z_j^0\bigr)^{l_j}$ 
\ \ for all 
$ z\in U(z_0).
$
\hfill (1.3)
\\[2ex]
\indent
We can readily see that this is well defined, since the sets of absolute convergence of the 
function series (1.2) and (1.3) coincide; moreover, $\overline{g}$ is 
${\mathbb R}\!$-ho\-lo\-mor\-p\-hic at the point $z_0.$

Since an ${\mathbb R}\!$-ho\-lo\-mor\-p\-hic function of $m$ independent variables
$z_j,\ j=1,\ldots, m,$ can be obtained from a holomorphic fuction of $2m$ 
independent variables $u_j$ and $v_j,\ j=1,\ldots, m,$ via the correspondence
\\[1ex]
\mbox{}\hfill
$
u_j\mapsto z_j, \ \ \ 
v_j\mapsto \overline{z}_j, 
\quad j=1,\ldots,m, 
$
\hfill (1.4)
\\[1.5ex]
we have the following assertions.
\vspace{0.75ex}

{\bf Proposition 1.1.}\!
{\it
Let functions $g_1\colon {\mathscr D}\to {\mathbb C}$ and 
$g_2\colon {\mathscr D}\to {\mathbb C},\, {\mathscr D}\subset {\mathbb C}^m,\!$ 
be ${\mathbb R}\!$-ho\-lo\-mor\-p\-hic 
at a point $z_0\in {\mathscr D}.$ Then the relations 
\\[1.75ex]
\mbox{}\hfill
$
\overline{g_1(z)+g_2(z)}=\overline{g}_1(z)+\overline{g}_2(z),
\qquad 
\overline{g_1(z)\cdot g_2(z)}=\overline{g}_1(z)\cdot \overline{g}_2(z),
\hfill
$
\\
\mbox{}\hfill {\rm (1.5)}
\\
\mbox{}\hfill
$
\overline{D_{{}_{\scriptstyle z_j}}\,g_1(z)}=
D_{{}_{\scriptstyle \overline{z}_j}}\,\overline{g}_1(z),
\quad \  
\overline{D_{{}_{\scriptstyle \overline{z}_j}}\,g_1(z)}=
D_{{}_{\scriptstyle z_j}}\,\overline{g}_1(z),
\ \ \ j=1,\ldots, m,
\hfill
$
\\[1.5ex]
are valid in some heighborhood $U(z_0)\subset {\mathscr D}$ of the point $z_0.$
}
\vspace{0.5ex}

{\bf Proposition 1.2} [10, p. 33].
{\it
If a function that is ${\mathbb R}\!$-ho\-lo\-mor\-p\-hic 
in a domain ${\mathscr D}\subset {\mathbb C}^m$ identically vanishes in some 
neighborhood $U\subset {\mathscr D},$ then this function identically vanishes 
in the entire domain ${\mathscr D}.$
}
\vspace{0.5ex}

{\bf Corollary 1.1.}
{\it
If two functions ${\mathbb R}\!$-ho\-lo\-mor\-p\-hic 
in a domain ${\mathscr D}\subset {\mathbb C}^m$ coincide in some 
neighborhood $U\subset {\mathscr D},$ then they coincide  
in the entire domain ${\mathscr D}.$
}
\vspace{0.5ex}

This corollary allows one to use the method of ${\mathbb R}\!$-ho\-lo\-mor\-p\-hic 
continuation for an ${\mathbb R}\!$-ho\-lo\-mor\-p\-hic function and hence consider 
multivalued ${\mathbb R}\!$-ho\-lo\-mor\-p\-hic functions. 
\vspace{0.5ex}

{\bf Definition 1.3.}\!
{\it 
An ${\mathbb R}\!$-holomorphic function 
$g\colon {\mathscr D}\to {\mathbb C},\ {\mathscr D}\subset {\mathbb C}^m,$
is said to be ${\mathbb R}\!$-regular at a point $z_0\in {\mathscr D}$ if 
}
\\[1.5ex]
\mbox{}\hfill
$
\displaystyle 
\text{rank} 
\left\|\!
\begin{array}{cccccc}
D_{{}_{\scriptstyle z_1}}g(z_0) & \ldots & D_{{}_{\scriptstyle z_m}}g(z_0) & 
D_{{}_{\scriptstyle \overline{z}_1}}g(z_0) & \ldots & 
D_{{}_{\scriptstyle \overline{z}_m}}g(z_0)
\\[1.25ex]
D_{{}_{\scriptstyle z_1}}\overline{g}(z_0) & \ldots & 
D_{{}_{\scriptstyle z_m}}\overline{g}(z_0) & 
D_{{}_{\scriptstyle \overline{z}_1}}\overline{g}(z_0) & \ldots & 
D_{{}_{\scriptstyle \overline{z}_m}}\overline{g}(z_0)
\end{array}
\!\right\| =2;
\hfill
$
\\[1.75ex]
{\it otherwise, it is said to be ${\mathbb R}\!$-singular}.
\vspace{0.5ex}

The possibility of ${\mathbb R}\!$-holomorphic continuation allows one to consider 
${\mathbb R}\!$-singular points, that is, points in a neighborhood of which an 
${\mathbb R}\!$-holomorphic function does not admit an 
${\mathbb R}\!$-holomorphic continuation.

Let an $\!{\mathbb R}\!$-holomorphic function 
$\!g\colon\! {\mathscr D}\!\to {\mathbb C}$ take the value $g(a)\!=g^a$ at a point 
$\!a\!\in {\mathscr D}\!\subset\! {\mathbb C}^m\!\!$ and satisfy the equation 
$\Phi(g,z)=0,$ where $\Phi$ is an ${\mathbb R}\!$-holomorphic function of its 
arguments in the neighborhood of the point $(g^a,a)\in V\subset {\mathbb C}^{m+1};$ 
moreover, $\Phi(g,a)\not\equiv 0.$ 
The point $a$ is referred to as an 
{\it algebraic critical ${\mathbb R}\!$-singular point} of the function $g$ if 
\\[1ex]
\mbox{}\hfill
$
\bigl|\partial_{{}_{\scriptstyle g}}\Phi(g^a,a)\bigr|^2\,-\,
\bigl|\partial_{\,{}_{\scriptstyle \overline{g}}}\,\Phi(g^a,a)\bigr|^2=\,0
\hfill
$
\\[1.25ex] 
and the function $g$ is not ${\mathbb R}\!$-holomorphic at the point $a.$

Suppose that an ${\mathbb R}\!$-holomorphic function 
$g\colon {\mathscr D}\to {\mathbb C}$ has the form $g(z)=1/f(z)$ and $f(a)=0;$ 
in this case, the following definitions will be used: 
1) if the function $f$ is ${\mathbb R}\!$-ho\-lo\-mor\-p\-hic at the point $a,$ then this 
point is referred to as an {\it ${\mathbb R}\!$-pole} of the function $g;$
2) if $a$ is an algebraic critical ${\mathbb R}\!$-singular point of the function $f,$ 
then this point  is referred to as a {\it critical ${\mathbb R}\!$-pole} of the function $g.$
Algebraic critical ${\mathbb R}\!$-singular points, ${\mathbb R}\!$-poles, and 
critical ${\mathbb R}\!$-poles are referred to as {\it algebraic ${\mathbb R}\!$-singular points}.

Let a point $a$ be a nonalgebraic $\!{\mathbb R}\!$-singular point of an 
${\mathbb R}\!$-holomorphic function $\!g\colon\! {\mathscr D}\to {\mathbb C}.$
In each plane $z_j$ we take the circle $|z_j-a_j|=r_j,\ j=1,\ldots, m,$
where $a=(a_1,\ldots,a_m).$ 
By $\triangle_r$ we denote the set of values that are taken by the function $g$ or to 
which it tends for its various ${\mathbb R}\!$-holomorphic continuations into 
the polydisk $|z_j-a_j|< r_j,\ j=1,\ldots, m.$ 
\vspace{0.35ex}
If $r_j\to 0,\ j=1,\ldots, m,$ then the set $\triangle_r$ tends to some limit set 
$\triangle_a g.$ If $\triangle_a g$ is a singleton, then the point $a$ is referred to as a 
{\it transcendental} ${\mathbb R}\!$-singular point of the function $g.$
If the set $\triangle_a g$ contains more than one point, then the point $a$ is referred to as 
a {\it $\triangle\!$-essentially} ${\mathbb R}\!$-singular point of the function $g.$

For example, the function 
$g\colon z_1\to z_1^2+z_1\overline{z}_1$ for all $z_1\in {\mathbb C}$ 
is ${\mathbb R}\!$-holomorphic 
on the entire complex plane ${\mathbb C}$ but is not holomorphic, since on ${\mathbb C}$ 
there is  no point in whose neighborhood the Cauchy-Riemann conditions are satisfied.
The point $z_1=0$ is

a) an algebraic ${\mathbb R}\!$-singular point for the function 
$g(z_1)=\sqrt{\overline{z}_1}\,z_1;$
\vspace{0.25ex}

b) an ${\mathbb R}\!$-pole for the function 
$g(z_1)=(1+z_1+\overline{z}_1)/z_1;$
\vspace{0.25ex}

c) a transcendental ${\mathbb R}\!$-singular point for the function 
$g(z_1)=\ln(\overline{z}_1+z_1^2);$
\vspace{0.25ex}

d) a $\triangle\!$-essentially ${\mathbb R}\!$-singular point for the function 
\vspace{0.25ex}
$g(z_1)=\exp (1/\overline{z}_1)\ (\triangle_0 g=\overline{{\mathbb C}}$ by analogy 
with the Sokhotskii theorem for an antiholomorphic function) and for 
\vspace{0.25ex}
$g(z_1)=z_1/\overline{z}_1$ (along any path 
$\!L_0\colon\! k\exp(i\omega_0),\, 0\leq k\!<\!+\infty,\!$ 
\vspace{0.25ex}
the function tends to $\exp (2i\omega_0),$
for $\!\omega_0\!\in\! [0; 2\pi)$ these limit values form the circle $|z_1|=1;$ therefore, 
$\triangle_0 g$ is not a singleton).
\\[3.75ex]
\indent
{\bf  \!\!1.3.\!\! The\! Cauchy\! existence and\! uniqueness theorem for an 
$\!\!{\mathbb R}\!\!$-holomorphic\! solution}
\\[1.75ex]
\indent
\!We assume that $X_{\xi l}\colon G\to {\mathbb C},\ \xi=1,\ldots, n,\ l=1,\ldots, 2m,$ 
are ${\mathbb R}\!$-holomorphic functions in the domain $G.$ Moreover, we consider 
system of total differential equations  (1.1) for the case in which it is completely solvable, 
i.e., the Frobenius conditions 
\\[1.5ex]
\mbox{}\hfill
$
\displaystyle
\partial_{{}_{\scriptstyle z_{{}_{\scriptsize\zeta }}}} X_{\tau j}(z,w) \, +\, 
\sum\limits_{\xi=1}^{n}
\bigl( X_{\xi\zeta}(z,w)\,
\partial_{{}_{\scriptstyle w_{{}_{\scriptsize\xi }}}} X_{\tau j}(z,w) \, +\, 
\overline{X}_{\xi,m+\zeta}(z,w)\,
\partial_{{}_{\scriptstyle
\overline{w}_{{}_{\scriptsize\xi }}}} X_{\tau j}(z,w)\bigr) =
\hfill
$
\\[1.5ex]
\mbox{}\hfill
$
\displaystyle
=\partial_{{}_{\scriptstyle z_{{}_{\scriptsize j}}}}
X_{\tau \zeta}(z,w) \, + \,
\sum\limits_{\xi=1}^{n}
\bigl( X_{\xi j}(z,w)\,
\partial_{{}_{\scriptstyle w_{{}_{\scriptsize\xi}}}}
X_{\tau \zeta}(z,w) \,  +\, 
\overline{X}_{\xi,m+j}(z,w)\,
\partial_{{}_{\scriptstyle
\overline{w}_{{}_{\scriptsize\xi }}}} X_{\tau \zeta}(z,w)\bigr),
\hfill
$
\\[2ex]
\mbox{}\hfill
$
\displaystyle
\partial_{{}_{\scriptstyle
\overline{z}_{{}_{\scriptsize\zeta}}}} X_{\tau,m+j}(z,w) \, +\, 
\sum\limits_{\xi=1}^{n}
\bigl( X_{\xi,m+\zeta}(z,w)\,
\partial_{{}_{\scriptstyle w_{{}_{\scriptsize\xi }}}} X_{\tau,m+j}(z,w) \, +\, 
\overline{X}_{\xi\zeta}(z,w)\,
\partial_{{}_{\scriptstyle
\overline{w}_{{}_{\scriptsize\xi }}}} X_{\tau,m+j}(z,w) \bigr) =
\hfill
$
\\
\mbox{}\hfill (1.6)
\\
\mbox{}\hfill
$
\displaystyle
=\partial_{{}_{\scriptstyle
\overline{z}_{{}_{\scriptsize j}}}} X_{\tau,m+\zeta}(z,w) \, +\, 
\sum\limits_{\xi=1}^{n}
\bigl( X_{\xi,m+j}(z,w)\,
\partial_{{}_{\scriptstyle w_{{}_{\scriptsize\xi }}}} X_{\tau,m+\zeta}(z,w) \, +\, 
\overline{X}_{\xi j}(z,w)\,
\partial_{{}_{\scriptstyle
\overline{w}_{{}_{\scriptsize\xi }}}} X_{\tau,m+\zeta}(z,w)\bigr),
\hfill
$
\\[2ex]
\mbox{}\hfill
$
\displaystyle
\partial_{{}_{\scriptstyle
z_{{}_{\scriptsize\zeta}}}} X_{\tau,m+j}(z,w) \, +\, 
\sum\limits_{\xi=1}^{n}
\bigl( X_{\xi\zeta}(z,w)\,
\partial_{{}_{\scriptstyle w_{{}_{\scriptsize\xi }}}} X_{\tau,m+j}(z,w) \, +\, 
\overline{X}_{\xi,m+\zeta}(z,w)\,
\partial_{{}_{\scriptstyle
\overline{w}_{{}_{\scriptsize\xi }}}} X_{\tau,m+j}(z,w) \bigr) =
\hfill
$
\\[1.5ex]
\mbox{}\hfill
$
\displaystyle
=\partial_{{}_{\scriptstyle
\overline{z}_{{}_{\scriptsize j}}}} X_{\tau\zeta}(z,w) \, +\, 
\sum\limits_{\xi=1}^{n}
\bigl( X_{\xi,m+j}(z,w)\,
\partial_{{}_{\scriptstyle w_{{}_{\scriptsize\xi }}}} X_{\tau\zeta}(z,w) \, + \, 
\overline{X}_{\xi j}(z,w)\,
\partial_{{}_{\scriptstyle
\overline{w}_{{}_{\scriptsize\xi }}}} X_{\tau\zeta}(z,w)\bigr)
\hfill
$
\\[1.5ex]
\mbox{}\hfill
\ \ for all $(z,w)\in G,
\ \ \tau=1,\ldots, n,
\, \ j=1,\ldots, m,
\ \, \zeta=1,\ldots, m,
\hfill
$
\\[1.5ex]
are satisfied [11 -- 13]. 
\vspace{0.75ex}

{\bf  Theorem 1.1.}\!
{\it
If the functions $X_{\xi l}\colon G\to {\mathbb C},\ \xi=1,\ldots, n,\ l=1,\ldots, 2m,$ 
are ${\mathbb R}\!$-ho\-lo\-mor\-p\-hic at a point $(z_0,w_0)\in G,$ then a 
completely solvable in the domain $G$ system of total differential equations {\rm (1.1)}  
has a unique solution 
$w=w(z)$ ${\mathbb R}\!$-holomorphic at the point $z_0$ and satisfying the 
initial condition $w(z_0)=w_0.$
}
\vspace{0.5ex}

{\sl Proof}. 
Taking into account properties (1.5), we construct the system 
conjugate to (1.1):
\\[2ex]
\mbox{}\hfill                               
$
d\,\overline{w} = \overline{X}_2(z,w)\,dz + \overline{X}_1(z,w)\,d\,\overline{z}\,,
$
\hfill (1.7)
\\[2ex]
for which the complete solvability conditions $(\overline{1.6})$ conjugate to (1.6) are 
satisfied in the $G.$  

The functions $X_{\xi l}\colon G\to {\mathbb C},$ which 
are ${\mathbb R}\!$-holomorphic in the domain $G,$ can be treated as functions
$h_{\xi l}\colon \Omega\to {\mathbb C}$ holomorphic in the domain
$\Omega\subset {\mathbb C}^{\, 2(m+n)}$ and such that 
\\[1.25ex]
\mbox{}\hfill
$h_{\xi l}(z,\overline{z},w,\overline{w}\,)=X_{\xi l}(z,w),
\ \ \ \xi=1,\ldots, n,\ \ l=1,\ldots, 2m.
\hfill
$ 
\\[1.25ex]
\indent
Using a correspondence similar to (1.4), on the basis of differential system 
$\!(1.1)\!\cup\! (1.7)\!$ un\-der 
con\-di\-ti\-ons $(1.6)\!\cup\! (\overline{1.6}\!\;)$ we construct the system 
\\[1.5ex]
\mbox{}\hfill
$
\displaystyle
dx_{\xi}=\sum\limits_{l=1}^{2m} h_{\xi l}(t,x)\,dt_{l},
\ \ \  
dx_{n+\xi}=\sum\limits_{j=1}^{m} 
\bigl(\overline{h}_{\xi, m+j}(t,x)\,dt_{j}+
\overline{h}_{\xi j}(t,x)\,dt_{m+j}\bigr),
\, \ \xi=1,\ldots, n,
$ 
\hfill (1.8)
\\[1ex]
with the independent variables $\!(t_1,\ldots, t_{2m})\!=\!t$ and 
\vspace{0.15ex}
the dependent variables $\!(x_1,\ldots, x_{2n})\!=\!x.$
This is a completely solvable system, and therefore (e.g., see [14, p. 26]), it 
\vspace{0.2ex}
has a unique solution $x=x(t)$ holomorphic at the point 
\vspace{0.25ex}
$\!t_0\!=\!(t_1^0,\ldots, t_{2m}^0)$ and satisfying the initial con\-di\-tion $x(t_0)=x_0,$
where the point $(t_0,x_0)\in \Omega,\   
(z_0,\overline{z}_0)\mapsto t_0,\ (w_0,\overline{w}_0)\mapsto x_0.$
\vspace{0.25ex}

Since $w=\varphi_1(z)$ and $\overline{w}=\varphi_2(z)$ are solutions of system 
$\!(1.1)\!\cup\! (1.7)\!$ ${\mathbb R}\!$-holomorphic at the point $z_0,$ it follows 
that the functions $w=\overline{\varphi}_2(z)$ and $\overline{w}=\overline{\varphi}_1(z),$
${\mathbb R}\!$-holomorphic at the point $z_0,$ are also solutions.
Therefore, system $\!(1.1)\!\cup\! (1.7)\!$ un\-der 
con\-di\-ti\-ons $(1.6)\!\cup\! (\overline{1.6}\!\;)$
has the unique solution $w=w(z),\ \overline{w}=\overline{w}(z)$ 
${\mathbb R}\!$-holomorphic at the point $z_0$ and 
satisfying the initial con\-di\-tions $w(z_0)=w_0$ and 
$\overline{w}(z_0)=\overline{w}_0.$ 
One can obtain it from the solution $x=x(t)$ of system (1.8) holomorphic at the point $t_0$
and satisfying the initial con\-di\-tion $x(t_0)=x_0$ with the help of the 
correspondence used when deriving system (1.8). 

Since system $\!(1.1)\!\cup\! (1.7)\!$
splits into systems (1.1) and (1.7), it follows that the original system (1.1) equipped with 
condition (1.6) has a unique ${\mathbb R}\!$-holomorphic solution at the point $z_0$
with the initial data $(z_0,w_0)\in G.$ \k
\vspace{0.5ex}

Theorem 1.1 is a counterpart of the well-known Cauchy theorem on a holomorphic solution 
for the case of  ${\mathbb R}\!$-holomorphic solutions; therefore, it is naturally referred to as 
the Cauchy existence and uniqueness theorem for an ${\mathbb R}\!$-holomorphic solution.  

By [15], the completely solvable system (1.8) has no holomorphic solutions 
(except for the holomorphic solution of the Cauchy problem with the initial condition 
$x(t_0)=x_0$ and $(t_0,x_0)\in \Omega)$ that are not holomorphic at the point $t_0$ and 
tend to $x_0$ as 
$t_l\to t_l^0$ along some paths $\gamma_l,\ l=1,\ldots, 2m.$ Just as in the proof of 
Theorem 1.1, hence we obtain the following property of 
${\mathbb R}\!$-holomorphic solution of the system (1.1).
\vspace{0.5ex}

{\bf Theorem 1.2} [15].
{\it
System {\rm (1.1)} completely solvable in the domain $G$ does not have an 
${\mathbb R}\!$-holomorphic solution that is not ${\mathbb R}\!$-holomorphic 
at $z_0$ and tends to $w_0$ as  $z\to z_0$ along some path $\gamma,$ where 
$(z_0,w_0)\in G.$
}
\\[2.5ex]
\indent
{\bf  1.4. ${\mathbb R}\!$-differentiable integrals and last multipliers}
\\[0.75ex]
\indent
For the unambiguous understanding of our notions we follow [16, p. 29; 17, p. 81; 18;  19, 
pp. 161 -- 178] and introduce the definitions. 
\vspace{0.25ex}

An $\!{\mathbb R}\!$-differentiable on a domain $\!G^{\,\prime}$  function:  
i)  $\!F\colon\! G^{\,\prime}\!\to {\mathbb C};$
ii) $\!f\colon G^{\,\prime}\!\to {\mathbb C};$
iii) $\!\mu\colon G^{\,\prime}\!\to {\mathbb C}$ 
is called  i) a {\it first integral}; ii) a {\it partial integral}; iii) a {\it last mul\-ti\-p\-lier} 
of the system of total differential equations (1.1) if and only if   
\\[0.75ex]
\indent
i)  ${\frak X}_{l}F(z,w)=0$ for all $(z,w)\in G^{\,\prime}, \ l=1,\ldots, 2m,\ \ 
G^{\,\prime}\subset G;$
\\[1ex]
\indent
ii)  ${\frak X}_{l}f(z,w)=\Phi_{l}(f; z,w)$ for all $ (z,w)\in G^{\,\prime},$ 
where $\Phi_{l}(0; z,w)\equiv 0, \ l=1,\ldots, 2m;$
\\[1ex]
\indent
iii) ${\frak X}_{l}\mu(z,w)={}-\mu(z,w)\,{\rm div}\,{\frak X}_{l}(z,w)$ for all 
$(z,w)\in G^{\,\prime},\ \ l=1,\ldots, 2m,$
\\[1.25ex]      
where the linear differential operators
\\[1.5ex]
\mbox{}\hfill                                  
$
\displaystyle
{\frak X}_{j}(z,w) =  
\partial_{z_{j}} +
\sum\limits_{\xi=1}^{n}\bigl(
X_{\xi j}(z,w)\partial_{w_{\xi}} +
\overline{X}_{\xi,m+j}(z,w)
\partial_{{}_{\scriptstyle \overline{w}_{\xi}}}\bigr)$
for all $(z,w)\in G,
\ \ \ j=1,\ldots, m,
\hfill
$
\\[1.5ex]
\mbox{}\hfill                                  
$
\displaystyle
{\frak X}_{m+j}(z,w)=
\partial_{{}_{\scriptstyle \overline{z}_{j}}} +
\sum\limits_{\xi=1}^{n}\bigl(
X_{\xi,m+j}(z,w)\partial_{w_{\xi}} +
\overline{X}_{\xi j}(z,w)
\partial_{{}_{\scriptstyle \overline{w}_{\xi}}}\bigr)$
for all $(z,w)\in G,
\ j=1,\ldots, m.
\hfill
$
\\[1.5ex]
\indent
The $\!{\mathbb R}\!$-differentiable first integral $\!F$  
(partial integral $f$  and last multiplier $\!\mu)$ of the sys\-tem of total differential equations (1.1) 
is called {\it $(s_1,s_2)\!$-non\-autonomous} [20; 21] if
\\[0.5ex]
\indent
(i)  $F \ (f$ and $\mu)$ is holomorphic of $m-s_2$ independent variables;
\\[0.75ex]
\indent
(ii)  $F \ (f$ and $\mu)$ is antiholomorphic of  $m-s_1$ independent variables.
\\[1ex]
\indent
The ${\mathbb R}\!$-differentiable first integral $F$  (partial integral $f$  and last multiplier $\mu)$ 
of the total differential system  (1.1) is called 
{\it $(n-k_1,n-k_2)\!$-cylindricality} [10; 20; 21] if
\\[0.5ex]
\indent
(i)   $F \ (f$ and $\mu)$ is holomorphic of  $n-k_2$ dependent variables;
\\[0.75ex]
\indent
(ii)  $F \ (f$ and $\mu)$ is antiholomorphic of  $n-k_1$ dependent variables.
\\[2.25ex]
\indent
{\bf 1.4.1. $\!\!{\mathbb R}\!$-differentiable partial integrals}.
\vspace{0.5ex}
Suppose the total differential system (1.1) has an ${\mathbb R}\!$-dif\-ferentiable 
$(s_1,s_{2})\!$-non\-autonomous  $(n-k_1,n-k_2)\!$-cy\-lin\-dricality  partial integral 
\\[2ex]
\mbox{}\hfill                              
$
f\colon (z,w)\to  f({}^{s}\!z,{}^{k}\!w)$
\ \ for all $(z,w)\in G^{\,\prime},
$
\hfill (1.9)
\\[2ex]
where $s=(s_1,s_2),\, k=(n-k_1,n-k_2).$ 
\vspace{0.5ex}
We can assume without loss of generality that $f$ is an antiholomorphic function of
\vspace{0.75ex}
$z_{s_1+1},\ldots,z_{m},\ w_{k_1+1},\ldots,w_n$ and $f$ is a holomorphic function of 
$z_{{}_{\scriptstyle j_{s_2+1}}}\!,\ldots, z_{{}_{\scriptstyle j_{m}}}\!, 
\vspace{1.25ex}
w_{{}_{\scriptstyle \zeta_{k_2+1}}}\!,\ldots, w_{{}_{\scriptstyle \zeta_{n}}}
(j_{\beta}\!\in\! \{\!1,\ldots, m\!\}, \beta\!=\!s_2\!+\!1,\ldots, m, 
\zeta_{\delta}\!\in\! \{\!1,\ldots, n\!\}, \delta\!=\!k_2\!+\!1,\ldots, n).\!$

Then, in accordance with the definition of a partial integral, 
\\[2ex]
\mbox{}\hfill                            
$
{\frak X}_{lsk} f({}^{s}\!z,{}^{k}\!w) = \Phi_{l}(f; z,w)$
\ for all $(z,w)\in G^{\,\prime},
\ \ \ l=1,\ldots, 2m,
$
\hfill (1.10)
\\[2ex]
where $\!\Phi_{l}(0; z,w)=0$ for all $(z,w)\in G^{\,\prime}, \ l=1,\ldots, 2m;$
the linear differential operators
\\[2ex]
\mbox{}\hfill
$
\displaystyle
{\frak X}_{\theta sk}(z,w)=
\partial_{z_{\theta}} +
\sum\limits_{\xi=1}^{k_1}X_{\xi\theta}(z,w)\partial_{w_{\xi}} +
\sum\limits_{\tau=1}^{k_2}\overline{X}_{\zeta_{\tau},m+\theta}(z,w)
\partial_{{}_{\scriptstyle \overline{w}_{\zeta_{\tau}}}}$
for all $(z,w)\in G,
\hfill
$
\\[1.75ex]
\mbox{}\hfill
$
\displaystyle
{\frak X}_{\eta sk}(z,w)=
\sum\limits_{\xi=1}^{k_1}X_{\xi\eta}(z,w)\partial_{w_{\xi}} +
\sum\limits_{\tau=1}^{k_2}\overline{X}_{\zeta_{\tau},m+\eta}(z,w)
\partial_{{}_{\scriptstyle \overline{w}_{\zeta_{\tau}}}}$
for all $(z,w)\in G,
\hfill
$
\\[1.75ex]
\mbox{}\hfill
$
\displaystyle
{\frak X}_{m+j_{g},sk}(z,w)=
\partial_{{}_{\scriptstyle \overline{z}_{j_{g}}}} +
\sum\limits_{\xi=1}^{k_1}X_{\xi,m+j_{g}}(z,w)\partial_{w_{\xi}} +
\sum\limits_{\tau=1}^{k_2}\overline{X}_{\zeta_{\tau}j_{g}}(z,w)
\partial_{{}_{\scriptstyle \overline{w}_{\zeta_{\tau}}}}$
for all $(z,w)\in G,
\hfill
$
\\[1.75ex]
\mbox{}\hfill
$
\displaystyle
{\frak X}_{m+j_{\nu},sk}(z,w)=
\sum\limits_{\xi=1}^{k_1}X_{\xi,m+j_{\nu}}(z,w)\partial_{w_{\xi}} +
\sum\limits_{\tau=1}^{k_2}\overline{X}_{\zeta_{\tau}j_{\nu}}(z,w)
\partial_{{}_{\scriptstyle \overline{w}_{\zeta_{\tau}}}}$
for all $(z,w)\in G,
\hfill
$
\\[1.5ex]
\mbox{}\hfill
$
\theta=1,\ldots, s_1, 
\ \ \
\eta=s_1+1,\ldots, m, \ \ \ 
g=1,\ldots, s_2,   \ \ \
\nu=s_2+1,\ldots, m,
\hfill
$
\\[2ex]
with $j_{g}\in\{1,\ldots,m\},\ j_{\nu}\in\{1,\ldots,m\},\ \zeta_{\tau}\in\{1,\ldots,n\}$
\vspace{0.75ex}
(if $J_{g}=\{j_{g}\colon g=1,\ldots, s_2\}$ and 
$J_{\nu}=\{j_{\nu}\colon \nu=s_2+1,\ldots, m\},$ then
$J_{g}\cap J_{\nu}=\O$ and ${\rm Card}\,J_{g}\cup J_{\nu}=m).$
\vspace{1.25ex}

System (1.10) implies that the functions from the sets
\\[1.75ex]
\mbox{}\hfill                            
$
\bigl\{
1,X_{1\theta}(z,w),\ldots,X_{k_1\theta}(z,w),
\overline{X}_{\zeta_{1},m+\theta}(z,w),\ldots,
\overline{X}_{\zeta_{{}_{\scriptstyle k_2}}{,}m+\theta}(z,w)\bigr\},
\ \ \theta=1,\ldots, s_1,
\hfill
$
\\[2ex]
\mbox{}\quad
$
\bigl\{
X_{1\eta}(z,w),\ldots,X_{k_1\eta}(z,w),
\overline{X}_{\zeta_{1},m+\eta}(z,w),\ldots,
\overline{X}_{\zeta_{{}_{\scriptstyle k_2}}{,}m+\eta}(z,w)\bigr\},
\ \ \eta=s_1+1,\ldots, m,
\hfill
$
\\
\mbox{}\hfill (1.11)
\\
\mbox{}\quad
$
\bigl\{
1,X_{1,m+j_{g}}(z,w),\ldots,X_{k_1,m+j_{g}}(z,w),
\overline{X}_{\zeta_{1}j_{g}}(z,w),\ldots,
\overline{X}_{\zeta_{{}_{\scriptstyle k_2}}\!,\,j_{g}}(z,w)\bigr\},
\ \ g=1,\ldots,s_2,
\hfill
$
\\[2ex]
\mbox{}\hfill
$
\bigl\{
X_{1,m+j_{\nu}}(z,w),\ldots,X_{k_1,m+j_{\nu}}(z,w),
\overline{X}_{\zeta_{1}j_{\nu}}(z,w),\ldots,
\overline{X}_{\zeta_{{}_{\scriptstyle k_2}}\!,\,j_{\nu}}(z,w)\bigr\},
\ \ \nu=s_2+1,\ldots, m,
\hfill
$
\\[2ex]
are linearly bound\footnote[3]{
Note that functions  (operators) are called 
{\it linearly bound} on the domain $G$ if 
these functions (operators) are linearly dependent in any point of the domain $G.$
} [22, p. 90; 23, pp. 113 -- 114] on the integral manifold 
\\[1.75ex]
\mbox{}\hfill                            
$
f({}^{s}\!z,{}^{k}\!w)=0.
$
\hfill (1.12)
\\[1.75ex]
\indent
Therefore the Wronskians of the functions from the sets (1.11) with respect to  
$\!z_{\alpha}, \ \overline{z}_{j_{\beta}},$ 
and with respect to $w_{\gamma},\ \overline{w}_{\zeta_{\delta}}\ 
\vspace{0.35ex}
(\alpha=s_1+1,\ldots,m,\ \beta=s_2+1,\ldots, m,\ \gamma=k_1+1,\ldots, n,$ 
$\delta=k_2+1,\ldots, n)$
vanish identically on the manifold (1.12), i.e., the system of identities
\\[2.25ex]
\mbox{}\hfill                                       
$
W_{{}_{\scriptstyle \chi}}
\bigl(1,{}^{\lambda}\!X^{\theta}(z,w)\bigr)= \Psi_{\theta\chi}(f;z,w)$
\ for all $(z,w)\in G,
\ \ \ \theta=1,\ldots, s_1,
\hfill
$
\\[2.5ex]
\mbox{}\hfill
$
W_{{}_{\scriptstyle \chi}}
\bigl({}^{\lambda}\!X^{\eta}(z,w)\bigr)=\Psi_{\eta\chi}(f;z,w)$
\ for all $(z,w)\in G,
\ \ \ \eta=s_1+1,\ldots,m,
$
\hfill (1.13)
\\[2.5ex]
\mbox{}\hfill
$
W_{{}_{\scriptstyle \chi}}
\bigl(1,{}^{\lambda}\!X^{m+j_{g}}(z,w)\bigr)=\Psi_{m+j_{g},\chi}(f;z,w)$
\ for all $(z,w)\in G,
\ \ \ g=1,\ldots, s_2,
\hfill
$
\\[2.5ex]
\mbox{}\hfill
$
W_{{}_{\scriptstyle \chi}}
\bigl({}^{\lambda}\!X^{m+j_{\nu}}(z,w)\bigr)=\Psi_{m+j_{\nu},\chi}(f;z,w)$
\ for all $(z,w)\in G,
\ \ \ \nu=s_2+1,\ldots,m,
\hfill
$
\\[2.25ex]
is consistent. 
\vspace{0.35ex}
Here $W_{{}_{\scriptstyle \chi}}$ are the Wronskians with respect to $\chi$  
(the variable $\chi$  ranges over  $z_{\alpha},$  $\alpha\!=\!s_1\!+\!1,\ldots,m,\ 
\vspace{0.35ex}
\overline{z}_{j_{\beta}},\, \beta=s_2+1,\ldots, m,\ w_{\gamma},\, \gamma=k_1+1,\ldots, n,\ 
\overline{w}_{\zeta_{\delta}},\, \delta\!=\!k_2+1,\ldots, n);$ the number $\lambda=k_1+k_2;$ 
the vector functions
\\[2ex]
\mbox{}\hfill
$
{}^{\lambda}\!X^{j}\colon (z,w)\to
\bigl(X_{1j}(z,w),\ldots, X_{k_1 j}(z,w),
\overline{X}_{\zeta_1,m+j}(z,w),\ldots,
\overline{X}_{{}_{\scriptstyle \zeta_{k_2},\,m+j}}(z,w)\bigr),
\hfill
$
\\[2ex]
\mbox{}\hfill
$
{}^{\lambda}\!X^{m+j}\colon (z,w)\to
\bigl(X_{1,m+j}(z,w),\ldots,X_{k_1,m+j}(z,w),
\overline{X}_{\zeta_1j}(z,w),\ldots,
\overline{X}_{{}_{\scriptstyle \zeta_{k_2}\,j}}(z,w)\bigr)
\hfill
$
\\[1.75ex]
\mbox{}\hfill
for all $(z,w)\in G, 
\quad
j=1,\ldots,m;
\hfill
$
\\[1.75ex]
$\Psi_{l\chi}\colon \!G\to {\mathbb C}$ are $\!{\mathbb R}\!$-differentiable functions 
\vspace{0.35ex}
of  $z$ and $w$ on the domain $\!G$ and 
$\!\Psi_{l\chi}\!(0;z,w)\!\equiv 0,$ $l=1,\ldots, 2m.$ 
Thus, the following theorem is valid.
\vspace{0.75ex}

{\bf Theorem 1.3.}\!
{\it 
For the system of total differential equations {\rm (1.1)} 
to have a partial in\-teg\-ral of the form {\rm (1.9)} it is necessary that 
{\rm (1.13)} be consistent.
}
\vspace{0.65ex}

{\bf Corollary 1.2.}\!
{\it 
For the total differential system {\rm (1.1)} to have a 
$(s_1,0)\!$-non\-autonomous  $(n-k_1,n)\!$-cy\-lin\-dricality holomorphic
partial in\-teg\-ral of the form 
{\rm (1.9)} it is necessary that  the system of identities
{\rm (1.13)} with $s_2=0,\, k_2=0$ be consistent.
}
\vspace{0.5ex}

{\bf Corollary 1.3.}\!
{\it 
For the total differential system {\rm (1.1)} to have a 
$(0,s_2)\!$-non\-au\-to\-no\-mous  $(n,n-k_2)\!$-cy\-lin\-dricality antiholomorphic
partial in\-teg\-ral of the form 
{\rm (1.9)} it is necessary that the system of identities
{\rm (1.13)} with $s_1=0,\, k_1=0$ be consistent.
}
\vspace{0.5ex}

{\bf Corollary 1.4.}\!
{\it 
For the system {\rm (1.1)} to have an 
autonomous  $(n-k_1,n-k_2)\!$-cy\-lin\-dricality 
$\!{\mathbb R}\!$-differentiable partial in\-teg\-ral $f\colon w\to  f({}^{k}\!w)$
for all $w\in \Omega^{\,\prime},\ \Omega^{\,\prime}\subset {\mathbb C}^{n},$
it is necessary that the system of identities {\rm (1.13)} with $s_1=0,\, s_2=0$ be consistent.
}

\newpage

Let the system (1.1) satisfy conditions (1.13). Let us write out the system of equations
\\[2ex]
\mbox{}\hfill                                       
$
\psi_{\theta s_1} +
{}^{\lambda}\!\varphi
\bigl[{}^{\lambda}\!X^{\theta}(z,w)\bigr]^{T}=
H_{\theta}(f;z,w),
\quad
{}^{\lambda}\!\varphi
\bigl[\partial_{{}_{\scriptstyle \chi}}^{p}\,
{}^{\lambda}\!X^{\theta}(z,w)\bigr]^{T}=
\partial_{{}_{\scriptstyle \chi}}^{p} H_{\theta}(f;z,w),
\ \ p=1,\ldots,\lambda,
\hfill
$
\\[2.75ex]
\mbox{}\hfill
$
{}^{\lambda}\!\varphi
\bigl[{}^{\lambda}\!X^{\eta}(z,w)\bigr]^{T}= H_{\eta}(f; z,w),
\quad
{}^{\lambda}\!\varphi
\bigl[\partial_{{}_{\scriptstyle \chi}}^{p}\,
{}^{\lambda}\!X^{\eta}(z,w)\bigr]^{T}=
\partial_{{}_{\scriptstyle \chi}}^{p}
H_{\eta}(f;z,w),
\ \ p=1,\ldots,\lambda-1,
\hfill
$
\\[2.75ex]
\mbox{}\hfill
$
\psi_{g s_2} +
{}^{\lambda}\!\varphi
\bigl[{}^{\lambda}\!X^{m+j_{g}}(z,w)\bigr]^{T}=
H_{m+j_{g}}(f; z,w),
$
\hfill (1.14)
\\[2.75ex]
\mbox{}\hfill
$
{}^{\lambda}\!\varphi
\bigl[\partial_{{}_{\scriptstyle \chi}}^{p}\,
{}^{\lambda}\!X^{m+j_{g}}(z,w)\bigr]^{T}=
\partial_{{}_{\scriptstyle \chi}}^{p}
H_{m+j_{g}}(f; z,w),
\ \ p=1,\ldots,\lambda,
\hfill 
$
\\[2.75ex]
\mbox{}\hfill
$\!\!
{}^{\lambda}\!\varphi
\bigl[{}^{\lambda}\!X^{m+j_{\nu}}(z,w)\bigr]^{T}\!\!=\!
H_{m+j_{\nu}}(f; z,w),
\ \,
{}^{\lambda}\!\varphi
\bigl[\partial_{{}_{\scriptstyle \chi}}^{p}
{}^{\lambda}\!X^{m+j_{\nu}}(z,w)\bigr]^{T}\!\!=\!
\partial_{{}_{\scriptstyle \chi}}^{p}
H_{m+j_{\nu}}(f; z,w),
p\!=\!1,\!\ldots\!,\lambda\!-\!1,\!\!\!
\hfill
$
\\[2.5ex]
\mbox{}\hfill
$
\theta=1,\ldots, s_1,
\ \ \ 
\eta=s_1+1,\ldots,m,
\ \ \
g=1,\ldots, s_2,
\ \ \
\nu=s_2+1,\ldots m,
\hfill
$
\\[2ex]
where the vector functions 
\\[2ex]
\mbox{}\hfill
$
{}^{k_1}\!\varphi\colon (z,w)\to
(\varphi_{1 k_1}({}^{s}\!z,{}^{k}\!w),\ldots,
\varphi_{k_1 k_1}({}^{s}\!z,{}^{k}\!w)), 
\ \
{}^{k_2}\!\varphi\colon (z,w)\to
(\varphi_{1 k_2}({}^{s}\!z,{}^{k}\!w),\ldots,
\varphi_{k_2 k_2}({}^{s}\!z,{}^{k}\!w)),
\hfill
$
\\[2.75ex]
\mbox{}\hfill
$
{}^{\lambda}\!\varphi\colon (z,w)\to
({}^{k_1}\!\varphi (z,w),{}^{k_2}\!\varphi (z,w))$
 \ for all $(z,w)\in G;
\hfill
$
\\[2ex]
$H_l\colon\! G\to\C$ are 
\vspace{0.35ex}
$\!{\mathbb R}\!$-differentiable functions of  $z$ and $w$ on 
the domain $\!G$ and
$\!H_{l}(0;z,w)\!\equiv\! 0,$ $l=1,\ldots,2m.$
Let us introduce the Pfaffian differential equation
\\[2.5ex]
\mbox{}\hfill                              
$
{}^{s_1}\!\psi({}^{s}\!z,{}^{k}\!w) d{}^{s_{\!1}}\!z +
{}^{s_2}\!\psi({}^{s}\!z,{}^{k}\!w) d\,\overline{{}^{s_{\!2}}\!z}+
{}^{k_1}\!\varphi({}^{s}\!z,{}^{k}\!w) d{}^{k_{\!1}}\!w+
{}^{k_2}\!\varphi({}^{s}\!z,{}^{k}\!w)
d\,\overline{{}^{k_{2}}\!w}=0,
$
\hfill (1.15)
\\[2.75ex]
where $d{}^{s_1}\!z\!=\!\mbox{colon}(dz_1,\ldots,dz_{s_1}),\,
d\,\overline{{}^{s_{2}}\!z}\!=\!\mbox{colon}\bigl(\!
d\,\overline{z}_{{}_{\scriptstyle j_1}},\ldots,
d\,\overline{z}_{{}_{\scriptstyle j_{s_2}}}\!\bigr),\,
d{}^{k_1}\!w\!=\!\mbox{colon}(dw_1,\ldots,dw_{k_1}\!),\!$
and $d\,\overline{{}^{k_{2}}\!w}=\mbox{colon}\bigl(
d\,\overline{w}_{{}_{\scriptstyle \zeta_1}},\ldots,
d\,\overline{w}_{{}_{\scriptstyle \zeta_{k_2}}}\bigr)$
are vector columns; the vector functions
\\[1.5ex]
\mbox{}\hfill
$
{}^{s_1}\!\psi\colon (z,w)\to
(\psi_{1 s_1}({}^{s}\!z,{}^{k}\!w),\ldots,
\psi_{s_1 s_1}({}^{s}\!z,{}^{k}\!w))$
\ for all $(z,w)\in G,
\hfill
$
\\[2.25ex]
\mbox{}\hfill
$
{}^{s_2}\!\psi\colon (z,w)\to
(\psi_{1 s_2}({}^{s}\!z,{}^{k}\!w),\ldots,
\psi_{s_2 s_2}({}^{s}\!z,{}^{k}\!w))$
\ for all $(z,w)\in G.
\hfill
$
\\[2.5ex]
\indent
{\bf Theorem 1.4.}\! 
{\it
A necessary and sufficient condition for the total differential system {\rm (1.1)} 
to have at least one partial integral of the form {\rm (1.9)} 
is that the functions ${}^{s_1}\!\psi,\ {}^{s_2}\!\psi,\ {}^{\lambda}\!\varphi,$  and     
$H_l,\ l=1,\ldots, 2m,$ exist so that they satisfy system {\rm (1.14)} and
\vspace{0.25ex}

{\rm (i)}\ the Pfaffian differential equation {\rm (1.15)} has an integrating factor{\rm;}
\vspace{0.25ex}

{\rm (ii)}\ the function {\rm (1.9)} is a general integral of the Pfaffian differential equation {\rm (1.15)}.  
}
\vspace{0.75ex}

{\sl Proof. Necessity}.\! 
Let the total differential system (1.1) have a 
$\!{\mathbb R}\!$-differentiable partial in\-teg\-ral of the form (1.9). 
Then the identity (1.10) holds. The vector functions 
\\[1.75ex]
\mbox{}\hfill
$
{}^{s_1}\psi\colon (z,w)\to
\partial_{{}_{\scriptstyle {}^{s_1}\!z}}
f({}^{s}\!z,{}^{k}\!w),
\ \ \
{}^{s_2}\psi\colon (z,w)\to
\partial_{{}_{\scriptstyle \overline{{}^{s_2}\!z}}}\,
f({}^{s}\!z,{}^{k}\!w)$
\ for all $(z,w)\in G^{\,\prime},
\hfill
$
\\[2.75ex]
\mbox{}\hfill
$
{}^{k_1}\varphi\colon (z,w)\to
\partial_{{}_{\scriptstyle {}^{k_1}\!w}}
f({}^{s}\!z,{}^{k}\!w),
\ \ \
{}^{k_2}\varphi\colon (z,w)\to
\partial_{{}_{\scriptstyle \overline{{}^{k_2}\!w}}}\,
f({}^{s}\!z,{}^{k}\!w)$
\  for all $(z,w)\in G^{\,\prime},
\hfill
$
\\[2.25ex]
where 
\vspace{0.35ex}
$\partial_{{}_{\scriptstyle {}^{s_{\!1}}\!z}}\!\!=\!
(\!\partial_{z_1},\ldots,\partial_{z_{s_1}}\!),
\partial_{{}_{\scriptstyle \overline{{}^{s_{\!2}}z}}}\! =\!
\bigl(\partial_{{}_{\scriptstyle \overline{z}_{j_1}}}\!,\!\ldots\!,
\partial_{{}_{\scriptstyle \overline{z}_{j_{s_2}}}}\!\bigr),
\partial_{{}_{\scriptstyle {}^{k_{\!1}}\!w}}\!\!=\!
(\partial_{w_1},\!\ldots\!,\partial_{w_{k_1}}\!), 
\partial_{{}_{\scriptstyle \overline{{}^{k_{\!2}}w}}} \!=\!
\bigl(\partial_{{}_{\scriptstyle \overline{w}_{\zeta_1}}}\!,\!\ldots\!,
\partial_{{}_{\scriptstyle \overline{w}_{\zeta_{k_2}}}}\!\!\bigr),$
is a solution to system (1.14) for  $H_{l}(f; z,w)=\Phi_{l}(f; z,w),\ l=1,\ldots, 2m,$
\vspace{0.35ex}
which can be shown by dif\-fe\-ren\-ti\-a\-ting (1.10)  
$\lambda$ times with respect to $\chi\ (\theta=1,\ldots, s_1, \ g=1,\ldots, s_2)$
\vspace{0.35ex}
and $\lambda-1$ times with respect to  $\chi\ (\eta=s_1+1,\ldots, m, \ \nu=s_2+1,\ldots, m).$
\vspace{0.35ex}
Therefore the ${\mathbb R}\!$-differentiable function\! (1.9) is a 
\vspace{0.75ex}
general integral of the Pfaffian differential equation\! (1.15). 

{\sl Sufficiency}. 
\vspace{0.35ex}
Let ${}^{s_1}\!\psi,\ {}^{s_2}\!\psi,\ {}^{\lambda}\!\varphi$ be a solution to 
the system (1.14), and let the corresponding Pfaffian differential equation 
(1.15) have an integrating factor $\mu\colon ({}^{s}\!z,{}^{k}\!w)\to \mu({}^{s}\!z,{}^{k}\!w)$
and the corresponding general integral (1.9).  Then 

\newpage

\mbox{}
\\[-1.75ex]
\mbox{}\quad\                        
$
\partial_{{}_{\scriptstyle {}^{s_1}\!z}}\,f({}^{s}\!z,{}^{k}\!w)-
\mu ({}^{s}\!z,{}^{k}\!w)\,{}^{s_1}\psi({}^{s}\!z,{}^{k}\!w)=0,
\quad
\partial_{{}_{\scriptstyle \overline{{}^{s_2}\!z}}}\,
f({}^{s}\!z,{}^{k}\!w) -
\mu ({}^{s}\!z,{}^{k}\!w)\,{}^{s_2}\psi({}^{s}\!z,{}^{k}\!w)=0,
\hfill
$
\\
\mbox{}\hfill (1.16)
\\
\mbox{}\quad\
$
\partial_{{}_{\scriptstyle {}^{k_1}\!w}}\,f({}^{s}\!z,{}^{k}\!w)-
\mu ({}^{s}\!z,{}^{k}\!w)\,{}^{k_1}\varphi({}^{s}\!z,{}^{k}\!w)=0,
\quad
\partial_{{}_{\scriptstyle \overline{{}^{k_2}\!w}}}\,
f({}^{s}\!z,{}^{k}\!w) -
\mu ({}^{s}\!z,{}^{k}\!w)\,{}^{k_2}\varphi({}^{s}\!z,{}^{k}\!w)=0.
\hfill
$
\\[2ex]
\indent
It follows from (1.14) and (1.16) that identity (1.10) is valid with
\\[1.75ex]
\mbox{}\hfill
$
\Phi_{l}(f; z,w)=\mu({}^{s}\!z,{}^{k}\!w)H_{l}(f; z,w)$
\ for all $(z,w)\in G^{\,\prime},
\ \ \  l=1,\ldots, 2m.
\hfill
$
\\[1.75ex]
\indent
Consequently, the function (1.9) is a partial integral of the system (1.1).\, \k
\vspace{1.25ex}

{\bf Theorem 1.5.} 
{\it
Let $h$ systems {\rm (1.14)} have $q$ not linearly bound solutions
\\[1.75ex]
\mbox{}\hfill                       
$
{}^{s_1}\psi^{\varepsilon}\colon (z,w)\to
{}^{s_1}\psi^{\varepsilon}({}^{s}\!z,{}^{k}\!w),
\quad
{}^{s_2}\psi^{\varepsilon}\colon (z,w)\to
{}^{s_2}\psi^{\varepsilon}({}^{s}\!z,{}^{k}\!w),
\hfill
$
\\[-0.15ex]
\mbox{}\hfill {\rm (1.17)}
\\[-0.15ex]
\mbox{}\hfill
$
{}^{\lambda}\varphi^{\varepsilon}\colon (z,w)\to
{}^{\lambda}\varphi^{\varepsilon}({}^{s}\!z,{}^{k}\!w)$
\ for all $(z,w)\in G^{\,\prime},
\quad \varepsilon=1,\ldots,q,
\hfill
$
\\[1.75ex]
for which the corresponding Pfaffian differential equations
\\[2.5ex]
\mbox{}\hfill                              
$\!\!
{}^{s_1}\!\psi^{\varepsilon}({}^{s}\!z,{}^{k}\!w)\,d{}^{s_{\!1}}\!z
+\, {}^{s_2}\!\psi^{\varepsilon}({}^{s}\!z,{}^{k}\!w)\,
d\,\overline{{}^{s_{\!2}}\!z}\, +\,
{}^{k_1}\!\varphi^{\varepsilon}({}^{s}\!z,{}^{k}\!w)\,
d{}^{k_{\!1}}\!w + \,
{}^{k_2}\!\varphi^{\varepsilon}({}^{s}\!z,{}^{k}\!w)\,
d\,\overline{{}^{k_{2}}\!w}=0,
\varepsilon\!=\!1,\ldots, q
\ {\rm(1.18)} \!\!
$
\\[2ex]
have the general $\!{\mathbb R}\!$-differentiable integrals     
\\[1.25ex]
\mbox{}\hfill
$
f_{\varepsilon}\colon (z,w)\to f_{\varepsilon}({}^{s}\!z,{}^{k}\!w)$
\ for all $(z,w)\in G^{\,\prime},
\quad\varepsilon=1,\ldots, q.
\hfill
$
\\[1.25ex]
Then these integrals are functionally independent. 
}
\vspace{0.75ex} 

{\sl Proof}. We have 
\\[1.25ex]
\mbox{}\hfill
$
\partial_{{}_{\scriptstyle {}^{s_1}\!z}}\,
f_{\varepsilon}({}^{s}\!z,{}^{k}\!w)=
\mu_{\varepsilon}({}^{s}\!z,{}^{k}\!w)\,
{}^{s_1}\psi^{\varepsilon}({}^{s}\!z,{}^{k}\!w),
\quad
\partial_{{}_{\scriptstyle \overline{{}^{s_2}\!z}}}\,
f_{\varepsilon}({}^{s}\!z,{}^{k}\!w) =
\mu_{\varepsilon}({}^{s}\!z,{}^{k}\!w)\,
{}^{s_2}\psi^{\varepsilon}({}^{s}\!z,{}^{k}\!w),
\hfill
$
\\[2ex]
\mbox{}\hfill
$
\partial_{{}_{\scriptstyle {}^{k_1}\!w}}\,
f_{\varepsilon}({}^{s}\!z,{}^{k}\!w) =
\mu_{\varepsilon}({}^{s}\!z,{}^{k}\!w)\,
{}^{k_1}\varphi^{\varepsilon}({}^{s}\!z,{}^{k}\!w),
\quad
\partial_{{}_{\scriptstyle \overline{{}^{k_2}\!w}}}\,
f_{\varepsilon}({}^{s}\!z,{}^{k}\!w)=
\mu_{\varepsilon}({}^{s}\!z,{}^{k}\!w)\,
{}^{k_2}\varphi^{\varepsilon}({}^{s}\!z,{}^{k}\!w)
\hfill
$
\\[2ex]
\mbox{}\hfill
for all 
$
(z,w)\in G^{\prime},
\quad
\varepsilon=1,\ldots, q,
\hfill
$
\\[1.75ex]
by virtue of (1.16).  Therefore, the Jacobi matrix
\\[2.5ex]
\mbox{}\hfill
$
J(f_{\varepsilon}({}^{s}\!z,\!{}^{k}\!w);{}^{s}\!z,{}^{k}\!w)\!=\!
\bigl\|
{}^{s_1}\!\Psi({}^{s}\!z,{}^{k}\!w)\,
{}^{s_2}\Psi({}^{s}\!z,{}^{k}\!w)\,
{}^{k_1}\!\Phi({}^{s}\!z,{}^{k}\!w)\,
{}^{k_2}\Phi({}^{s}\!z,{}^{k}\!w)
\bigr\|,
\hfill
$
\\[2.75ex]
where
\vspace{1.5ex}
${}^{s_1}\!\Psi\!\! =\!\!
\bigl\|\mu_{\varepsilon}\psi_{\varepsilon \theta s_1}\!\bigr\|\!$
 is a $(q\!\times\! s_1)\!$-matrix,
${}^{s_2}\Psi\!=\!\!
\bigl\|\mu_{\varepsilon}\psi_{\varepsilon g s_2}\!\bigr\|\!$ is a $(q\!\times\! s_2)\!$-matrix,
${}^{k_1}\!\Phi\! =\!\!
\bigl\|\mu_{\varepsilon}\varphi_{\varepsilon \xi k_1}\!\bigr\|\!$
is a $(q\times k_1)\!$-matrix,
\vspace{1ex}
and ${}^{k_2}\Phi=
\bigl\|\mu_{\varepsilon}\varphi_{\varepsilon \tau k_2}\bigr\|$ is a 
$(q\times k_2)\!$-matrix. 

We have  ${\rm rank}\,J=q$ since the solutions (1.17) are not linearly bound.
\vspace{0.5ex}

Consequently, the general ${\mathbb R}\!$-differentiable integrals of the Pfaffian 
equations (1.18) are functionally independent. The proof of the theorem is complete. \k
\vspace{0.5ex}

The Theorem 1.5 (taking into account the Theorem 1.4) let us to find a quantity of 
functionally independent 
$(s_1,s_{2})\!$-non\-autonomous  $(n-k_1,n-k_2)\!$-cy\-lin\-dricality  
${\mathbb R}\!$-dif\-ferentiable partial integrals of 
the total differential system (1.1). 
 \vspace{0.75ex}

For example, the system of total differential equations 
\\[1.5ex]
\mbox{}\hfill                      
$
\begin{array}{l}
dw_1=(w_1^2+w_2\,\overline{w}_2)\,dz +
(w_1w_2+w_2\,\overline{w}_2+
(2+\overline{z}\,)\,\overline{w}\,{}_2^2\,)\,d\,\overline{z}\,,
\\[2.25ex]
dw_2=(w_2\,\overline{w}_1-(1+z)w_2^2\,)\,dz +
\overline{w}_1(w_2+\overline{w}_2)\,d\,\overline{z}
\end{array}
$
\hfill (1.19)
\\[2ex]
has the vector functions (see (1.11))
\\[1.75ex]
\mbox{}\hfill
$
P_1\colon (z,w)\to 
\bigl(w_1w_2+w_2\,\overline{w}_2+
(2+\overline{z}\,)\,\overline{w}\,{}_2^2,
w_1\,\overline{w}_2-
(1+\overline{z}\,)\,\overline{w}\,{}_2^2\,\bigr)$
\ for all $(z,w)\in {\mathbb C}^3,
\hfill
$
\\[2.5ex]
\mbox{}\hfill
$
P_2\colon (z,w)\to 
\bigl((w_1^2+w_2\,\overline{w}_2, w_1(w_2+\overline{w}_2)\bigr)$
\ for all $(z,w)\in {\mathbb C}^3,
\hfill
$
\\[1.75ex]
and the Wronskians (see (1.13))
\\[2ex]
\mbox{}\hfill
$
W_{{}_{\scriptstyle z}}
\bigl(P_1(z,w)\bigr)= 0,
\quad
W_{{}_{\scriptstyle \overline{z}}}\bigl(P_1(z,w)\bigr)=
{}-\overline{w}\,{}_2^2\,(w_2+\overline{w}_2)(w_1+\overline{w}_2),
\hfill
$
\\[2ex]
\mbox{}\hfill
$
W_{{}_{\scriptstyle w_2}}\bigl(P_1(z,w)\bigr)=
{}-(w_1+\overline{w}_2)
(w_1\,\overline{w}_2-(1+\overline{z}\,)\,\overline{w}\,{}_2^2\,),
\quad
W_{{}_{\scriptstyle \overline{w}_1}}\bigl(P_1(z,w)\bigr)=0,
\hfill
$
\\[2ex]
\mbox{}\hfill
$
W_{{}_{\scriptstyle z}}\bigl(P_2(z,w)\bigr)= 0,
\quad
W_{{}_{\scriptstyle \overline{z}}}\bigl(P_2(z,w)\bigr)= 0,
\quad
W_{{}_{\scriptstyle w_2}}\bigl(P_2(z,w)\bigr)=
w_1(w_1-\overline{w}_2)(w_1+\overline{w}_2),
\hfill
$
\\[2ex]
\mbox{}\hfill
$
W_{{}_{\scriptstyle \overline{w}_1}}\bigl(P_2(z,w)\bigr)= 0$
\ \ for all \, $(z,w)\in {\mathbb C}^3
\hfill
$
\\[1.75ex]
vanish identically on the manifold $w_1+\overline{w}_2=0$ (see (1.12)).

Therefore a necessary condition for system of total differential equations (1.19) 
to have an $\!{\mathbb R}\!$-differentiable autonomous 
\vspace{0.35ex}
(1,1)-cylindricality partial integral is complied (Theorem 1.3).

The functions
$
\varphi_1\colon (z,w)\to 1$
for all $(z,w)\in {\mathbb C}^3,
\ \varphi_2\colon (z,w)\to 1$
for all $(z,w)\in {\mathbb C}^3$
is a solution to system (1.14)  for 
\\[1.5ex]
\mbox{}\hfill
$
H_1\colon (z,w)\to (w_1+\overline{w}_2)(w_1+w_2)$
\ for all $(z,w)\in {\mathbb C}^3,
\hfill
$
\\[2ex]
\mbox{}\hfill
$
H_2\colon (z,w)\to (w_1+\overline{w}_2)(w_2+\overline{w}_2)$
\ for all $(z,w)\in {\mathbb C}^3.
\hfill
$
\\[1.5ex]
\indent
The corresponding Pfaffian differential equation
\\[0.75ex]
\mbox{}\hfill
$
dw_1+d\,\overline{w}_2 =0
\hfill
$
\\[1.5ex]
has the integrating factor $\mu\colon w\to 1$ for all $w\in {\mathbb C}^2$ and 
the general integral (Theorem 1.4)
\\[1.5ex]
\mbox{}\hfill                  
$
f\colon (w_1,w_2)\to w_1 +\overline{w}_2$
\ for all $(w_1,w_2)\in {\mathbb C}^2.
$
\hfill (1.20)
\\[1.5ex]
\indent
Thus the system of total differential equations (1.19) has the ${\mathbb R}\!$-differentiable 
au\-to\-no\-mo\-us (1,1)-cylindricality partial integral (1.20).
\\[1.5ex]
\indent
{\bf 1.4.2. ${\mathbb R}\!$-differentiable first integrals}.
Suppose the system of total differential equations (1.1) has a 
$(s_1,s_2)\!$-nonautonomous and $(n-k_1,n-k_2)\!$-cylindricality  
${\mathbb R}\!$-differentiable on the domain $G^{\,\prime}$ first integral 
\\[1.5ex]
\mbox{}\hfill                              
$
F\colon (z,w)\to  F({}^{s}\!z,{}^{k}\!w)$
\ for all $(z,w)\in G^{\,\prime}.
$
\hfill (1.21)
\\[1.5ex]
\indent
Then, in accordance with the criteria of a first integral,
\\[1.75ex]
\mbox{}\hfill                          
$
{\frak X}_{lsk} F({}^{s}\!z,{}^{k}\!w) = 0$
\ for all $(z,w)\in G^{\,\prime},
\quad 
l=1,\ldots, 2m.
\hfill
$
\\[1.75ex]
\indent
Therefore the Wronskians of the functions (1.11) vanish identically on the domain $G,$ i.e., 
the system of identities (1.13) for  $\Psi_{l\chi}\equiv 0, \ l=1,\ldots, 2m$ is consistent in $G.$
\vspace{0.25ex}

We obtain the following statements.
\vspace{0.5ex}

{\bf Theorem 1.6.} 
{\it
For the differential system {\rm (1.1)} to have 
\vspace{0.35ex}
a first integral of the form {\rm (1.21)} it is necessary that {\rm (1.13)} 
with $\Psi_{l\chi}\equiv 0, \ l=1,\ldots, 2m$  be consistent in $G.$
}
\vspace{0.5ex}

{\bf Corollary 1.5.}\!
{\it 
For the total differential system {\rm (1.1)} to have a 
$(s_1,0)\!$-non\-autonomous  $(n-k_1,n)\!$-cy\-lin\-dricality holomorphic
first in\-teg\-ral of the form {\rm (1.21)} it is necessary that  the system of identities
{\rm (1.13)} with $\Psi_{l\chi}\equiv 0, \ l=1,\ldots, 2m,$ and $s_2=0,\, k_2=0$ be consistent.
}
\vspace{0.75ex}

{\bf Corollary 1.6.}\!
{\it 
For the total differential system {\rm (1.1)} to have a 
$(0,s_2)\!$-non\-au\-to\-no\-mous  $(n,n-k_2)\!$-cy\-lin\-dricality antiholomorphic
first in\-teg\-ral of the form 
{\rm (1.21)} it is necessary that the system of identities
{\rm (1.13)} with $\Psi_{l\chi}\equiv 0, \, l=1,\ldots, 2m,$ and $s_1=0,\, k_1=0$ 
\vspace{0.75ex}
be consistent.
}

{\bf Corollary 1.7.}\!
{\it 
For the system {\rm (1.1)} to have an 
autonomous  $(n-k_1,n-k_2)\!$-cy\-lin\-dricality 
$\!{\mathbb R}\!$-differentiable first in\-teg\-ral $F\colon w\to  F({}^{k}\!w)$
for all $w\in \Omega^{\,\prime},\ \Omega^{\,\prime}\subset {\mathbb C}^{n},$
it is necessary that the system of identities {\rm (1.13)} with 
$\Psi_{l\chi}\equiv 0, \, l=1,\ldots, 2m,$ and $s_1=0,\, s_2=0$ 
\vspace{0.75ex}
be consistent.
}

The  proof of the following assertions is similar to those of Theorems 1.4 and 1.5.
\vspace{0.75ex}

{\bf Theorem 1.7.}\!\!
{\it 
For the system of total differential equations {\rm (1.1)} to 
have at least one first integral of the form {\rm (1.21)}  
\vspace{0.25ex}
it is necessary and sufficient that there exist functions  
${}^{s_1}\!\psi,\, {}^{s_2}\!\psi,\, {}^{\lambda}\!\varphi\!$ 
\vspace{0.25ex}
satisfying to system {\rm (1.14)} for  $\!H_{l}\!\equiv\! 0,\, l\!=\!1,\ldots, 2m,\!$ 
that the function {\rm (1.21)} 
is a general integral of the Pfaffian differential equation {\rm (1.15)}.
}
\vspace{1ex}

{\bf Theorem 1.8.} 
{\it 
Let functional system {\rm (1.14)} with  $H_{l}\equiv 0, \ l=1,\ldots, 2m$  
\vspace{0.35ex}
has $q$ not linearly bound so\-lu\-ti\-ons {\rm (1.17)} such that 
the corresponding Pfaffian differential equations {\rm (1.18)} have the general integrals     
\\[1.75ex]
\mbox{}\hfill                              
$
F_{\varepsilon}\colon (z,w)\to
F_{\varepsilon}({}^{s}\!z,{}^{k}\!w)$
\ for all $(z,w)\in G^{\,\prime},
\quad 
\varepsilon=1,\ldots, q.
\hfill 
$
\\[1.75ex]
Then these integrals are functionally independent. 
}
\vspace{0.75ex}

The Theorem 1.8 (taking into account the Theorem 1.7) let us to find a quantity of 
functionally independent 
$(s_1,s_{2})\!$-non\-autonomous  $(n-k_1,n-k_2)\!$-cy\-lin\-dricality  
${\mathbb R}\!$-dif\-ferentiable first integrals of 
the total differential system (1.1). 
 \vspace{1ex}

As an example, the system of total differential equations
\\[2ex]
\mbox{}\hfill                              
$
dw_1 =\dfrac{2}{\overline{z}} \ w_2\,dz -
\Bigl(\, \dfrac{1}{\overline{z}}\,w_1+2w_2^2+
2z\,w_2\,\overline{w}_1\Bigr)d\,\overline{z}\,,
\quad
dw_2={}-dz+ \overline{z}\,(w_2+z\,\overline{w}_1)\,d\,\overline{z}
$
\hfill (1.22)
\\[2.25ex]
has the functions (see (1.11))
\\[2ex]
\mbox{}\hfill
$
P_1\colon (z,w_1,w_2)\to \Bigl(1,{}-\dfrac{1}{z} \ \overline{w}_1-
2\,\overline{w}\,{}_2^2-
2\,\overline{z}\,w_1\,\overline{w}_2,\,
z(\,\overline{z}\,w_1+\overline{w}_2)\Bigr)$
\ \ for all $(z,w_1,w_2)\in\Omega,
\hfill
$
\\[2.5ex]
\mbox{}\hfill
$
P_2\colon (z,w_1,w_2)\to 
\Bigl(\,\dfrac{2\,\overline{w}_2}{z}\,, {}-1\Bigr)$
\ \ for all $(z,w_1,w_2)\in\Omega,
\ \ \ \Omega\subset {\mathbb C}^3.
\hfill
$
\\[2.25ex]
\indent
The Wronskians of the vector functions $P_1$  
\vspace{0.35ex}
and $P_2$  with respect to   $\overline{z},\ w_1,\ w_2$ vanish identically  on a domain
$\Omega\subset \{(z,w_1,w_2)\colon z\ne 0\}\subset {\mathbb C}^3.$
\vspace{0.5ex}

Therefore a necessary condition for the total differential system (1.22) to have an 
$\!{\mathbb R}\!$-dif\-fe\-ren\-ti\-ab\-le (1,0)-nonautonomous 
\vspace{0.35ex}
(2,0)-cylindricality first integral is complied (Theorem 1.6).

The scalar functions 
\\[1.5ex]
\mbox{}\hfill
$
\psi_1\colon (z,w_1,w_2)\to \overline{w}_1,
\ \
\varphi_1\colon (z,w_1,w_2)\to z,
\ \
\varphi_2\colon (z,w_1,w_2)\to 2\,\overline{w}_2$ for all $(z,w_1,w_2)\in \Omega
\hfill
$
\\[2ex]
is a solution to system of equations (see (1.14) with $H_{l}\equiv 0, \ l=1, 2)$
\\[2ex]
\mbox{}\hfill
$
\psi_1 -
\Bigl(\,\dfrac{1}{z} \ \overline{w}_1 +
2\,\overline{w}\,{}_2^2 +
2\,\overline{z}\,w_1\,\overline{w}_2\Bigr)\varphi_1 +
z(\,\overline{w}_2+\overline{z}\,w_1)\varphi_2=0,
\hfill
$
\\[2ex]
\mbox{}\hfill
$
{}-2w_1\,\overline{w}_2\,\varphi_1 + zw_1\,\varphi_2=0,
\ \ \
{}-2\,\overline{z}\,\overline{w}_2\,\varphi_1 +
z\,\overline{z}\,\varphi_2=0,
\ \ \
\dfrac{2}{z} \ \overline{w}_2\,\varphi_1 - \varphi_2=0.
\hfill
$
\\[1.75ex]
\indent
The corresponding Pfaffian differential equation
\\[1.75ex]
\mbox{}\hfill
$
\overline{w}_1\,dz + z\,d\,\overline{w}_1 +
2\,\overline{w}_2\,d\,\overline{w}_2 =0
\hfill
$
\\[1.75ex]
has the general integral (Theorem 1.7)
\\[1.75ex]
\mbox{}\hfill                                       
$
F\colon (z,w_1,w_2)\to \,
z\,\overline{w}_1+\overline{w}\,{}_2^2$
\  \ for all $(z,w_1,w_2)\in \Omega.
$
\hfill (1.23)
\\[1.75ex]
\indent
The Poisson bracket
\\[2.25ex]
\mbox{}\hfill
$
\bigl[{\frak X}_1(z,w), {\frak X}_2(z,w)\bigr] =
\Bigl[\partial_{z} +
\dfrac{2}{\overline{z}} \ w_2\,\partial_{w_1} -
\partial_{w_2} -
\Bigl(\,\dfrac{1}{z} \ \overline{w}_1 +
2\,\overline{w}\,{}_2^2 +
2\,\overline{z}\,w_1\,\overline{w}_2\Bigr)
\partial_{{}_{\scriptstyle \overline{w}_1}} +
z(\,\overline{z}\,w_1+\overline{w}_2)
\partial_{{}_{\scriptstyle \overline{w}_2}}, 
\hfill
$
\\[2.5ex]
\mbox{}\hfill
$
\partial_{{}_{\scriptstyle \overline{z}}} -
\Bigl(\,\dfrac{1}{\overline{z}} \ w_1+2w_2^2+
2z\,w_2\,\overline{w}_1\Bigr)\partial_{w_1} +
\overline{z}\,(w_2+z\,\overline{w}_1)\partial_{w_2} +
\dfrac{2}{z} \ \overline{w}_2\,
\partial_{{}_{\scriptstyle \overline{w}_1}} -
\partial_{{}_{\scriptstyle \overline{w}_2}}\Bigr]=
\hfill
$
\\[2.5ex]
\mbox{}\hfill
$
=\bigl(1+2z\,\overline{w}_2
(\,\overline{z}\,w_1+\overline{w}_2\,)\bigr)
\bigl(2w_2\,\partial_{w_1} -
\overline{z}\,\partial_{w_2}\bigr) -
\bigl(1+2\,\overline{z}\,w_2
(w_2+ z\,\overline{w}_1)\bigr)
\bigl(2\,\overline{w}_2\,
\partial_{{}_{\scriptstyle \overline{w}_1}}-
z\,\partial_{{}_{\scriptstyle \overline{w}_2}}\bigr)
\hfill
$
\\[2.25ex]
is not the null operator on the domain  $\Omega,$ i.e., system (1.22) is not completely solvable.
\vspace{0.35ex}

Thus the ${\mathbb R}\!$-differentiable (1,0)-nonautonomous (2,0)-cylindricality 
first integral (1.23) is an integral basis on the domain  $\Omega$ of the total differential system (1.22).

\newpage

{\bf 1.4.3. $\!{\mathbb R}\!$-differentiable last multipliers}.
Suppose the system of total dif\-fe\-ren\-tial equations (1.1) has a 
$(s_1,s_2)\!$-non\-auto\-no\-mous  and $(n-k_1,n-k_2)\!$-cylindricality 
${\mathbb R}\!$-dif\-fe\-ren\-ti\-ab\-le on the domain $G^{\,\prime}$ last mul\-tip\-lier 
\\[1.75ex]
\mbox{}\hfill                              
$
\mu\colon (z,w)\to  \mu ({}^{s}\!z,{}^{k}\!w)$
\ \ for all $(z,w)\in G^{\,\prime}.
$
\hfill (1.24)
\\[1.75ex]
\indent
Then, in accordance with the criteria of a last multiplier,
\\[1.75ex]
\mbox{}\hfill                            
$
{\frak X}_{lsk} \mu ({}^{s}\!z,{}^{k}\!w)
+ \mu ({}^{s}\!z,{}^{k}\!w)\,{\rm div} {\frak X}_{l}(z,w) = 0$
\ \ for all $(z,w)\in G^{\,\prime},
\quad 
l=1,\ldots, 2m.
$
\hfill (1.25)
\\[2ex]
\indent
Using (1.25), we get 
\\[2ex]
\mbox{}\hfill                                       
$
W_{{}_{\scriptstyle \chi}}
\bigl(1,{}^{\lambda}\!X^{\theta}(z,w),
{\rm div} {\frak X}_{\theta}(z,w)\bigr)=0$
\ for all $(z,w)\in G,$  \ $\theta=1,\ldots, s_1,
\hfill
$
\\[2.25ex]
\mbox{}\hfill
$
W_{{}_{\scriptstyle \chi}}
\bigl({}^{\lambda}\!X^{\eta}(z,w),
{\rm div} {\frak X}_{\eta}(z,w)\bigr)=0$
\ for all $(z,w)\in G,$  \ $\eta=s_1+1,\ldots, m,
\hfill
$
\\[0.2ex]
\mbox{}\hfill (1.26)
\\[0.2ex]
\mbox{}\hfill
$
W_{{}_{\scriptstyle \chi}}
\bigl(1,{}^{\lambda}\!X^{m+j_{g}}(z,w),
{\rm div} {\frak X}_{m+j_{g}}(z,w)\bigr)=0$
\ for all $(z,w)\in G,$  \ $g=1,\ldots, s_2,
\hfill
$
\\[2.25ex]
\mbox{}\hfill
$
W_{{}_{\scriptstyle \chi}}
\bigl({}^{\lambda}\!X^{m+j_{\nu}}(z,w),
{\rm div} {\frak X}_{m+j_{\nu}}(z,w)\bigr)=0$
\ for all $(z,w)\in G,$  \ 
$\nu=s_2+1,\ldots, m.
\hfill
$
\\[2ex]
\indent
The  proof of the following statements is similar to those of Theorems 1.3, 1.4, and 1.5.
\vspace{0.75ex}

{\bf Theorem 1.9.}\!
{\it 
For the system of total differential equations {\rm (1.1)} to have a 
last mul\-ti\-p\-li\-er of the form {\rm (1.24)} 
it is necessary that {\rm (1.26)} be consistent on the domain $G.$
}
\vspace{0.75ex}

{\bf Corollary 1.8.}\!
{\it 
For the total differential system {\rm (1.1)} to have a 
$(s_1,0)\!$-non\-autonomous  $(n-k_1,n)\!$-cy\-lin\-dricality holomorphic
last multiplier of the form {\rm (1.24)} it is necessary that  the system of identities
{\rm (1.26)} with $s_2=0,\, k_2=0$ be consistent.
}
\vspace{0.75ex}

{\bf Corollary 1.9.}\!
{\it 
For the total differential system {\rm (1.1)} to have a 
$(0,s_2)\!$-non\-au\-to\-no\-mous  $(n,n-k_2)\!$-cy\-lin\-dricality antiholomorphic
last multiplier of the form 
{\rm (1.24)} it is necessary that the system of identities
{\rm (1.26)} with $s_1=0,\, k_1=0$ 
\vspace{0.75ex}
be consistent.
}

{\bf Corollary 1.10.}\!
{\it 
For the system 
{\rm (1.1)} to have an 
au\-to\-no\-mo\-us  $(n-k_1,n-k_2)\!$-cy\-lin\-dricality 
$\!{\mathbb R}\!$-differentiable last multiplier 
$
\mu\colon w\to  \mu({}^{k}\!w)$ 
for all $w\in \Omega^{\,\prime},
\ \Omega^{\,\prime}\subset {\mathbb C}^{n},
$
it is necessary that the system of identities {\rm (1.26)} with $s_1=0,\, s_2=0$ 
\vspace{1ex}
be consistent.
}

{\bf Theorem 1.10.}\!\! 
{\it
For the system of total differential equations {\rm (1.1)} to have at least one last multiplier 
of the form {\rm (1.24)} it is necessary and sufficient 
\vspace{0.35ex}
that there exist functions ${}^{s_1}\!\psi,\ {}^{s_2}\!\psi,\ {}^{\lambda}\!\varphi$ 
satisfying sys\-tem {\rm (1.14)} 
\vspace{0.5ex}
with 
\\[1.25ex]
\mbox{}\hfill                                                 
$
H_l\colon (z,w)\to {}-{\rm div}\, {\frak X}_{l}(z,w)$
\ for all $ (z,w)\in G,\ \ l=1,\ldots, 2m,
$ 
\hfill {\rm (1.27)}
\\[1.5ex]
such that the Pfaffian differential equation {\rm (1.15)} has the integrating factor 
$\nu({}^{s}\!z,{}^{k}\!w)=1$ for all $ (z,w)\in G^{\,\prime};$  
in this case the last multiplier is given by
\\[2ex]
\mbox{}\hfill          
$
\displaystyle
\mu\colon (z,w)\to 
\exp \int {}^{s_1}\!\psi({}^{s}\!z,{}^{k}\!w)\, d{}^{s_{\!1}}\!z +
{}^{s_2}\!\psi({}^{s}\!z,{}^{k}\!w)\, d\,\overline{{}^{s_{\!2}}\!z}\, +\,
{}^{k_1}\!\varphi({}^{s}\!z,{}^{k}\!w)\, d{}^{k_{\!1}}\!w +
{}^{k_2}\!\varphi({}^{s}\!z,{}^{k}\!w)\,d\,\overline{{}^{k_{2}}\!w}
\hfill
$
\\[1.75ex]
\mbox{}\hfill          
for all $ (z,w)\in G^{\,\prime}.
\hfill
$
}
\\[2.25ex]
\indent
{\bf Theorem 1.11.} 
\vspace{0.25ex}
{\it Let system {\rm (1.14)} with {\rm (1.27)} has $q$ not linearly bound solutions 
\vspace{0.25ex}
{\rm (1.17)} for which the corresponding Pfaff equations {\rm (1.18)} have the integrating factors  
\vspace{0.25ex}
$\nu_{\varepsilon}({}^{s}\!z,{}^{k}\!w)=1$
for all $ (z,w)\in G^{\,\prime},\ \varepsilon\!=\!1,\ldots, q.$   
Then the last multiplies of the total differential system {\rm (1.1)}
\\[2ex]
\mbox{}\hfill                       
$
\displaystyle
\mu_{\varepsilon}\! \colon\! (z, w)\to
\exp\!
\int {}^{s_1}\!\psi^{\varepsilon}({}^{s}\!z,{}^{k}\!w)\,
d{}^{s_{\!1}}\!z +
{}^{s_2}\!\psi^{\varepsilon}({}^{s}\!z,{}^{k}\!w)\,
d\,\overline{{}^{s_{\!2}}\!z} \,+\,
{}^{k_1}\!\varphi^{\varepsilon}({}^{s}\!z,{}^{k}\!w)\,
d{}^{k_{\!1}}\!w \,+\,
{}^{k_2}\!\varphi^{\varepsilon}({}^{s}\!z,{}^{k}\!w)\,
d\,\overline{{}^{k_{2}}\!w}
\hfill                       
$
\\[1.75ex]
\mbox{}\hfill                       
for all $(z,w)\in G^{\,\prime},
\quad \varepsilon=1,\ldots, q
\hfill
$
\\[1.5ex]
are functionally independent.
}
\vspace{1.25ex}

\newpage

The system of total differential equations 
\\[2ex]
\mbox{}\hfill                     
$
\begin{array}{l}
dw_1=w_1(1+2\,\overline{w}_2)\,dz +w_1(1+2w_2)\,d\,\overline{z},
\\[2ex]
dw_2=w_2(w_1-1)\,dz -w_2(w_2+\overline{w}_1)\,d\,\overline{z}
\end{array}
$
\hfill (1.28)
\\[2.25ex]
has the functions $({\rm div}\,{\frak X}_1(z,w)=1+2\,\overline{w}_2,  \
{\rm div}\,{\frak X}_2(z,w)=1+2w_2$ for all $(z,w)\in {\mathbb C}^3)$
\\[2.25ex]
\mbox{}\hfill
$
P_1\colon (z,w_1,w_2)\to \bigl(w_1(1+2\,\overline{w}_2), \, 1+2\,\overline{w}_2\bigr)$
\ for all $(z,w_1,w_2)\in {\mathbb C}^3
\hfill
$
\\[0.75ex]
and
\\[0.75ex]
\mbox{}\hfill
$
P_2\colon (z,w_1,w_2)\to 
\bigl(w_1(1+2w_2), \, 1+2w_2\bigr)$
\ for all $(z,w_1,w_2)\in {\mathbb C}^3.
\hfill
$
\\[2.25ex]
\indent
The Wronskians of the vector functions $P_1$  and $P_2$  with respect to 
$z,\ \overline{z},\ w_2,\ \overline{w}_1,$  
and $\overline{w}_2$ vanish identically on the ${\mathbb C}^3.$
\vspace{0.35ex}

Therefore a necessary condition for the total differential system  (1.28) 
to have an $\!{\mathbb R}$\!-dif\-fe\-ren\-ti\-a\-b\-le autonomous 
\vspace{0.35ex}
(1,2)-cylindricality last multiplier is complied (Theorem 1.9).

The scalar function 
\\[2ex]
\mbox{}\hfill
$
\varphi\colon (z,w_1,w_2)\to {}-\dfrac{1}{w_1}$
\ \ for all $(z,w_1,w_2)\in {\mathbb C}\times\Omega,
\hfill
$
\\[2.25ex]
where 
\vspace{0.5ex}
$\Omega$ is a domain from the set $\{(w_1,w_2)\colon w_1\ne 0\},$  
is a solution to system of equations (see (1.14) with  
$H_l(z,w_1,w_2)={}-{\rm div}\,{\frak X}_{l}(z,w_1,w_2)$ for all 
$(z,w_1,w_2)\in {\mathbb C}^3,\ l=1, 2)$
\\[2.25ex]
\mbox{}\hfill
$
w_1(1+2\,\overline{w}_2)\,\varphi= {}-(1+2\,\overline{w}_2),
\quad
2w_1\,\varphi={}-2,
\quad
w_1(1+2w_2)\,\varphi={}-(1+2w_2).
\hfill
$
\\[2ex]
\indent
Thus the total differential system (1.28) has the last multiplier (Theorem 1.10)
\\[2.25ex]
\mbox{}\hfill
$
\displaystyle
\mu\colon (z,w_1,w_2) \to \dfrac{1}{w_1}$
\ \  for all $(z_1,w_1,w_2)\in {\mathbb C}\times\Omega.
\hfill
$
\\[5.25ex]
\indent
{\bf  1.5. ${\mathbb R}\!$-regular solutions of an algebraic equation have no
movable nonalgebraic ${\mathbb R}\!$-singular point}
\\[0.75ex]
\indent
${\mathbb R}\!$-holomorphic solutions of  a completely solvable total differential equation  
may have 
\linebreak
${\mathbb R}\!$-sin\-gu\-lar points. 
In addition, we can distinguish two classes of ${\mathbb R}\!$-singular points 
of solutions: an ${\mathbb R}\!$-singular point of solutions of a completely solvable 
total differential equation 
whose position depends on the initial data determining a particular solution is referred to as 
a {\it movable} $\!{\mathbb R}\!$-singular point; 
if the position is independent of the initial data, then the point is called a 
{\it fixed} ${\mathbb R}\!$-singular point.
\vspace{0.5ex}

Let us consider the algebraic total differential equation
\\[2ex]
\mbox{}\hfill                               
$
\displaystyle
Q(z,w)\,dw = 
\sum\limits_{j=1}^m
\bigl(P_j(z,w)\,dz_j + P_{m+j}(z,w)\,d\,\overline{z}_j\bigr)\,,
$
\hfill (1.29)
\\[2ex]
where the functions $Q\colon G\to {\mathbb C}$ and 
$P_l\colon G\to {\mathbb C}, \ l=1,...,2m, \ G={\mathscr D}\times {\mathbb C},$
are ${\mathbb R}\!$-po\-ly\-no\-mi\-als in $w$ (po\-ly\-no\-mi\-als in $w$ and $\overline{w}\,)$
whose coefficients are 
${\mathbb R}\!$-holomorphic in $z$ in a domain ${\mathscr D}\subset {\mathbb C}^m$
and do not have common factors.
\vspace{1ex}

{\bf Definition 1.4.} 
{\it
Equation {\rm (1.29)} completely solvable in the domain $G$ 
is said to be nondegenerate if the rank of the matrix
\\[2ex]
\mbox{}\hfill
$
\displaystyle 
P(z,w)=
\left\|\!
\begin{array}{cccccc}
P_1(z,w) & \ldots & P_m(z,w) & P_{m+1}(z,w) & \ldots & P_{2m}(z,w)
\\[1.25ex]
\overline{P}_{m+1}(z,w) & \ldots & \overline{P}_{2m}(z,w) & 
\overline{P}_{1}(z,w) & \ldots & \overline{P}_{m}(z,w)
\end{array}
\!\right\|
\hfill
$
\\[2ex]
is equal to $2$ almost everywhere in $G$ and is said to be degenerate otherwise.
}
\vspace{0.75ex}

\newpage

By Definition 1.3, all ${\mathbb R}\!$-holomorphic solutions of a degenerate completely 
solvable equation (1.29) are ${\mathbb R}\!$-singular, and all 
${\mathbb R}\!$-singular solutions $w = w(z)$ of a nondegenerate completely solvable equation 
(1.29) satisfy the condition $P(z, w) < 2.$
\vspace{0.75ex}

{\bf Theorem 1.12.}\!\! 
{\it
${\mathbb R}\!$-holomorphic solutions of a nondegenerate completely solvable 
total differential equation 
{\rm (1.29)} have no movable nonalgebraic ${\mathbb R}\!$-singular points.
}
\vspace{0.35ex}

{\sl Proof}. 
Suppose the contrary: let $z_0\in {\mathscr D}$ be a nonalgebraic movable 
${\mathbb R}\!$-singular point for some solution $w = w(z)$ of 
the total differential equation (1.29), and let $\gamma \subset {\mathscr D}$ 
be the path along which the point $z$ 
tends to $z_0$ so that the solution $w = w(z)$ is ${\mathbb R}\!$-holomorphic on 
$\gamma$ everywhere except for the point $z_0.$ 
We have two possible cases: 

1) $z_0$ is a transcendental ${\mathbb R}\!$-singular point; 

2) $z_0$ is a $\triangle\!$-essentially ${\mathbb R}\!$-singular point.

In the first case, the solution $w = w(z)$ tends to some value $w_0\in \overline{{\mathbb C}}$ 
along the path $\gamma$ as $z\to z_0.$

If $w_0\in {\mathbb C}$, then we have two possibilities: 
a) the point $w_0$ is not a root of the equation
\\[1.5ex]
\mbox{}\hfill
$
Q(z_0,w) = 0;	
$
\hfill (1.30)
\\[1ex]
b) the point $w_0$ is a root of the equation (1.30).

By Theorem 1.1, in case 
a) the completely solvable total differential equation (1.29) has a solution 
$w = \widetilde{w}(z)\ {\mathbb R}\!$-holomorphic 
in a neighborhood of the point $z_0$ and satisfying the initial condition 
$\widetilde{w}(z_0) = w_0.$ 
Therefore, by Theorem 1.2, the solution $w = w(z)$ coincides with the solution 
$w = \widetilde{w}(z);$ 
consequently, $w = w(z)$ is ${\mathbb R}\!$-holomorphic at the point $z_0.$

Let us consider case b). 
Since $z_0$ is not a movable ${\mathbb R}\!$-singular point of the nondegenerate equation 
(1.29), we have ${\rm rank}\, P(z_0, w_0) = 2.$ 

We have the following three cases:
\vspace{0.35ex}

${\rm b}_1)$ there exist indices $k\in \{1,..., m\}$ and $\tau\in \{1,..., m\},\ k < \tau,$ such that
\\[1.75ex]
\mbox{}\hfill
$
P_{1k}(z_0,w_0)\,\overline{P}_{1,m+\tau}(z_0,w_0)-
P_{1\tau}(z_0,w_0)\,\overline{P}_{1,m+k}(z_0,w_0)\ne 0;
$
\hfill (1.31)
\\[2ex]
\indent
${\rm b}_2)$ there exist indices $k\in \{1,..., m\}$ and $\tau\in \{m+1,..., 2m\}$ such that
\\[1.75ex]
\mbox{}\hfill
$
P_{1k}(z_0,w_0)\,\overline{P}_{1,\tau-m}(z_0,w_0)-
P_{1\tau}(z_0,w_0)\,\overline{P}_{1,m+k}(z_0,w_0)\ne 0;
\hfill 
$
\\[2ex]
\indent
${\rm b}_3)$ there exist indices 
$k\in \{m+1,..., 2m\}$ and $\tau\in \{m+1,..., 2m\},\ k<\tau,$ such that
\\[1.75ex]
\mbox{}\hfill
$
P_{1k}(z_0,w_0)\,\overline{P}_{1,\tau-m}(z_0,w_0)-
P_{1\tau}(z_0,w_0)\,\overline{P}_{1,k-m}(z_0,w_0)\ne 0.
\hfill 
$
\\[2ex]
\indent
In case ${\rm b}_1),$ we rewrite the total differential equation (1.29) in the form
\\[2ex]
\mbox{}\ \ \ \                                
$
\displaystyle
P_k(z,w)\,dz_k + P_{\tau}(z,w)\,dz_{\tau}\, =\, 
Q(z,w)\,dw - P_{m+k}(z,w)\,d\,\overline{z}_k - 
P_{m+\tau}(z,w)\,d\,\overline{z}_{\tau} \ -
\hfill
$
\\[0.5ex]
\mbox{}\hfill (1.32)
\\
\mbox{}\hfill
$
\displaystyle
-\ \, \sum\limits_{j=1,\, j\ne k,\,  j\ne \tau}^m
\bigl(P_j(z,w)\,dz_j + P_{m+j}(z,w)\,d\,\overline{z}_j\bigr)\,.
\hfill 
$
\\[2.25ex]
\indent
By taking the conjugate of (1.32), we obtain the total differential equation
\\[2ex]
\mbox{}\ \                               
$
\displaystyle
\overline{P}_{m+k}(z,w)\,dz_k + \overline{P}_{m+\tau}(z,w)\,dz_{\tau}\, =\, 
\overline{Q}(z,w)\,d\,\overline{w} - \overline{P}_{k}(z,w)\,d\,\overline{z}_k - 
\overline{P}_{\tau}(z,w)\,d\,\overline{z}_{\tau} \ -
\hfill
$
\\[0.5ex]
\mbox{}\hfill (1.33)
\\
\mbox{}\hfill
$
\displaystyle
-\ \, \sum\limits_{j=1,\, j\ne k,\,  j\ne \tau}^m
\bigl(\overline{P}_{m+j}(z,w)\,dz_j + \overline{P}_{j}(z,w)\,d\,\overline{z}_j\bigr)\,.
\hfill 
$
\\[2ex]
\indent
Treating $Q$ and $P_{\,l}$ as the functions 
\\[1.5ex]
\mbox{}\hfill
$
Q(z,w) = q(z,\overline{z},w,\overline{w}\,)$ 
\ \ and \ \
$P_l(z,w) = p_l(z,\overline{z},w,\overline{w}\,),\ \ l = 1,..., 2m,
\hfill
$ 
\\[1.5ex]
holomorphic in $(z,\overline{z},w,\overline{w}),$ to 
differential system $(1.32)\cup (1.33)$ we assign 
the completely solvable system of total differential equations
\\[2ex]
\mbox{}\hfill                               
$
\displaystyle
p_k(t,x,y)\,dt_k + p_{\tau}(t,x,y)\,dt_{\tau}\, =\, 
q(t,x,y)\,dx - p_{m+k}(t,x,y)\,dt_{m+k} - 
p_{m+\tau}(t,x,y)\,dt_{m+\tau} \ -
\hfill
$
\\[2.5ex]
\mbox{}\hfill
$
\displaystyle
-\ \, \sum\limits_{j=1,\, j\ne k,\,  j\ne \tau}^m
\bigl(p_j(t,x,y)\,dt_j + p_{m+j}(t,x,y)\,dt_{m+j}\bigr),
\hfill 
$
\\
\mbox{}\hfill (1.34)
\\[0.75ex]
\mbox{}\hfill                           
$
\displaystyle
\overline{p}_{m+k}(t,x,y)\,dt_k + \overline{p}_{m+\tau}(t,x,y)\,dt_{\tau}\, =\, 
\overline{q}(t,x,y)\,dy - \overline{p}_{k}(t,x,y)\,dt_{m+k} - 
\overline{p}_{\tau}(t,x,y)\,dt_{m+\tau} \ -
\hfill
$
\\[2.5ex]
\mbox{}\hfill                           
$
\displaystyle
-\ \, \sum\limits_{j=1,\, j\ne k,\,  j\ne \tau}^m
\bigl(\overline{p}_{m+j}(t,x,y)\,dt_j + \overline{p}_{j}(t,x,y)\,dt_{m+j}\bigr)\,.
\hfill 
$
\\[2ex]
 \indent
Taking into account the complex analog of the results from [24, pp. 75 -- 80] 
and condition (1.31), we find that there exists a unique holomorphic solution 
\\[1.25ex]
\mbox{}\hfill
$
t_k = t_k(t_1,\ldots,t_{k-1},t_{k+1},\ldots, t_{\tau-1},t_{\tau+1},\ldots, t_{2m},x,y),
\hfill
$ 
\\[1.75ex]
\mbox{}\hfill
$
t_{\tau} = t_{\tau}(t_1,\ldots,t_{k-1},t_{k+1},\ldots, t_{\tau-1},t_{\tau+1},\ldots, t_{2m},x,y),
\hfill
$
\\[1.5ex]
of system (1.34) passing through the point $(t_0,x_0,y_0).$ 
Since system $(1.32)\cup (1.33)$ is self-adjoint, it follows that the equation 
(1.32) has ${\mathbb R}\!$-holomorphic integral manifolds 
\\[1.5ex]
\mbox{}\hfill
$
z_k- z_k(z, w) = 0
$ 
\ \ and \ \ 
$z_{\tau} - z_{\tau}(z, w) = 0
\hfill
$ 
\\[1.5ex]
passing through the point 
$(z_0, w_0).$ 
These manifolds are not determined by the equations 
$z_k=z_k^0$ and $z_{\tau} = z_{\tau}^0,$ respectively, since the function $Q$ 
is not identically zero at the point $z_0.$ 
Consequently, $z_0$ is an algebraic ${\mathbb R}\!$-singular point of the solution $w = w(z).$

Likewise, for cases ${\rm b}_2)$ and ${\rm b}_3)$
we can prove that $z_0$ cannot be a nonalgebraic ${\mathbb R}$-sin\-gu\-lar point of the 
solution $w = w(z).$

Let $w_0 = \infty.$ Performing the transformation $\xi = w^{{}-1},$ from the equation (1.29) 
we obtain a nondegenerate completely solvable equation; all functions occurring in this equation 
are ${\mathbb R}\!$-polynomials in $\xi$ (polynomials in $\xi$ and $\overline{\xi}\,)$ 
whose coefficients are ${\mathbb R}\!$-holomorphic in $z$ in the domain 
${\mathscr D}\subset {\mathbb C}^m$ and have no common factors. 
Just as in the case $w_0\in {\mathbb C}$, we find that for the solution $\xi=\xi(z)$ 
of this equation the point $z_0$ is either an ${\mathbb R}\!$-holomorphic point or a 
critical algebraic ${\mathbb R}\!$-singular point. 
Therefore, the solution $w = w(z)$ of the equation in question has either an 
${\mathbb R}\!$-pole or a critical ${\mathbb R}\!$-pole at the point $z_0.$

Thus, $z_0$ is not a transcendental ${\mathbb R}\!$-singular point of the solution 
$w = w(z)$ of the com\-p\-le\-te\-ly solvable total differential equation (1.29).
\vspace{0.25ex}

Let us now consider the case in which $z_0$ is a $\!\triangle\!$-essential 
${\mathbb R}\!$-singular point. If there exists at least one path $\gamma\subset {\mathscr D}$ 
that infinitely approaches the point $z_0$ and along which the solution $w = w(z)$ 
tends to some limit, then, just as above, we can prove that this is either an ordinary point or an 
algebraic point. Therefore, we assume that along any path $\gamma\subset {\mathscr D}$ 
the solution $w = w(z)$ does not tend to any limit as $z\to z_0.$ 
\vspace{0.25ex}
On one of such paths we choose a sequence of points 
$\bigl\{z^{(p)}\bigr\}_{p=1}^{{}+\infty}$ converging to the point $z_0$ as $p\to {}+\infty.$ 
\vspace{0.35ex}
The corresponding sequence of values of $w(z)$ has the form 
$\bigl\{w^{(p)}\bigr\}_{p=1}^{{}+\infty}.$ 
\vspace{0.35ex}
Since any sequence of complex numbers contains a subsequence converging to some 
\vspace{0.75ex}
number $w_0\in \overline{{\mathbb C}},$ it follows that without loss of generality 
we can assume that the sequence $\bigl\{w^{(p)}\bigr\}_{p=1}^{{}+\infty}$ 
itself converges to $w_0.$
\vspace{0.25ex}

Let $w_0\in {\mathbb C}.$ We consider two possibilities: 
a) the point $w_0$ is not a root of the equation (1.30); 
b) the point $w_0$ is a root of the equation (1.30).
\vspace{0.35ex}

By virtue of Theorem 1.1, in case a) the total differential equation (1.29) has the solution
\\[1.75ex]
\mbox{}\hfill  
$
w=\widetilde{w}{}^{(p)}(z)$
\ \ for \ \ 
$
\widetilde{w}{}^{(p)}\bigl(z^{(p)}\bigr)=w^{(p)},
$
\hfill (1.35)
\\[1.75ex]
where $\!\widetilde{w}{}^{(p)}\!(z)\!$ is a function $\!{\mathbb R}\!$-holomorphic in a 
neighborhood of the point $\!z_0\!$ provided that $\!p\!\!$ is a sufficiently large number. 
\!Therefore, by virtue of Theorem 1.1 and \!Corollary 1.1, the so- lu\-tion $\!w\! =\! w(z)\!$ 
coincides with the solution (1.35) and hence is $\!{\mathbb R}\!$-holomorphic at the $z_0.$

Case b). 
The point $z_0$ is not a fixed ${\mathbb R}\!$-singular point of solutions 
of the nondegenerate equation (1.29); therefore, ${\rm rank}\, P(z_0,w_0) = 2.$ 
Just as in the first case, we consider three possibilities, ${\rm b}_1),\ {\rm b}_2),$ and ${\rm b}_3).$

In case $\!{\rm b}_1\!),\!$ we construct the differential system $\!(1.32)\!\cup (1.33)\!$ and, using (1.31), 
conclude that the total dif\-fe\-ren\-ti\-al equation (1.29) has ${\mathbb R}\!$-holomorphic integral manifolds 
\\[1.5ex]
\mbox{}\hfill
$
z_k - z_k^{(p)}(z,w) = 0
$ 
\ \ and \ \ 
$
z_{\tau} - z_{\tau}^{(p)}(z,w) = 0
\hfill
$ 
\\[1.5ex]
passing through 
the point $\bigl(z^{(p)},w^{(p)}\bigr)$ and such that the functions 
$z_k^{(p)}(z,w)$ and $z_{\tau}^{(p)}(z,w)$ are ${\mathbb R}\!$-holomorphic in a 
\vspace{0.35ex}
neighborhood of the point $(z_0,w_0)$ for sufficiently large $p.$ 
We have 
\\[1.5ex]
\mbox{}\hfill
$
\lim\limits_{p\to {}+\infty} z_k^{(p)} = z_k^0
$ 
\ \ \, and \, \ \ 
$\lim\limits_{p\to {}+\infty} z_{\tau}^{(p)} = z_{\tau}^0.
\hfill
$
\\[1.5ex]
\indent
Let $\gamma_0$ be the path in the complex plane $w$ corresponding to the solution $w = w(z)$ 
as the point $z$ goes along the path $\gamma.$ The path $\gamma$ can be chosen so that relations 
of the form (1.31) are valid on $\gamma$ and $\gamma_0$ including the point $(z_0,w_0).$ 
Then the functions $z_k^{(p)}(z,w)$ and $z_{\tau}^{(p)}(z,w)$ are ${\mathbb R}\!$-holomorphic 
\vspace{0.35ex}
along the path $\gamma_0\times\gamma$ for sufficiently large $p.$ 

Therefore, $z_{k}^{(p)}(z_0,w_0) = z_{k}^0$ and $z_{\tau}^{(p)}(z_0,w_0) = z_{\tau}^0$ 
\vspace{0.35ex}
for sufficiently large $p.$ 

Since the functions $z_{k}^{(p)}(z,w)$ and 
$z_{\tau}^{(p)}(z,w)$ are ${\mathbb R}\!$-holomorphic in a neighborhood of the point 
\vspace{0.25ex}
$(z_0,w_0)$ and the total dif\-fe\-ren\-ti\-al 
 equation (1.29) has a unique ${\mathbb R}\!$-holomorphic solution 
\vspace{0.25ex}
with the initial data $(z_0,w_0),$ we find that the identities 
\\[1.5ex]
\mbox{}\hfill
$
z_{k}^{(p)}(z,w)\equiv z_{k}(z,w)
$ 
\ \ and \ \ 
$
z_{\tau}^{(p)}(z,w)\equiv z_{\tau}(z,w)
\hfill
$
\\[1.5ex]
are valid for all sufficiently large $p.$ 
\vspace{0.25ex}
Consequently, the solution $w = w(z)$ is ${\mathbb R}\!$-holomorphic along the path $\gamma$ 
except for the point $z$ and satisfies the equations 
\\[1.5ex]
\mbox{}\hfill
$
z - z_k^{(p)}(z,w) = 0
$ 
\ \ and \ \ 
$z - z_{\tau}^{(p)}(z,w) = 0.
\hfill
$
\\[1.5ex]
\indent
Therefore, $z_0$ is an algebraic point for this solution.
\vspace{0.25ex}

In a similar way, we can show that in cases ${\rm b}_2)$ and ${\rm b}_3)$ 
the point $z_0$ cannot be a nonalgebraic ${\mathbb R}\!$-singular point of the solution $w = w(z).$
\vspace{0.35ex}

Now let $w_0 =\infty.$ 
Then, by setting $\xi=w^{{}-1}$ in the total differential equation (1.29), we find that the solution 
$\xi=\xi(z)$ of the obtained equation has an algebraic ${\mathbb R}\!$-singularity at the point $z_0.$ 
Therefore, $z_0$ is an algebraic point for the solution $w = w(z)$ of 
the completely solvable total differential equation (1.29). 
The proof of the theorem is complete. \k
\\[5.75ex]
\centerline{
\large\bf  
2. System of first-order partial differential equations}
\\[1.5ex]
\indent
{\bf 2.1. $\!{\mathbb R}\!$-differentiable integrals and last multipliers}
\\[0.75ex]
\indent
Consider a linear homogeneous system of  first-order partial differential equations
\\[1.75ex]
\mbox{}\hfill                           
$
{\frak A}_{j} (z) \, u \ = \ 0, 
\quad j=1,\ldots,m,
$
\hfill (2.1)
\\[1.5ex]
with not linearly bound [25, p. 105] differential operators 
\\[1.75ex]
\mbox{}\hfill                          
$
\displaystyle
{\frak A}_{j}(z)=
\sum\limits_{\xi=1}^{n}\bigl( u_{j\xi}(z)\partial_{z_{\xi}} +
u_{j,n+\xi}(z)\partial_{{}_{\scriptstyle \overline{z}_{\xi}}}\bigr)$
\ for all $z\in G, \quad j=1,\ldots,m,
\hfill                      
$
\\[1.75ex]
where the scalar functions 
\vspace{0.25ex}
$u_{jp}\colon G\to {\mathbb C},\ 
j=1,\ldots, m,\ p=1,\ldots, 2n,$ are ${\mathbb R}\!$-differentiable 
in a domain $G\subset {\mathbb C}^{n},$
the $\overline{z}_j$  are the complex conjugates of  $z_j, \ j=1,\ldots,m.$ 
\vspace{0.5ex}

We begin with definitions. 
\vspace{0.35ex}
An $\!{\mathbb R}\!$-differentiable on a domain $\!G^{\,\prime}\!\subset\! G\!$  function:  
i)  $\!\!F\colon\! G^{\,\prime}\!\to\! {\mathbb C};\!\!$
ii) $\!f\colon G^{\,\prime}\!\to {\mathbb C};$
\vspace{0.75ex}
iii) $\!\mu\colon G^{\,\prime}\!\to {\mathbb C}$ 
is called  \ i) a {\it first integral}; ii) a {\it partial integral}; iii) a {\it last mul\-ti\-p\-lier} 
of the partial differential system (2.1)  iff \,
i)  $\!{\frak A}_{j}F(z)=0$ for all $z\in G^{\,\prime}, \, j=1,\ldots, m;$
\\[0.75ex]
\indent
ii)  ${\frak A}_{j}f(z)=\Phi_{j}(f; z)$ for all $ z\in G^{\,\prime},$ 
where $\Phi_{j}(0; z)\equiv 0, \ j=1,\ldots, m;$
\\[1ex]
\indent
iii) $\!{\frak A}_{j}\mu(z)={}-\mu(z)\,{\rm div}\,{\frak u}^{j}(z)$ for all 
$z\in G^{\,\prime},$ where the vector functions
\\[1.75ex]
\mbox{}\hfill                                  
$
u^{j}\colon z\to (u_{j1}(z),\ldots,u_{j,2n}(z))
$
\ for all $z\in G,\ \ j=1,\ldots, m.
\hfill
$
\\[1.75ex]      
\indent
The ${\mathbb R}\!$-differentiable first integral $F$  (partial integral $f$  and last multiplier $\mu)$ 
of the partial differential system (2.1) is called 
{\it $(n-k_1,n-k_2)\!$-cylindricality} [20; 26; 27] if
\\[0.5ex]
\indent
(i)   $F \ (f$ and $\mu)$ is holomorphic of  $n-k_2$ variables;
\\[0.5ex]
\indent
(ii)  $F \ (f$ and $\mu)$ is antiholomorphic of  $n-k_1$ variables.
\\[2ex]
\indent
{\bf 2.1.1. {\boldmath $\!\!(n-k_1,n-k_2)\!$}-cy\-lin\-dricality  partial integrals}.                                  
\vspace{0.5ex}
Suppose the system (2.1) has an ${\mathbb R}\!$-dif\-ferentiable 
$(n-k_1,n-k_2)\!$-cy\-lin\-dricality  partial integral 
\\[1.75ex]
\mbox{}\hfill                              
$
f \colon z\to  f({}^{k}\!z)$
\ for all $z\in G^{\,\prime},
$
\hfill (2.2)
\\[1.75ex]
where $k=(n-k_1,n-k_2).$ 
\vspace{0.5ex}
Without loss of generality it can be assumed that 
the function $f$ is an antiholomorphic function of
\vspace{0.75ex}
$z_{k_1+1},\ldots, z_{n}$ and the function $f$ is a holomorphic function of 
\vspace{0.75ex}
$z_{{}_{\scriptstyle \zeta_{k_2+1}}},\ldots, z_{{}_{\scriptstyle \zeta_{n}}},\ 
\zeta_{\delta}\in \{1,\ldots, n\},\  \delta=k_2+1,\ldots, n.$

Then, in accordance with the definition of a partial integral for the system (2.1), 
\\[1.75ex]
\mbox{}\hfill                              
$
{\frak A}_{j}^{k} f({}^{k}\!z) = \Phi_{j}(f; z)$
\ for all $z\in G^{\,\prime},
\ \ j=1,\ldots, m,
$
\hfill (2.3)
\\[1.75ex]
where the linear differential operators of first order
\\[1.75ex]
\mbox{}\hfill                                  
$
\displaystyle
{\frak A}_{j}^{k}(z)=
\sum\limits_{\xi=1}^{k_1}u_{{}_{\scriptstyle j\xi}}(z)\,\partial_{{}_{\scriptstyle z_{\xi}}} +
\sum\limits_{\tau=1}^{k_2}u_{{}_{\scriptstyle j\zeta_{\tau}}}(z)\,
\partial_{{}_{\scriptstyle \overline{z}_{\zeta_{\tau}}}}
$
\ for all $z\in G,
\hfill
$
\\[1.75ex]
the indexes $\zeta_{\tau}\in\{1,\ldots,n\},\ \tau=1,\ldots, k_2,$ 
the functions 
\\[1.75ex]
\mbox{}\hfill                                  
$
\Phi_{j}(0; z)=0
$ 
\ for all $z\in G,\ \ \ j=1,\ldots, m.
\hfill
$
\\[1.75ex]
\indent
Let the system of identities (2.3) hold. Then the functions from the sets
\\[1.5ex]
\mbox{}\hfill                            
$
\bigl\{
u_{j1}(z),\ldots,u_{jk_1}(z),
u_{j\zeta_{1}}(z),\ldots,u_{j\zeta_{{}_{\scriptstyle k_2}}}(z)
\bigr\},
\ \ j=1,\ldots, m,
$
\hfill (2.4)
\\[1.5ex]
are linearly bound on the integral manifold 
\\[1.5ex]
\mbox{}\hfill                       
$
f({}^{k}\!z)=0.
$
\hfill (2.5)
\\[1.5ex]
\indent
Therefore the Wronskians of the functions from the sets 
(2.4) with respect to  
\vspace{0.35ex}
$z_{\gamma}, \ \overline{z}_{\zeta_{\delta}},$ 
$\gamma=k_1+1,\ldots, n,\ \delta=k_2+1,\ldots, n$
vanish identically on the manifold (2.5), i.e., 
\\[2ex]
\mbox{}\hfill                                       
$
W_{{}_{\scriptstyle z_{\gamma}}}
\bigl({}^{\lambda}u^{j}(z)\bigr)\!=
{\stackrel{*}{\Psi}}_{j\gamma}(f;z)$
\ for all $z\in G,
\ \ j=1,\ldots, m,
\ \gamma=k_1+1,\ldots, n,
\hfill
$
\\
\mbox{}\hfill (2.6)
\\
\mbox{}\hfill
$
W_{{}_{\scriptstyle \overline{z}_{\zeta_{\delta}}}}
\bigl({}^{\lambda}u^{j}(z)\bigr)\!=
{\stackrel{**}{\Psi}}_{j\zeta_{\delta}}(f;z)$
\ for all $z\in G,
\ \ j=1,\ldots, m,
\ \delta=k_2+1,\ldots, n,
\hfill
$
\\[2.25ex]
where $W_{{}_{\scriptstyle z_{\gamma}}}$ and
$W_{{}_{\scriptstyle \overline{z}_{\zeta_{\delta}}}}$ are 
the Wronskians with respect to $z_{\gamma}$ and to $\overline{z}_{\zeta_{\delta}}$
respectively, 
the func\-ti\-ons
${\stackrel{*}{\Psi}}_{j\gamma}\colon G\to {\mathbb C},\
{\stackrel{**}{\Psi}}_{j \zeta_{\delta}}\colon G\to {\mathbb C}$
are $\!{\mathbb R}\!$-differentiable 
\vspace{0.35ex}
on the domain $\!G$ and 
$\!{\stackrel{*}{\Psi}}_{j\gamma}\!(0;z)\!\equiv\! 0,$
${\stackrel{**}{\Psi}}_{j \zeta_{\delta}}\!(0;z)\equiv 0,
\ \gamma=k_1+1,\ldots, n,\ \delta=k_2+1,\ldots, n, \ j=1,\ldots, m,$
\vspace{1.25ex}
the number $\lambda=k_1+k_2,$ 
the functions
\vspace{0.5ex}
${}^{\lambda}\!u^{j}\colon z\to
\bigl(u_{j1}(z),\ldots,u_{j k_1}(z),
u_{j\,\zeta_1}(z),\ldots,
u_{{}_{\scriptstyle j\,\zeta_{k_2}}}(z)\bigr)$
for all $z\in G,\ j=1,\ldots, m.$

Thus the following statements are valid.
\vspace{0.75ex}

{\bf Theorem 2.1.}\!
{\it 
For the system of partial differential equations {\rm (2.1)} to have an 
${\mathbb R}\!$-dif\-fe\-ren\-ti\-able
$(n-k_1,n-k_2)\!$-cy\-lin\-dricality partial in\-teg\-ral of the form 
{\rm (2.2)} it is necessary that  the system of identities
{\rm (2.6)} be consistent.
}

\newpage

{\bf Corollary 2.1.}\!\!
{\it 
For the linear homogeneous system of partial differential 
equations {\rm (2.1)} to have an $\!(n\!-\!k_1,n)\!$-cy\-lin\-dricality 
\vspace{0.25ex}
holomorphic partial in\-teg\-ral of the form 
{\rm (2.2)} it is necessary that  the system of identities
{\rm (2.6)} with $k_2=0$ be consistent.
}
\vspace{0.75ex}

{\bf Corollary 2.2.}\!
{\it 
For the linear homogeneous system of partial differential equations {\rm (2.1)} to have an 
$(n,n-k_2)\!$-cy\-lin\-dricality antiholomorphic
\vspace{0.25ex}
partial in\-teg\-ral of the form 
{\rm (2.2)} it is necessary that the system of identities
{\rm (2.6)} with $k_1=0$ be consistent.
}
\vspace{0.75ex}

Let the system (2.1) satisfy conditions (2.6). Let us write out the system of equations
\\[2ex]
\mbox{}\hfill                                       
$
{}^{\lambda}\!\varphi
\bigl({}^{\lambda}u^{j}(z)\bigr)^{T}= H_{j}(f;z),
\ \ j=1,\ldots, m,
\hfill
$
\\[2ex]
\mbox{}\hfill
$
{}^{\lambda}\!\varphi
\bigl(\partial_{{}_{\scriptstyle z_{\gamma}}}^{l}
{}^{\lambda}u^{j}(z)\bigr)^{T}=
\partial_{{}_{\scriptstyle z_{\gamma}}}^{l}\!
H_{j}(f;z),
\ \
l=1,\ldots, \lambda-1, \
\gamma=k_1+1,\ldots, n, \
\ j=1,\ldots, m,
$
\hfill (2.7)
\\[2.25ex]
\mbox{}\hfill
$
{}^{\lambda}\!\varphi
\bigl(
\partial_{{}_{\scriptstyle \overline{z}_{\zeta_{\delta}}}}^{l}
{}^{\lambda}u^{j}(z)\bigr)^{T}=
\partial_{{}_{\scriptstyle \overline{z}_{\zeta_{\delta}}}}^{l}
H_{j}(f;z),
\ \
l=1,\ldots, \lambda-1, \
\delta=k_2+1,\ldots, n, \
\ j=1,\ldots, m,
\hfill
$
\\[2ex]
where the vector functions 
\\[2ex]
\mbox{}\hfill
$
{}^{k_1}\!\varphi\colon z\to
(\varphi_{1 k_1}\!({}^{k}\!z),\ldots, \varphi_{k_1 k_1}\!({}^{k}\!z)),
\qquad 
{}^{k_2}\!\varphi\colon\! z\!\to\!
(\varphi_{1 k_2}\!({}^{k}\!z),\ldots, \varphi_{k_2 k_2}\!({}^{k}\!z)),\!\!
\hfill
$
\\[2.5ex]
\mbox{}\hfill
$
{}^{\lambda}\!\varphi\colon z\to ({}^{k_1}\!\varphi (z),{}^{k_2}\!\varphi (z))
$
\ \ for all $z\in G,
\hfill
$
\\[2ex]
the scalar functions $H_{j}\colon G\to {\mathbb C}$ are 
${\mathbb R}\!$-differentiable on the domain $\!G$ and
$H_{j}(0;z) = 0$ for all $z\in G, \ j=1,\ldots,m.$ 
\vspace{0.35ex}

Let us introduce the Pfaffian differential equation
\\[2ex]
\mbox{}\hfill                              
$
{}^{k_1}\!\varphi({}^{k}\!z)\, d{}^{k_{1}}\!z\ +\
{}^{k_2}\!\varphi({}^{k}\!z)\,d\,\overline{{}^{k_{2}}\!z}\ =\ 0,
$
\hfill (2.8)
\\[2ex]
where the vector columns $d{}^{k_1}\!z=\mbox{colon}(dz_1,\ldots,dz_{k_1}), \
d\,\overline{{}^{k_{2}}\!z}=\mbox{colon}\bigl(
d\,\overline{z}_{{}_{\scriptstyle \zeta_1}},\ldots,
d\,\overline{z}_{{}_{\scriptstyle \zeta_{k_2}}}\bigr).$
\vspace{0.75ex}

{\bf Theorem 2.2.}\! 
\vspace{0.25ex}
{\it
A necessary and sufficient condition for the partial differential system {\rm (2.1)} 
to have at least one 
\vspace{0.35ex}
$\!{\mathbb R}\!$-differentiable partial integral of the form {\rm (2.2)} 
is that the functions ${}^{\lambda}\!\varphi\colon G\to {\mathbb C}^{\lambda}$  and     
$H_j\colon G\to {\mathbb C},\ j=1,\ldots, m,$ exist so that they satisfy system {\rm (2.7)} and
\vspace{0.5ex}

{\rm (i)}\ the Pfaff equation {\rm (2.8)} has an integrating factor{\rm;}
\vspace{0.5ex}

{\rm (ii)}\ the function {\rm (2.2)} is a general integral of the Pfaffian equation {\rm (2.8)}.  
}
\vspace{0.75ex}

{\sl Proof. Necessity}.\! 
Let the partial differential system (2.1) have a 
$\!{\mathbb R}\!$-differentiable partial in\-teg\-ral of the form (2.2). 
Then the identity (2.3) holds. The vector functions 
\\[2ex]
\mbox{}\hfill
$
{}^{k_1}\varphi\colon z\to\,
\partial_{{}_{\scriptstyle {}^{k_1}\!z}}
f({}^{k}\!z)$
\  for all $z\in G^{\,\prime},
\qquad
{}^{k_2}\varphi\colon z\to\,
\partial_{{}_{\scriptstyle \overline{{}^{k_2}\!z}}}\,
f({}^{k}\!z)$
\  for all $z\in G^{\,\prime},
\hfill
$
\\[2.25ex]
where 
\vspace{0.35ex}
$\partial_{{}_{\scriptstyle {}^{k_{\!1}}\!z}}\!\!=\!
(\partial_{z_1},\!\ldots\!,\partial_{z_{k_1}}\!), \
\partial_{{}_{\scriptstyle \overline{{}^{k_{\!2}}z}}} \!=\!
\bigl(\partial_{{}_{\scriptstyle \overline{z}_{\zeta_1}}}\!,\!\ldots\!,
\partial_{{}_{\scriptstyle \overline{z}_{\zeta_{k_2}}}}\!\!\bigr),$
is a solution to system (2.7),
\vspace{0.35ex}
which can be shown by dif\-fe\-ren\-ti\-a\-ting (2.3)  
\vspace{0.35ex}
$\lambda-1$ times with respect to  $z_{\gamma},\  \gamma=k_1+1,\ldots, n,$ 
and $\lambda-1$ times with respect to 
$\overline{z}_{{}_{\scriptstyle \zeta_{\delta}}},\ \delta=k_2+1,\ldots, n.$
\vspace{0.35ex}
Therefore the scalar function (2.2) is a 
general integral of the Pfaffian differential equation (2.8). 
\vspace{0.75ex}

{\sl Sufficiency}. 
\vspace{0.35ex}
Let ${}^{\lambda}\!\varphi$ be a solution to system (2.7), and let the corresponding Pfaff equation 
(2.8) have an integrating factor $\mu\colon {}^{k}\!z\to \mu({}^{k}\!z)$
and the corresponding general integral (2.2). 

Then 
\\[1.5ex]
\mbox{}\hfill                        
$
\partial_{{}_{\scriptstyle {}^{k_1}\!z}}\,f({}^{k}\!z)\, -\,
\mu ({}^{k}\!z)\,{}^{k_1}\varphi({}^{k}\!z)=0,
\qquad
\partial_{{}_{\scriptstyle \overline{{}^{k_2}\!z}}}\,f({}^{k}\!z) \, -\, 
\mu ({}^{k}\!z)\,{}^{k_2}\varphi({}^{k}\!z)=0.
$
\hfill (2.9)
\\[2ex]
\indent
It follows from (2.7) and (2.9) that identity (2.3) is valid with
\\[1.75ex]
\mbox{}\hfill
$
\Phi_{j}(f; z)=\mu({}^{k}\!z)\,H_{j}(f; z)$
\ for all $z\in G^{\,\prime},\ \ \  j=1,\ldots, m.
\hfill
$
\\[1.75ex]
\indent
Consequently the function (2.2) is a partial integral of the system (2.1).\, \k

\newpage

Consider the linear homogeneous system of partial differential equations 
\\[2ex]
\mbox{}\hfill                              
$
{\frak A}_1^{}(z_1^{},z_2^{})\;\!u=0,
\qquad  
{\frak A}_2^{}(z_1^{},z_2^{})\;\!u=0,
$
\hfill (2.10)
\\[2ex]
where the linear differential operators of first order
\\[2.25ex]
\mbox{}\hfill
$
{\frak A}_1^{}(z_1^{},z_2^{}) =
z_1^{}(z_2^{}+\overline{z}_1)\,\partial_{z_1^{}} + 
z_2^{}(z_2^{}+\overline{z}_1)\,\partial_{z_2^{}} + 
(z_1^2+z_2^2+\overline{z}\;\!{}_1^2+\overline{z}\;\!{}_2^2)\,
\partial_{{}_{\scriptstyle \overline{z}_1}} +
(z_1^2-z_2^2+\overline{z}\;\!{}_1^2-\overline{z}\;\!{}_2^2)\,
\partial_{{}_{\scriptstyle \overline{z}_2}} ,
\hfill
$
\\[2.5ex]
\mbox{}\hfill
$
{\frak A}_2^{}(z_1^{},z_2^{}) =
z_1^{}(\overline{z}_1+\overline{z}_2)\,\partial_{z_1^{}} + 
z_2^{}(\overline{z}_1+\overline{z}_2)\,\partial_{z_2^{}} +
(z_1^2-z_2^2+\overline{z}\;\!{}_1^2-\overline{z}\;\!{}_2^2)\,
\partial_{{}_{\scriptstyle \overline{z}_1}} +
(z_1^2+z_2^2+\overline{z}\;\!{}_1^2+\overline{z}\;\!{}_2^2)\,
\partial_{{}_{\scriptstyle \overline{z}_2}} 
\hfill
$
\\[2.25ex]
\mbox{}\hfill
 for all $(z_1^{},z_2^{})\in {\mathbb C}^2.
\hfill
$
\\[1.5ex]
\indent
Let us find a (0,2)-cy\-lin\-dricality holomorphic partial integral of system (2.10). 
\vspace{0.25ex}
The Wronskians of the sets of functions 
\vspace{0.5ex}
$U_1=\{z_1^{}(z_2^{}+\overline{z}_1),z_2^{}(z_2^{}+\overline{z}_1)\}$ and 
$U_2=\{z_1^{}(\overline{z}_1+\overline{z}_2), z_2^{}(\overline{z}_1+\overline{z}_2)\}$
with respect to   $\overline{z}_1^{}$ and $\overline{z}_2^{}$ vanish identically  
on the space ${\mathbb C}^2\colon$
\\[2ex]
\mbox{}\hfill
$
W_{{}_{\scriptstyle \overline{z}_1}}(z_1^{}(z_2^{}+\overline{z}_1),z_2^{}(z_2^{}+\overline{z}_1))=
\left|\!\!
\begin{array}{cc}
z_1^{}(z_2^{}+\overline{z}_1) & z_2^{}(z_2^{}+\overline{z}_1)
\\[0.5ex]
z_1^{} & z_2^{}
\end{array}
\!\!\right| =0
$
\ for all $(z_1^{},z_2^{})\in {\mathbb C}^2,
\hfill
$
\\[2.5ex]
\mbox{}\hfill
$
W_{{}_{\scriptstyle \overline{z}_2}}(z_1^{}(z_2^{}+\overline{z}_1),z_2^{}(z_2^{}+\overline{z}_1))=0
$
\ for all $(z_1^{},z_2^{})\in {\mathbb C}^2,
\hfill
$
\\[2.5ex]
\mbox{}\hfill
$
W_{{}_{\scriptstyle \overline{z}_1}}(z_1^{}(\overline{z}_1+\overline{z}_2), z_2^{}(\overline{z}_1+\overline{z}_2))=
\left|\!\!
\begin{array}{cc}
z_1^{}(\overline{z}_1+\overline{z}_2) & z_2^{}(\overline{z}_1+\overline{z}_2)
\\[0.5ex]
z_1^{} & z_2^{}
\end{array}
\!\!\right| =0
$
\ for all $(z_1^{},z_2^{})\in {\mathbb C}^2,
\hfill
$
\\[2.75ex]
\mbox{}\hfill
$
W_{{}_{\scriptstyle \overline{z}_2}}(z_1^{}(\overline{z}_1+\overline{z}_2), z_2^{}(\overline{z}_1+\overline{z}_2))=
\left|\!\!
\begin{array}{cc}
z_1^{}(\overline{z}_1+\overline{z}_2) & z_2^{}(\overline{z}_1+\overline{z}_2)
\\[0.5ex]
z_1^{} & z_2^{}
\end{array}
\!\!\right| =0
$
\ for all $(z_1^{},z_2^{})\in {\mathbb C}^2.
\hfill
$
\\[2ex]
\indent
Therefore the necessary conditions for the partial differential system (2.10) to have an 
holomorphic partial integral is complied (Corollary 2.1).

Let us write the functional system {\rm (2.7)}:
\\[1.75ex]
\mbox{}\hfill
$
z_1^{}(z_2^{}+\overline{z}_1)\,\varphi_1^{} + z_2^{}(z_2^{}+\overline{z}_1)\,\varphi_2^{}=
(z_1^{}+z_2^{})(z_2^{}+\overline{z}_1),
\qquad
z_1^{}\,\varphi_1^{} + z_2^{}\,\varphi_2^{}=z_1^{}+z_2^{},
\hfill
$
\\[2.5ex]
\mbox{}\hfill
$
z_1^{}(\overline{z}_1+\overline{z}_2)\,\varphi_1^{} + 
z_2^{}(\overline{z}_1+\overline{z}_2)\,\varphi_2^{}=
(z_1^{}+z_2^{})(\overline{z}_1+\overline{z}_2),
\qquad
z_1^{}\,\varphi_1^{} + z_2^{}\,\varphi_2^{}=z_1^{}+z_2^{},
\hfill
$
\\[2.25ex]
where 
\vspace{1.25ex}
$H_1^{}\colon (z_1^{},z_2^{})\to (z_1^{}+z_2^{})(z_2^{}+\overline{z}_1),\
H_2^{}\colon (z_1^{},z_2^{})\to (z_1^{}+z_2^{})(\overline{z}_1+\overline{z}_2)$
for all $(z_1^{},z_2^{})\in {\mathbb C}^2.$

The functions 
\vspace{0.75ex}
$
\varphi_1^{}\colon (z_1^{},z_2^{})\to 1
$
for all $(z_1^{},z_2^{})\in {\mathbb C}^2,
\
\varphi_2^{}\colon (z_1^{},z_2^{})\to 1
$ 
for all $(z_1^{},z_2^{})\in {\mathbb C}^2
$
is a solution to this system.
The corresponding Pfaffian differential equation
\\[1.75ex]
\mbox{}\hfill
$
dz_1^{} + dz_2^{} =0
\hfill
$
\\[1.75ex]
has the integrating factor $\mu\colon (z_1^{},z_2^{})\to 1$ for all $(z_1^{},z_2^{})\in {\mathbb C}^2$ and 
the general integral 
\\[1.75ex]
\mbox{}\hfill                                       
$
f\colon (z_1^{},z_2^{}) \to z_1^{}+z_2
$
\ for all $(z_1^{},z_2^{})\in {\mathbb C}^2.
$
\hfill (2.11)
\\[1.75ex]
\indent
By Theorem 2.2, the function (2.11) is a holomorphic partial integral of system (2.10).
\vspace{1.25ex}

{\bf Theorem 2.3.} 
{\it
Let $h$ systems {\rm (2.7)} have $q$ not linearly bound solutions
\\[1.75ex]
\mbox{}\hfill                       
$
{}^{\lambda}\varphi^{\varepsilon}\colon z\to
{}^{\lambda}\varphi^{\varepsilon}({}^{k}\!z)$
\ for all $z\in G^{\,\prime},
\quad \varepsilon=1,\ldots,q,
$
\hfill {\rm (2.12)}
\\[1.75ex]
for which the corresponding Pfaffian differential equations
\\[2.5ex]
\mbox{}\hfill                              
$
{}^{k_1}\!\varphi^{\varepsilon}({}^{k}\!z)\,
d{}^{k_{1}}\!z \ + \
{}^{k_2}\!\varphi^{\varepsilon}({}^{k}\!z)\,
d\,\overline{{}^{k_{2}}\!z}\, =\, 0,
\quad
\varepsilon=1,\ldots, q
$
\hfill {\rm(2.13)}
\\[2ex]
have the general $\!{\mathbb R}\!$-differentiable integrals     
$f_{\varepsilon}\colon z\to f_{\varepsilon}({}^{k}\!z)$
for all $z\in G^{\,\prime},\ \varepsilon=1,\ldots, q.$
Then these integrals are functionally independent. 
}
\vspace{0.75ex}

{\sl Proof}. Using (2.9), we have 
\\[2ex]
\mbox{}\hfill
$
\partial_{{}_{\scriptstyle {}^{k_1}\!z}}\,f_{\varepsilon}({}^{k}\!z) =
\mu_{\varepsilon}({}^{k}\!z)\,
{}^{k_1}\varphi^{\varepsilon}({}^{k}\!z),
\quad
\partial_{{}_{\scriptstyle \overline{{}^{k_2}\!z}}}\,f_{\varepsilon}({}^{k}\!z)=
\mu_{\varepsilon}({}^{k}\!z)\,
{}^{k_2}\varphi^{\varepsilon}({}^{k}\!z)$
\ for all 
$
z\in G^{\prime},
\quad
\varepsilon=1,\ldots, q.
\hfill
$
\\[2ex]
\indent
Therefore the Jacobi matrix
$
J(f_{\varepsilon}({}^{k}\!z);{}^{k}\!z)=
\bigl\|
{}^{k_1}\!\Phi({}^{k}\!z)\,
{}^{k_2}\Phi({}^{k}\!z)
\bigr\|,
$
where
\vspace{1.5ex}
${}^{k_1}\!\Phi =
\bigl\|\mu_{\varepsilon}\varphi_{\xi k_1}^{\varepsilon }\!\bigr\|\!$
is a $(q\times k_1)\!$-matrix,
\vspace{1ex}
${}^{k_2}\Phi=\bigl\|\mu_{\varepsilon}\varphi_{\tau k_2}^{\varepsilon}\bigr\|$ is a 
$(q\times k_2)\!$-matrix. 
Since the solutions (2.12) are not linearly bound, it follows that 
${\rm rank}\,J(f_{\varepsilon}({}^{k}\!z);{}^{k}\!z)=q.$ 
\vspace{0.5ex}
Thus the general integrals of the Pfaffian differential equations (2.13) are functionally independent. \k
\vspace{1.75ex}

{\bf 2.1.2. {\boldmath $\!(n-k_1,n-k_2)\!$}-cylindricality first integrals}.
Let the function
\\[1.75ex]
\mbox{}\hfill                              
$
F\colon z\to  F({}^{k}\!z)$
\ for all $z\in G^{\,\prime}
$
\hfill (2.14)
\\[1.75ex]
be an ${\mathbb R}\!$-differentiable $(n-k_1,n-k_2)\!$-cylindricality first integral of system (2.1).

Then, in accordance with the criteria of an $\!{\mathbb R}\!$-differentiable first integral,
\\[1.75ex]
\mbox{}\hfill                          
$
{\frak A}_{j}^{k}\,F({}^{k}\!z) = 0$
\ for all $z\in G^{\,\prime},
\quad 
j=1,\ldots, m.
\hfill
$
\\[1.75ex]
\indent
Hence the Wronskians of the functions from the sets (2.4) vanish identically on the domain $G,$ i.e., 
the system of identities (2.6) with  
\\[1.5ex]
\mbox{}\hfill    
$
{\stackrel{*}{\Psi}}_{j\gamma}(z)\, =\,
{\stackrel{**}{\Psi}}_{j \zeta_{\delta}}(z)\,=0,
\ \ \gamma=k_1+1,\ldots, n, \
\delta=k_2+1,\ldots, n,\ 
j=1,\ldots, m,
$
\hfill (2.15)
\\[1.5ex]
is consistent in $G.$ Indeed, we obtain the following assertions.
\vspace{0.5ex}

{\bf Theorem 2.4.}\! 
{\it
For the partial differential system {\rm (2.1)} to have 
an $\!\!(n-k_1,n-k_2)\!$-cy\-lind\-ri\-ca\-li\-ty first integral of the form {\rm (2.14)} 
\vspace{0.5ex}
it is necessary that {\rm (2.6)} with {\rm (2.15)}  be consistent.
}

{\bf Corollary 2.3.}\!
{\it 
For the linear homogeneous system of partial differential equations {\rm (2.1)}
to have a holomorphic $(n-k_1,n)\!$-cy\-lin\-dricality 
first in\-teg\-ral of the form {\rm (2.14)} it is necessary that  the system of identities
{\rm (2.6)} with  {\rm (2.15)} and $k_2=0$ be consistent.
}
\vspace{0.75ex}

{\bf Corollary 2.4.}\!
{\it 
For the linear homogeneous system of partial differential equations {\rm (2.1)}
to have an an\-ti\-ho\-lo\-mor\-p\-hic $(n,n-k_2)\!$-cy\-lin\-dricality 
first in\-teg\-ral of the form {\rm (2.14)} it is necessary that the system of identities
{\rm (2.6)} with {\rm (2.15)} and $k_1=0$ be consistent.
}
\vspace{0.75ex}

The  proof of the following statements is similar to those of Theorems 2.2 and 2.3.
\vspace{0.75ex}

{\bf Theorem 2.5.}\!\!
{\it 
For the system {\rm (2.1)} to have at least one first integral of the form {\rm (2.14)}  
\vspace{0.25ex}
it is necessary and sufficient that there exist functions ${}^{\lambda}\!\varphi$ 
\vspace{0.25ex}
satisfying to system {\rm (2.7)} with  $\!H_{j}\!\equiv\! 0,$ $j\!=\!1,\ldots, m,\!$ 
that the function {\rm (2.14)} is a general integral of the Pfaffian equation {\rm (2.8)}.
}
\vspace{1ex}

{\bf Theorem 2.6.}\! 
\vspace{0.25ex}
{\it 
Let the system {\rm (2.7)} with  $H_{j}\equiv 0,\, j\!=\!1,\ldots, m$  
has $q$ not linearly bound so\-lu\-ti\-ons {\rm (2.12)} such that 
\vspace{0.35ex}
the corresponding Pfaff equations {\rm (2.13)} have the general integrals     
$\!F_{\varepsilon}\colon z\to F_{\varepsilon}({}^{k}\!z)\!$
for all $\!z\!\in\! G^{\,\prime},\, \varepsilon\!=\!1,\ldots, q.\!$
\vspace{1.75ex}
Then these integrals are functionally in\-de\-pen\-dent. 
}

As an example, consider the linear homogeneous system of partial differential equations 
\\[2ex]
\mbox{}\hfill                              
$
{\frak A}_1^{}(z_1^{},z_2^{})\;\!u=0,
\qquad  
{\frak A}_2^{}(z_1^{},z_2^{})\;\!u=0,
$
\hfill (2.16)
\\[2ex]
where the linear differential operators of first order
\\[2.25ex]
\mbox{}\hfill
$
{\frak A}_1^{}(z_1^{},z_2^{}) =
z_1^{}\overline{z}_2\,\partial_{z_1^{}} + 
(z_2^2+\overline{z}\;\!{}_1^2)\,\partial_{z_2^{}} +
(z_1-z_2^2+\overline{z}\;\!{}_1^2+\overline{z}\;\!{}_2^2)\,
\partial_{{}_{\scriptstyle \overline{z}_1}} -
z_1^2\,\partial_{{}_{\scriptstyle \overline{z}_2}} 
$
\ for all $(z_1^{},z_2^{})\in {\mathbb C}^2,
\hfill
$
\\[2.5ex]
\mbox{}\hfill
$
{\frak A}_2^{}(z_1^{},z_2^{}) =
\overline{z}\;\!{}_2^2\,\partial_{z_1^{}} + 
(z_1+z_2^2+\overline{z}\;\!{}_1^2+\overline{z}\;\!{}_2^2)\,\partial_{z_2^{}} +
(z_2^2+\overline{z}\;\!{}_1^{})\,\partial_{{}_{\scriptstyle \overline{z}_1}} -
z_1^{}\overline{z}_2\,\partial_{{}_{\scriptstyle \overline{z}_2}} 
$
\ for all $(z_1^{},z_2^{})\in {\mathbb C}^2.
\hfill
$
\\[2.25ex]
\indent
Let us find a (1,1)-cy\-lin\-dricality first integral of system (2.16) . 
\vspace{0.25ex}
The Wronskians of the sets of functions 
$U_1=\{z_1^{}\overline{z}_2,{}-z_1^2\}$ and 
$U_2=\{\overline{z}\;\!{}_2^2, {}-z_1^{}\overline{z}_2\}$
\vspace{0.5ex}
with respect to   $z_2^{}$ and $\overline{z}_1^{}$ vanish identically  on the space ${\mathbb C}^2.$
Therefore the necessary conditions for system (2.16) to have an 
$\!{\mathbb R}\!$-dif\-fe\-ren\-ti\-ab\-le (1,1)-cylindricality first integral is complied (Theorem 2.4).

Let us write the functional system {\rm (2.7)} with  $H_{1}\equiv 0, \, H_{2}\equiv 0\colon$ 
\\[2ex]
\mbox{}\hfill
$
z_1^{}\overline{z}_2\,\varphi_1^{} -z_1^2\,\varphi_2^{}=0,
\qquad
\overline{z}\;\!{}_2^2\,\varphi_1^{} -z_1^{}\overline{z}_2\,\varphi_2^{}=0,
\hfill
$
\\[2.25ex]
\mbox{}\hfill
$
\partial_{z_2^{}}^{}(z_1^{}\overline{z}_2)\,\varphi_1^{} +
\partial_{z_2^{}}^{}({}-z_1^2)\,\varphi_2^{}=0,
\qquad
\partial_{z_2^{}}^{}(\overline{z}\;\!{}_2^2)\,\varphi_1^{} +
\partial_{z_2^{}}^{}({}-z_1^{}\overline{z}_2)\,\varphi_2^{}=0,
\hfill
$
\\[2.25ex]
\mbox{}\hfill
$
\partial_{{}_{\scriptstyle \overline{z}_1}}^{}(z_1^{}\overline{z}_2)\,\varphi_1^{} +
\partial_{{}_{\scriptstyle \overline{z}_1}}^{}({}-z_1^2)\,\varphi_2^{}=0,
\qquad
\partial_{{}_{\scriptstyle \overline{z}_1}}^{}(\overline{z}\;\!{}_2^2)\,\varphi_1^{} +
\partial_{{}_{\scriptstyle \overline{z}_1}}^{}({}-z_1^{}\overline{z}_2)\,\varphi_2^{}=0.
\hfill
$
\\[2.25ex]
\indent
This system is reduced to the equation 
$
\overline{z}_2\,\varphi_1^{} -z_1^{}\,\varphi_2^{}=0.
$
The scalar functions 
\\[1.5ex]
\mbox{}\hfill
$
\varphi_1^{}\colon (z_1^{},z_2^{})\to z_1^{}
$
\ for all $(z_1^{},z_2^{})\in {\mathbb C}^2,
\quad
\varphi_2^{}\colon (z_1^{},z_2^{})\to \overline{z}_2^{}
$ 
\ for all $(z_1^{},z_2^{})\in {\mathbb C}^2
\hfill
$
\\[2ex]
is a solution to this equation.
The corresponding Pfaffian differential equation
\\[1.75ex]
\mbox{}\hfill
$
z_1^{}\,dz_1^{} + \overline{z}_2^{}\,d\,\overline{z}_2^{} =0
\hfill
$
\\[1.75ex]
has the general integral (Theorem 2.5)
\\[1.75ex]
\mbox{}\hfill                                       
$
F\colon (z_1^{},z_2^{}) \to \,
z_1^{2}+\overline{z}\;\!{}_2^2
$
\ for all $(z_1^{},z_2^{})\in {\mathbb C}^2.
$
\hfill (2.17)
\\[1.75ex]
\indent
Since the Poisson bracket
\\[2.25ex]
\mbox{}\hfill
$
\bigl[{\frak A}_1^{}(z_1^{},z_2^{}), {\frak A}_2^{}(z_1^{},z_2^{})\bigr] =
{}-\overline{z}_2^{}(z_1^2+\overline{z}\;\!{}_2^2)\,\partial_{z_1^{}} +
({}-2z_1^{} z_2^{} +
2z_1^{}\overline{z}_1^{}+
z_1^{}\overline{z}_2^{} 
+ 2\,\overline{z}\;\!{}_1^{2}
- 2\;\!z_1^{2}\,\overline{z}_2^{}\ -
\hfill
$
\\[2.5ex]
\mbox{}\hfill
$
-\ 2\;\!z_2^{}\;\!\overline{z}\;\!{}_2^{2} 
-4z_2^{2}\,\overline{z}_1^{}
+2\,\overline{z}_1^{}\overline{z}\;\!{}_2^{2}
+2\,\overline{z}\;\!{}_1^{3}\;\!)\,\partial_{z_2^{}} 
+
(z_1^{} + 2z_1^{}z_2^{}
-z_2^{2}
-\overline{z}\;\!{}_1^{2}
+2\;\!z_1^{}\;\!\overline{z}\;\!{}_2^{2}
+4\;\!z_2^{3}
\ +
\hfill
$
\\[2.5ex]
\mbox{}\hfill
$
+\ 4z_2^{}\;\!\overline{z}\;\!{}_1^{2} 
+ 2\;\!z_2^{}\;\!\overline{z}\;\!{}_2^{2} 
-2\;\!z_2^{2}\,\overline{z}_1^{})\,\partial_{{}_{\scriptstyle \overline{z}_1}} 
+ z_1^{}(z_1^2+\overline{z}\;\!{}_2^2)\,\partial_{{}_{\scriptstyle \overline{z}_2}}
$
\ for all $(z_1^{},z_2^{})\in {\mathbb C}^2,
\hfill
$
\\[2.25ex]
is not a linear combination on the space ${\mathbb C}^2$ of the operators ${\frak A}_1^{}$
and ${\frak A}_2^{},$ we see that the 
linear homogeneous partial differential system (2.16) is not complete.
\vspace{0.35ex}

Thus the ${\mathbb R}\!$-differentiable (1,1)-cylindricality 
first integral (2.17) is an integral basis on the  space ${\mathbb C}^2$ of the 
incomplete system of partial differential equations (2.16).
\vspace{1.25ex}

{\bf 2.1.3. {\boldmath $\!(n-k_1,n-k_2)\!$}-cylindricality last multipliers}.
\vspace{0.35ex}
Suppose the system (2.1) has an $(n-k_1,n-k_2)\!$-cylindricality 
$\!{\mathbb R}\!$-dif\-fe\-ren\-ti\-ab\-le on the domain $G^{\,\prime}$ last mul\-tip\-lier 
\\[1.5ex]
\mbox{}\hfill                              
$
\mu\colon z\to\  \mu ({}^{k}\!z)$
\ for all $z\in G^{\,\prime}.
$
\hfill (2.18)
\\[1.5ex]
\indent
Then, in accordance with the criteria of an $\!{\mathbb R}\!$-dif\-fe\-ren\-ti\-ab\-le 
last multiplier,
\\[1.5ex]
\mbox{}\hfill                            
$
{\frak A}_{j}^{k}\,\mu ({}^{k}\!z)
\, +\,  \mu ({}^{k}\!z)\,{\rm div}\,u^{j}(z) = 0$
\ for all $z\in G^{\,\prime},
\quad 
j=1,\ldots, m.
$
\hfill (2.19)
\\[1.75ex]
\indent
Using the methods of Subsubsection 2.1.1, we get the following statements.
\vspace{0.75ex}

{\bf Theorem 2.7.}\!
{\it 
For the partial differential system {\rm (2.1)} to have an 
$\!{\mathbb R}\!\!$-dif\-fe\-ren\-ti\-ab\-le  
last mul\-ti\-p\-li\-er of the form {\rm (2.18)} it is necessary that 
the system of identities}
\\[2ex]
\mbox{}\hfill                                       
$
W_{{}_{\scriptstyle z_{\gamma}}}
\bigl({}^{\lambda}u^{j}(z), {\rm div} \,u^{j}(z)\bigr)=0$
\ for all $z\in G,
\ \ j=1,\ldots, m,
\ \gamma=k_1+1,\ldots, n,
\hfill
$
\\
\mbox{}\hfill (2.20)
\\
\mbox{}\hfill
$
W_{{}_{\scriptstyle \overline{z}_{\zeta_{\delta}}}}
\bigl({}^{\lambda}u^{j}(z), {\rm div} \,u^{j}(z)\bigr)=0$
\ for all $z\in G,
\ \ j=1,\ldots, m,
\ \delta=k_2+1,\ldots, n,
\hfill
$
\\[1.75ex]
{\it be consistent on the domain $G.$
}
\vspace{0.5ex}

{\bf Corollary 2.5.}\!
{\it 
For the system {\rm (2.1)} to have a  holomorphic $\!(n-k_1,n)\!$-cy\-lin\-dricality
last mul\-ti\-p\-li\-er of the form {\rm (2.18)} it is necessary that {\rm (2.20)} 
with $k_2=0$ be consistent.
}
\vspace{0.5ex}

{\bf Corollary 2.6.}\!
{\it 
For the system {\rm (2.1)} to have an antiholomorphic 
$(n,n-k_2)\!$-cy\-lin\-dricality last multiplier of the form 
{\rm (2.18)} it is necessary that {\rm (2.20)} with $k_1=0$ 
\vspace{0.5ex}
be consistent.
}

{\bf Theorem 2.8.}\! 
{\it
For the system {\rm (2.1)} to have at least one last multiplier 
of the form {\rm (2.18)} it is necessary and sufficient 
\vspace{0.35ex}
that there exist function ${}^{\lambda}\!\varphi$ 
satisfying sys\-tem {\rm (2.7)} with 
\\[1.25ex]
\mbox{}\hfill                                      
$
H_j\colon z\to {}-{\rm div}\, u^{j}(z)$
\ for all $z\in G,\ \ j=1,\ldots, m,
$ 
\hfill {\rm (2.21)}
\\[1.5ex]
such that the Pfaffian equation {\rm (2.8)} has the integrating factor 
\vspace{0.35ex}
$\nu\colon {}^{k}\!z\to 1$ for all $z\in G^{\,\prime};$  
in this case the last multiplier is given by
\\[1.5ex]
\mbox{}\hfill                                        
$
\mu\colon z\to\ \exp g({}^{k}\!z)$ 
\ for all $z\in G^{\,\prime},
$
\hfill {\rm (2.22)}
\\[2ex]
where the function
\vspace{1.25ex}
$g\colon z\to \int {}^{k_1}\!\varphi({}^{k}\!z)\, d{}^{k_{1}}\!z\, +\,
{}^{k_2}\!\varphi({}^{k}\!z)\,d\,\overline{{}^{k_{2}}\!z}$ 
for all $z\in G^{\,\prime}.$
}

{\sl Proof. Necessity}. 
\vspace{0.25ex}
Let the function (2.18) be an $\!(n-k_1,n-k_2)\!$-cylindricality 
$\!{\mathbb R}\!$-dif\-fe\-ren\-ti\-ab\-le
last multiplier of the partial differential system (2.1).
Then the vector functions 
\\[1.75ex]
\mbox{}\hfill
$
{}^{k_1}\varphi\colon z\to\,
\partial_{{}_{\scriptstyle {}^{k_1}\!z}}
\ln \mu ({}^{k}\!z)$
\  for all $z\in G^{\,\prime},
\qquad
{}^{k_2}\varphi\colon z\to\,
\partial_{{}_{\scriptstyle \overline{{}^{k_2}\!z}}}\,
\ln \mu ({}^{k}\!z)$
\  for all $z\in G^{\,\prime}
\hfill
$
\\[1.75ex]
are a solution to system (2.7).
This implies that the function  
$\nu\colon {}^{k}\!z\to 1$ for all $z\in G^{\,\prime}$
is an integrating factor of the Pfaffian differential equation (2.8). 
\vspace{0.5ex}

{\sl Sufficiency}. 
Let ${}^{\lambda}\!\varphi$ be a solution to system (2.7) with (2.21)
and let $\nu\colon {}^{k}\!z\to 1$ be an integrating factor 
of the corresponding Pfaffian differential equation  (2.8). Then 
\\[1.75ex]
\mbox{}\hfill                      
$
\partial_{{}_{\scriptstyle {}^{k_1}\!z}}\,g({}^{k}\!z)\, -\,{}^{k_1}\varphi({}^{k}\!z)=0,
\qquad
\partial_{{}_{\scriptstyle \overline{{}^{k_2}\!z}}}\,\, g({}^{k}\!z) \, -\,{}^{k_2}\varphi({}^{k}\!z)=0.
\hfill
$
\\[1.75ex]
\indent
Using {\rm (2.7)} with (2.21), we have the identity (2.19) is valid.
This yields that the scalar function (2.22) is a last multiplier of the partial differential system (2.1).\, \k
\vspace{1.25ex}

For example, consider the linear homogeneous system of partial differential equations 
\\[1.75ex]
\mbox{}\hfill                              
$
{\frak A}_1^{}(z_1^{},z_2^{})\;\!u=0,
\qquad  
{\frak A}_2^{}(z_1^{},z_2^{})\;\!u=0,
$
\hfill (2.23)
\\[2ex]
where the linear differential operators 
\vspace{0.75ex}
$
{\frak A}_1^{}(z_1^{},z_2^{}) =
z_2^{}\;\!\overline{z}{}_2^{}\,\partial_{z_1^{}} + 
\overline{z}\;\!{}_1^2\,\partial_{z_2^{}} +
z_1^{}\;\!\overline{z}{}_2^{}\,\partial_{{}_{\scriptstyle \overline{z}_1}}\! +
\overline{z}{}_1^{}\overline{z}{}_2^{}\,\partial_{{}_{\scriptstyle \overline{z}_2}} 
$
for all $(z_1^{},z_2^{})\in {\mathbb C}^2,
$
$
{\frak A}_2^{}(z_1^{},z_2^{}) =
\overline{z}\;\!{}_2^{2}\,\partial_{z_1^{}} + 
\overline{z}\;\!{}_1^2\,\partial_{z_2^{}} +
z_2^{}\;\!\overline{z}{}_2^{}\,\partial_{{}_{\scriptstyle \overline{z}_1}}\! +
z_1^{}\;\!\overline{z}{}_2^{}\,\partial_{{}_{\scriptstyle \overline{z}_2}} 
$
\vspace{1ex}
for all $(z_1^{},z_2^{})\in {\mathbb C}^2.$

Let us find a (2,1)-cy\-lin\-dricality antiholomorphic last multiplier of system (2.23). 
\vspace{0.75ex}

The divergences 
\vspace{0.75ex}
$
{\rm div}\;\! u^1_{}(z_1^{},z_2^{}) =
\partial_{z_1^{}}(z_2^{}\;\!\overline{z}{}_2^{}) + 
\partial_{z_2^{}}(\overline{z}\;\!{}_1^2) +
\partial_{{}_{\scriptstyle \overline{z}_1}} (z_1^{}\;\!\overline{z}{}_2^{}) +
\partial_{{}_{\scriptstyle \overline{z}_2}} (\overline{z}{}_1^{}\overline{z}{}_2^{})=
\overline{z}{}_1^{}
$
and
$
{\rm div}\;\! u^2_{}(z_1^{},z_2^{}) =
\partial_{z_1^{}}(\overline{z}\;\!{}_2^{2}) + 
\partial_{z_2^{}}(\overline{z}\;\!{}_1^2) +
\partial_{{}_{\scriptstyle \overline{z}_1}} (z_2^{}\;\!\overline{z}{}_2^{}) +
\partial_{{}_{\scriptstyle \overline{z}_2}} (z_1^{}\;\!\overline{z}{}_2^{}) =z_1^{}
$
for all $(z_1^{},z_2^{})\in {\mathbb C}^2.$
\vspace{1.25ex}

The Wronskians of the sets of functions 
$U_1=\{\overline{z}_1^{}\overline{z}_2, \overline{z}_1^{}\}$ and 
$U_2=\{z_1^{}\overline{z}_2^{}, z_1^{}\}$
\vspace{0.65ex}
with respect to  $z_1^{},\, z_2^{},$ and $\overline{z}_1^{}$ vanish identically  on the space ${\mathbb C}^2.$
\vspace{0.5ex}

Therefore the necessary conditions for system (2.23) to have a 
(2,1)-cy\-lin\-dricality antiho\-lo\-mor\-p\-hic last multiplier is complied (Corollary 2.6).
Let us write the system {\rm (2.7)} with  (2.21):
\\[2ex]
\mbox{}\hfill
$
\overline{z}_1^{}\overline{z}_2\,\varphi_1^{}={}-\overline{z}_1^{},
\quad \
z_1^{}\overline{z}_2^{}\,\varphi_1^{} ={}-z_1^{},
\quad \
\overline{z}_2^{}\,\varphi_1^{} ={}-1.
$
\hfill (2.24)
\\[2ex]
The function 
\vspace{0.5ex}
$
\varphi_1^{}\colon\! (z_1^{},z_2^{})\to {}-\dfrac{1}{\overline{z}_2^{}}
$ 
for all $(z_1^{},z_2^{})\in G^{\,\prime},$ 
where a domain 
$G^{\,\prime}\!\subset\! \{(z_1^{},z_2^{})\colon z_2^{}\!\ne\! 0\},\!\!$
is a solution to the sys\-tem (2.24).
By Theorem 2.8, the function
\\[1.75ex]
\mbox{}\hfill
$
\mu\colon (z_1^{},z_2^{})\to\, \dfrac{1}{\overline{z}_2^{}}
$
\ for all $(z_1^{},z_2^{})\in G^{\,\prime}
\hfill
$
\\[1.75ex]
is a (2,1)-cy\-lin\-dricality antiho\-lo\-mor\-p\-hic last multiplier 
\vspace{1.25ex}
on the domain $G^{\,\prime}$ of system (2.23).

{\bf Theorem 2.9.} 
\vspace{0.25ex}
{\it Let the  system {\rm (2.7)} with {\rm (2.21)} has $q$ not linearly bound solutions 
{\rm (2.12)} for which the corresponding 
\vspace{0.25ex}
Pfaffian differential equations {\rm (2.13)} have the integrating factors  
$\nu_{\varepsilon}({}^{k}\!z)=1$
for all $z\in G^{\,\prime},\ \varepsilon=1,\ldots, q.$   
Then the last multiplies 
\\[2ex]
\mbox{}\hfill                       
$
\displaystyle
\mu_{\varepsilon}\colon z\to
\exp \int {}^{k_1}\!\varphi^{\varepsilon}({}^{k}\!z)\,d{}^{k_{1}}\!z \,+\,
{}^{k_2}\!\varphi^{\varepsilon}({}^{k}\!z)\,d\,\overline{{}^{k_{2}}\!z}$
\ for all $z\in G^{\,\prime},
\ \ \varepsilon=1,\ldots, q
\hfill
$
\\[1.5ex]
of system {\rm (2.1)} are functionally independent.
}
\vspace{0.75ex}

The idea of the proof of Theorem 2.9 is similar to that one in Theorem 2.3.

\newpage

{\bf  
2.2. First integrals of  linear homogeneous system with $\!{\mathbb R}\!$-linear coefficients
}
\\[1ex]
\indent
Let us consider a linear homogeneous system of first-order partial differential equations
\\[1.75ex]
\mbox{}\hfill                           
$
{\frak L}_{j} (z) \, w \ = \ 0, 
\quad j=1,\ldots,m,
$
\hfill (2.25)
\\[1.75ex]
where the coefficients of the linear differential operators 
\\[1.75ex]
\mbox{}\hfill                          
$
\displaystyle
{\frak L}_{j} (z) = \sum\limits_{\xi=1}^{n}\,
\bigl(a_{{}_{\scriptstyle j\xi}}(z)\,\partial_{{}_{\scriptstyle z_{\xi}}}+
a_{{}_{\scriptstyle j, n+\xi}}(z)\,\partial_{{}_{\scriptstyle \overline{z}_{\xi}}}\bigr)$ 
\  for all 
$z\in {\mathbb C}^{n},
\quad j=1,\ldots,m,
\hfill                      
$
\\[1.75ex]
are the ${\mathbb R}\!$-linear [2, p. 21] functions
\\[1.75ex]
\mbox{}\hfill                          
$
\displaystyle
a_{jk}\colon z \to \sum\limits_{\tau=1}^{n}
\bigl(a_{jk\tau}\,z_{\tau}+a_{jk,n+\tau}\,\overline{z}_{\tau}\bigr)$ \, 
for all 
$z\in  {\mathbb C}^{n}
\ \ (a_{{}_{\scriptstyle jk l }}\in {\mathbb C}, \ \,
l, k\!=\!1,\ldots, 2n,\,  j\!=\!1,\ldots, m).
\hfill                      
$
\\[1.75ex]
\indent
Assume that the system (2.25) is related by the conditions in terms of the Poisson brackets 
\\[1.75ex]
\mbox{}\hfill                                       
$
\bigl[ {\frak L}_{j}(z), {\frak L}_{\zeta}(z)\bigr]  = {\frak O}
\ \ \ \text{for all}\  z\in {\mathbb C}^n,
\ \ j=1,\ldots, m,\ \zeta=1,\ldots, m,
$
\hfill (2.26)
\\[1.75ex]
where ${\frak O}$ is the null operator,
\vspace{0.35ex}
i.e., the system (2.25) is {\it jacobian} [17, p. 523; 19, pp. 38 -- 40].

An integral basis of the jacobian system (2.25)
is $2n-m$ (the proof  is similar to that one in [28, p. 70])
functionally in\-de\-pen\-dent $\!{\mathbb R}\!$-differentiable first integrals of system (2.25).

In this Subsection we study Darboux's problem of finding first in\-teg\-rals 
in case that partial integrals are known [29]. Using method of partial integrals 
[18; 19, pp. 161 -- 311; 30 -- 33], 
we obtain the spectral method [9] for building first integrals of the jacobian system (2.25).
\\[1.5ex]
\indent                       
{\bf 2.2.1. $\!\!{\mathbb R}\!$-linear partial integral}.
The $\R\!$-linear function
\\[2ex]
\mbox{}\hfill
$
\displaystyle
p\colon z\to 
\sum\limits_{\xi=1}^{n}
\bigl(
b_{{}_{\scriptstyle \xi}} z_{{}_{\scriptstyle \xi}} +
b_{{}_{\scriptstyle n+\xi}}
\overline{z}_{{}_{\scriptstyle \xi}}\bigr)
\ \ \text{for all}\  \,z\in {\mathbb C}^n
\quad 
(b_{l}\in {\mathbb C},\ l=1,\ldots, 2n)
\hfill
$
\\[2ex]
is a {\it partial integral} of the system (2.25) if and only if
\\[1.75ex]
\mbox{}\hfill
$
{\frak L}_{{}_{\scriptstyle j}}\, p(z)= p(z)\lambda^j
\ \ \text{for all}\  z\in {\mathbb C}^n,
\quad
\lambda^j\in {\mathbb C}, \ j=1,\ldots, m.
\hfill
$
\\[1.75ex]
This system of identities is equivalent to the linear homogeneous system 
\\[1.75ex]
\mbox{}\hfill                                      
$
\bigl(A_{j} - \lambda^{j} E\bigr)\,b= 0,
\ \ \ j=1,\ldots, m,
$
\hfill (2.27)
\\[1.75ex]
where $A_{j}= \|a_{jk l}\|,\,  j=1,\ldots, m$ 
\vspace{0.75ex}
are $2n\times 2n\!$-matrices,    
$E$ is the $2n\times 2n$ identity matrix, 
$b=\mbox{colon}(b_1,\ldots,b_{2n})$ is a vector column. 
\vspace{0.75ex}

The conditions (2.26) for the partial differential system (2.25) are equivalent
\\[1.5ex]
\mbox{}\hfill
$
A_{j}A_{\zeta} = A_{\zeta}A_{j},
\ \ \ 
j=1,\ldots,m,
\ \ \zeta=1,\ldots, m.
\hfill
$
\\[1.75ex]
\indent
Then there exists 
\vspace{0.25ex}
a relation [34, pp. 193 -- 194; 35] between eigenvectors and 
eigenvalues of the matrices $A_{j},\ j=1,\ldots, m.$
\vspace{0.75ex}

{\bf  Lemma 2.1.}
\vspace{0.25ex}
{\it
Let $\nu\in {\mathbb C}^{2n}$ be a common eigenvector of the matrices 
$A_{j},\ j=1,\ldots, m.$ Then the  function
$p\colon z\to \nu\gamma$ for all \ $z\in {\mathbb C}^n,$
\vspace{0.25ex}
where $\gamma=\mbox{\rm colon}(z_1,\ldots,z_n,\overline{z}_{1},\ldots,\overline{z}_n),$
is an  ${\mathbb R}\!$-linear partial integral of 
the system of partial differential equations {\rm (2.25)}. 
}
\vspace{0.5ex}

{\sl Proof}. If $\nu$ is a common eigenvector of the matrices 
$A_{j},\ j=1,\ldots, m,$ then $\nu$ is a solution to system (2.27),
where $\lambda^{j}$ is an eigenvalue of the matrix $A_{j}$ 
corresponding to the eigenvector $\nu.$ Therefore we obtain the identities
\\[1.5ex]
\mbox{}\hfill
$
{\frak L}_{j}\, \nu\gamma = \lambda^{j}\,\nu\gamma$
\ for all $z\in {\mathbb C}^n,
\ \ \ j=1,\ldots,m.
\hfill
$
\\[1.5ex]
Thus  the $\!{\mathbb R}\!$-li\-ne\-ar function $\!p\colon z\to \nu\gamma$ for all $\!z\!\in\! {\mathbb C}^n\!$
\vspace{2.25ex}
is a partial integral of system (2.25). \k

{\bf 2.2.2. $\!\!{\mathbb R}\!$-differentiable first integrals}
\\[1ex]
\indent    
{\bf Theorem 2.10.}
\vspace{0.35ex}
{\it 
Let $\nu^{\theta},\, \theta=1,\ldots, m+1,$ be common eigenvectors of  the matrices  
$A_{j},$ $j=1,\ldots, m.$
Then the system {\rm (2.25)} has the ${\mathbb R}\!$-differentiable first integral
\\[2ex]
\mbox{}\hfill                                         
$
\displaystyle
F\colon z\to
\prod\limits_{\theta=1}^{m+1}
\bigl(\nu^{\theta}\gamma\bigr)^{h_{\theta}}$ \ 
for all $z\in\Omega,
\ \ \ \Omega\subset {\rm D}(F),
$
\hfill {\rm (2.28)}
\\[1.5ex]
where $h_{1},\ldots,h_{m+1}$ is a nontrivial solution to the system  
$
\sum\limits_{\theta=1}^{m+1}\,\lambda_{\theta}^{j}\,h_{\theta} =0,
\ j=1,\ldots,m,
$
and $\lambda_{\theta}^{j}$ are the eigenvalues of the matrices $A_{j},\ j=1,\ldots, m,$ 
\vspace{0.35ex}
corresponding to the common eigenvectors $\nu^{\theta},\ \theta=1,\ldots, m+1,$ respectively.
}
\vspace{0.75ex}

{\sl Proof}.
Suppose $\nu^{\theta}$ are common eigenvectors of  the matrices  $A_{j}$ 
\vspace{0.35ex}
corresponding to the eigenvalues $\lambda_{\theta}^{j},\ j=1,\ldots, m,\ \theta=1,\ldots,m+1,$ respectively.
\vspace{0.5ex}

By Lemma 2.1, it follows that the ${\mathbb R}\!$-linear functions 
\\[1.5ex]
\mbox{}\hfill
$
p_{\theta}\colon z\to\nu^{\theta}\gamma$ 
\ for all $z\in {\mathbb C}^n,
\quad \theta=1,\ldots, m+1,
\hfill
$
\\[1.5ex]
are partial integrals of the system of partial differential equations (2.25). Hence,
\\[2ex]
\mbox{}\hfill                                         
$
{\frak L}_{{}_{\scriptstyle j}}\,\nu^{\theta}\gamma =
\lambda_{\theta}^{j}\,\nu^{\theta}\gamma
$
\ for all $z\in {\mathbb C}^n,
\ \
j=1,\ldots, m,
\quad \theta=1,\ldots, m+1.
$
\hfill (2.29)
\\[1.75ex]
\indent
We form the function
\\[1.5ex]
\mbox{}\hfill
$
\displaystyle
F\colon z\to
\prod\limits_{\theta=1}^{m+1}
\bigl(\nu^{\theta}\gamma\bigr)^{h_{\theta}}
$ 
\ \ for all \ $ z\in\Omega,
\hfill
$  
\\[1.5ex]
where $\Omega$ is a domain (open arcwise connected set) in ${\mathbb C}^n$ 
\vspace{0.35ex}
and $h_{\theta},\,\theta=1,\ldots, m+1,$ are complex numbers 
with $\sum\limits_{\theta=1}^{m+1}|h_{\theta}|\ne 0.$
The Lie derivative of $F$ by virtue of (2.25) is equal to 
\\[2ex]
\mbox{}\hfill
$
\displaystyle
{\frak L}_{{}_{\scriptstyle j}} F(z)=
\prod\limits_{\theta=1}^{m+1}
\bigl(\nu^{\theta}\gamma\bigr)^{h_{\theta}-1}\,
\sum\limits_{\theta=1}^{m+1} h_{\theta}\,
\prod\limits_{l=1,l\ne \theta}^{m+1}(\nu^{l}\gamma) \
{\frak L}_{{}_{\scriptstyle j}}\,\nu^{\theta}\gamma
$
\ \ for all $z\in\Omega, 
\quad j=1,\ldots, m.
\hfill
$
\\[1.75ex]
\indent
Using (2.29), we have
\\[1.5ex]
\mbox{}\hfill
$
\displaystyle
{\frak L}_{{}_{\scriptstyle j}}F(z)=
\sum\limits_{\theta=1}^{m+1}\lambda_{\theta}^{j}\,h_{\theta}\,F(z)
$
\ \ for all $z\in\Omega, 
\quad j=1,\ldots, m.
\hfill
$
\\[1.5ex]
\indent
If $\sum\limits_{\theta=1}^{m+1}\lambda_{\theta}^{j}\,h_{\theta}=0,\
j=1,\ldots m,$ then the function (2.28) is an ${\mathbb R}\!$-differentiable  first integral 
of  the linear homogeneous system of partial differential equations (2.25). \k
\vspace{1.25ex}

{\bf Corollary 2.7.}
\vspace{0.35ex}
{\it
Let $\nu^{\theta}$ be common eigenvectors of  the matrices  $A_{j}$ 
corresponding to the eigenvalues $\lambda_{\theta}^{j},\ j=1,\ldots, m,\ \theta=1,\ldots,m+1,$ respectively.
\vspace{0.35ex}
Then the linear homogeneous sys\-tem of partial differential equations
 {\rm (2.25)} has the ${\mathbb R}\!$-differentiable first integral
\\[1.75ex]
\mbox{}\hfill
$
F_{12\ldots m(m+1)}\colon z\to
\prod\limits_{\theta=1}^{m}
\bigl(\nu^{\theta}\gamma\bigr)^{{}-\delta_{\theta}}
\bigl(\nu^{m+1}\gamma\bigr)^{\delta}$
\ \ for all $z\in\Omega,
\quad \Omega\subset {\rm D}(F_{12\ldots m(m+1)}),
\hfill
$
\\[1.75ex]
where the determinants $\delta_{\theta},\, \theta=1,\ldots, m$ 
\vspace{0.35ex}
are obtained by replacing the $\theta\!$-th column of the determinant 
$\delta=\big|\lambda_{\theta}^{j}\big|$  by  
${\rm colon}\bigl(\lambda_{m+1}^{1}, \ldots,\lambda_{m+1}^{m}\bigr),$ respectively.
}
\vspace{1.5ex}

For example, the linear homogeneous system of first-order partial differential equations
\\[2ex]
\mbox{}\hfill                              
$
{}-z_1\,\partial_{{}_{\scriptstyle z_1}}w+z_2\,\partial_{{}_{\scriptstyle z_2}}w+
\overline{z}_2\,\partial_{{}_{\scriptstyle \overline{z}_1}}w+
\overline{z}_1\,\partial_{{}_{\scriptstyle \overline{z}_2}}w=0,
\ \ 
2(\overline{z}_1+\overline{z}_2)\,\partial_{{}_{\scriptstyle z_2}}w +
z_2\,\partial_{{}_{\scriptstyle \overline{z}_1}}w+
z_2\,\partial_{{}_{\scriptstyle \overline{z}_2}}w=0
$
\hfill (2.30)
\\[2ex]
has the commuting matrices 
\\[1.5ex]
\mbox{}\hfill
$
A_1 =
\left\|\!
\begin{array}{rrcc}
   -1 & 0 & 0 & 0
\\
  0 &   1 & 0 & 0
\\
  0 &    0 &  0 & 1
\\
0 & 0 & 1 & 0
\end{array}\!
\right\|
$
\ \ and \ \ 
$
A_2 =
\left\|\!
\begin{array}{crcc}
0 &  0 &  0 & 0
\\
 0 & 0 & 1 & 1
\\
  0 &  2 &  0 & 0
\\
  0 &  2 &   0 & 0
\end{array}\!
\right\|.
\hfill
$
\\[2ex]
\indent
Therefore the system of partial differential equations (2.30) is jacobian. 
\vspace{0.35ex}

The matrices $\!A_1\!$ and $\!A_2\!$ have the eigenvalues 
\vspace{0.5ex}
$\lambda_{1}^{1}\!=1,\,
\lambda_{2}^1=\lambda_{3}^{1}=-1,\,
\lambda_{4}^1\!=1,\!$
and
$\!\lambda_{1}^{2}\!=-2,\!\!$
$\lambda_{2}^2=\lambda_{3}^{2}=0,
\lambda_{4}^2=2
$
\vspace{0.5ex}
corresponding to the eigenvectors 
$\nu^1\!=\!(0,-1,1,1),\  \nu^2\!=\!(1,0,0,0),\!$
$\nu^3=(0,0,1,{}-1),\ \nu^4\!=\!(0,1,1,1),$ respectively. 
\vspace{0.5ex}

The solution to the  system 
$
h_{11}-h_{12}-h_{13}=0,
\
{}-2h_{11}=0
$
\vspace{0.5ex}
is $h_{11}=0,\ h_{12}=-1, \ h_{13}=1.$

The solution to the linear homogeneous system 
$
h_{21}-h_{22}+h_{24}=0,
\ {}-2h_{21}+2h_{24}=0$
is $h_{21}=1,\ h_{22}=2, \ h_{24}=1.$
\vspace{0.5ex}

The ${\mathbb R}\!$-differentiable functions (by Theorem 2.10)
\\[1.75ex]
\mbox{}\hfill                                    
$
F_1\colon (z_1,z_2)\to
\dfrac{\overline{z}_1-\overline{z}_2}{z_1}$
\ \ for all 
$(z_1,z_2)\in \Omega
$
\hfill (2.31)
\\[1ex]
and                                
\\[1.25ex]
\mbox{}\hfill                                   
$
F_2\colon (z_1,z_2)\to
z_1^2\bigl(z_2^2-(\overline{z}_1+\overline{z}_2)^2\bigr)$
\ for all 
$(z_1,z_2)\in {\mathbb C}^2,
$
\hfill (2.32)
\\[2ex]
where a domain $\Omega\subset \{(z_1,z_2)\colon z_1\ne 0\},$ are 
first integrals of the system (2.30).
\vspace{0.5ex}

The ${\mathbb R}\!$-differentiable first integrals (2.31) and (2.32)
are an integral basis of the  jacobian linear homogeneous system of first-order partial differential equations (2.30). 
\vspace{1ex}

From the entire set of partial differential equations (2.25), we extract the equation
\\[1.75ex]
\mbox{}\hfill                                                           
$
{\frak L}_{\zeta}(z)\,w=0,
\quad
\zeta\in \{1,\ldots,m\},
\hfill (2.25.\zeta)
$
\\[1.75ex]
such that the matrix $A_{\zeta}$ has the smallest number of elementary divisors [34, p. 147]. 
\vspace{0.75ex}

{\bf Definition 2.1.}
{\it
Let $\nu^{0l}$ be an eigenvector of the matrix $A_{\zeta}$ corresponding to the eigen\-va\-lue $\lambda_{l}^{\zeta}$ with 
elementary divisor of multiplicity $s_{l}.$ 
A non-zero vector $\nu^{\eta l}\in {\mathbb C}^{2n}$ is called a
\textit{\textbf{generalized eigenvector of order}} {\boldmath $\eta$} for $\lambda_{l}^{\zeta}$ if and only if
\\[1.5ex]
\mbox{}\hfill                          
$
(A_{\zeta}-\lambda_{l}^{\zeta} E)\,\nu^{\eta l}=\eta\, \nu^{\eta-1, l},
\quad 
\eta=1,\ldots, s_{l}-1,
$
\hfill {\rm (2.33)}
\\[1.5ex]
where $E$ is the $2n\times 2n$ identity matrix.
}
\vspace{1ex}

Using Lemma 2.1 and (2.33), we obtain
\\[1.5ex]
\mbox{}\hfill                                
$
{\frak L}_{\zeta}\,\nu^{0l}\gamma =
\lambda_{l}^{\zeta}\,\nu^{0l}\gamma
$
\ for all $z\in {\mathbb C}^n,
\hfill
$
\\[-0.15ex]
\mbox{}\hfill (2.34)
\\[-0.15ex]
\mbox{}\hfill 
$
{\frak L}_{\zeta}\,\nu^{\eta l}\gamma =
\lambda_{l}^{\zeta}\,\nu^{\eta l}\gamma +
\eta\,\nu^{\eta-1{,}\,l}\gamma$
\ for all $z\in {\mathbb C}^n,
\ \ \eta=1,\ldots, s_{l}-1.
\hfill
$
\\[1.75ex]
\indent
The following lemma is needed for the sequel.
\vspace{1ex}

{\bf Lemma 2.2.}                                           
{\it
Let $\nu^{0l}$ be a common eigenvector of the matrices $A_j$ corresponding to the eigenvalues
$\lambda^j_{l}, \ j=1,\ldots,m,$ respectively. Let $\nu^{\,\eta l},\ \eta =1,\ldots, s_l-1$ be generalized eigenvectors
of the matrix $A_\zeta$ 
\vspace{0.35ex}
corresponding to the eigenvalue $\lambda^{\zeta}_{l}$ with elementary divisor of multiplicity 
$s_l\ (s_{l}\geq 2).$ If the partial diffrential equation $(2.25.\zeta)$ hasn't the first integrals 
\\[2ex]
\mbox{}\hfill                               
$
\displaystyle
F_{j\eta l}^{\,\zeta}\colon z\to
{\frak A}_j\, \Psi_{\eta l}^{\zeta}(z)$ 
\ for all $z\in \Omega, 
\ \ j=1,\ldots, m, \ \, j\ne\zeta,
\ \ \  \eta =1,\ldots, s_l-1,
$
\hfill {\rm (2.35)}
\\[2ex]
then 
\\[1.75ex]
\mbox{}\hfill                             
$
{\frak L}_{\zeta}\, \Psi_{\eta l}^{\zeta}(z)=
\left[\!\!
\begin{array}{lll}
1\! & \text{for all}\ \ z\in \Omega, & \eta =1,
\\[1ex]
0\! & \text{for all}\ \ z\in \Omega, & \eta =2,\ldots, s_{l}-1,
\end{array}
\right.
$
\hfill {\rm (2.36)}
\\[2ex]
\mbox{}\hfill                               
$
{\frak L}_{j}\, \Psi_{\eta l}^{\zeta}(z)=\mu_{\eta l}^{j\zeta}={\rm const}$
\ for all $z\in \Omega,
\ \ j=1,\ldots, m,\ \ j\ne \zeta,
\ \ \  \eta =1,\ldots, s_l-1,
\hfill
$ 
\\[2.5ex]
where $\Psi_{\eta l}^{\zeta}\colon \Omega\to {\mathbb C},\  \eta =1,\ldots, s_l-1,$ is a solution to the system 
\\[2ex]
\mbox{}\hfill                             
$
\nu^{\,\eta l}\gamma=
{\displaystyle \sum\limits_{\delta=1}^{\eta} }
\binom{\eta -1}{\delta-1}\Psi_{\delta l}^{\zeta}(z)\, \nu^{\,\eta-\delta,l}\gamma,
\quad 
\eta=1,\ldots, s_l-1,
\quad
\Omega\subset \{z\colon \nu^{0l}\gamma\ne 0\}.
$
\hfill {\rm (2.37)}
\\[2.25ex]
}
\indent
{\sl Proof}. 
\vspace{0.5ex}
The system (2.37) has the determinant $(\nu^{0l}\gamma)^{s_{l}-1}.$ Therefore
there exists the solution $\Psi_{\eta l}^{\zeta},\ \eta=1,\ldots, s_l-1$ 
on a domain $\Omega\subset \{z\colon \nu^{0l}\gamma\ne 0\}$ 
\vspace{0.75ex}
of system (2.37).  

The proof of the lemma is by induction on $s_l.$
\vspace{0.5ex}

For $s_l=2$ and $s_l=3,$ the assertion (2.36) follows from (2.34).
\vspace{0.35ex}

Assume that (2.36) for $s_l=\varepsilon$ is true. Using (2.34) and (2.37), we get
\\[1.75ex]
\mbox{}\hfill
$
\displaystyle
{\frak L}_{\zeta}\,\nu^{\varepsilon l}\gamma =
\lambda_{l}^{\zeta}\sum\limits_{\delta=1}^{\varepsilon}
{\textstyle\binom{\varepsilon -1}{\delta-1}}\,\Psi_{\delta l}^{\zeta}(z)\,
\nu^{\varepsilon-\delta{,}\,l}\gamma +  
(\varepsilon -1)\sum\limits_{\delta=1}^{\varepsilon -1}
{\textstyle\binom{\varepsilon -2}{\delta-1}}\,\Psi_{\delta l}^{\zeta}(z)\,
\nu^{\varepsilon-\delta-1{,}\,l}\gamma  \ +
\hfill
$
\\[2.25ex]
\mbox{}\hfill
$
\displaystyle
+\ \nu^{\varepsilon -1{,}\,l}\gamma \,+\,
\nu^{0l}\gamma\ {\frak L}_{\zeta}\Psi_{\varepsilon l}^{\zeta}(z)$
\ \ for all $z\in \Omega.
\hfill
$
\\[2.5ex]
Combining (2.37) for 
\vspace{0.35ex}
$\eta=\varepsilon\!-\!1$ and $\eta=\varepsilon,\!$ (2.34) for $\eta=\varepsilon,$ and 
$\nu^{0l}\gamma \not\equiv 0\!$ in $\!\C^n\!,\!$ we obtain 
\\[1.5ex]
\mbox{}\hfill
$
{\frak L}_{\zeta}\,\Psi_{\varepsilon l}^{\zeta}(z)=0$ 
\ for all $z\in\Omega.
\hfill
$
\\[1.75ex]
\indent
So by the principle of mathematical induction, the statement (2.36) is true for every 
natural number $s_{l}\geq 2$ and $\zeta\in \{1,\ldots, m\}.$ 
\vspace{0.5ex}

Taking into account (2.32) and (2.35), we have the statement (2.36) 
\vspace{1ex}
is true for $j\ne \zeta.\!$ \k

{\bf Theorem 2.11.}
{\it
\vspace{0.5ex}
Let the assumptions of Lemma {\rm 2.2} with 
$l=1,\ldots r\ \Bigl(\, \sum\limits_{l=1}^{r}s_{l}\geq m+1\Bigr)$ hold.
Then the jacobian system {\rm (2.25)} has the ${\mathbb R}\!$-differentiable first integral
\\[2ex]
\mbox{}\hfill                                                                  
$
\displaystyle
F\colon z\to
\prod\limits_{\xi=1}^{\alpha}\bigl(\nu^{0\xi} \gamma\bigr)^{h_{0 \xi}}
\exp\sum\limits_{q=1}^{\varepsilon_{\xi}}\,h_{q\xi} \Psi_{q\xi}^{\zeta}(z)$
\ \ for all $z\in \Omega,
\quad \Omega\subset {\rm D}(F),
$
\hfill {\rm (2.38)}
\\[2ex]
where 
\vspace{0.5ex}
$\sum\limits_{\xi=1}^{\alpha}\varepsilon_{\xi}=m-\alpha+1, \
\varepsilon_{\xi}\leq s_{\xi}\!-\!1,\, \xi=1,\ldots,\alpha,\, \alpha\leq r,\!$ and 
$h_{q\xi},\, q\!=\!0,\ldots,\varepsilon_{\xi},\, \xi\!=\!1,\ldots, \alpha$ is a nontrivial solution to
the linear homogeneous system of equations 
\\[1.5ex]
\mbox{}\hfill
$
\displaystyle
\sum\limits_{\xi=1}^{\alpha}
\bigl(\lambda_{\xi}^{j}\,h_{0\xi}
+ \sum\limits_{q=1}^{\varepsilon_{\xi}}
\mu_{q\xi}^{j\zeta}\,h_{q\xi}\big) = 0,
\ \ \ j=1,\ldots, m.
\hfill
$
\\[2ex]
}
\indent
{\sl Proof}.
The Lie derivative of the function (2.38) by virtue of (2.25) is equal to 
\\[1.5ex]
\mbox{}\hfill
$
\displaystyle
{\frak L}_{j}\,F(z) =
\sum\limits_{\xi=1}^{\alpha}
\bigl(\lambda_{\xi}^{j}h_{0\xi} +
\sum\limits_{q=1}^{\varepsilon_{\xi}}\mu_{q\xi}^{j\zeta}h_{q\xi}\bigr)F(z)$
\ \ for all $z\in\Omega,
\quad 
j=1,\ldots, m.
\hfill
$
\\[2ex]
\indent
If
\vspace{0.5ex}
$
\sum\limits_{\xi=1}^{\alpha}
\bigl(\lambda_{\xi}^{j}h_{0\xi} +
\sum\limits_{q=1}^{\varepsilon_{\xi}}
\mu_{q\xi}^{j\zeta}h_{q\xi}\bigr) = 0,\ j=1,\ldots, m,$ then 
the ${\mathbb R}\!$-differentiable function (2.38) is a first integral of the 
\vspace{1.25ex}
jacobian system of partial differential equations (2.25). \k

As an example, the jacobian linear homogeneous system of partial differential equations
\\[1.5ex]
\mbox{}\ \ \                                        
$
z_2\,\partial_{{}_{\scriptstyle z_1}}w+(2z_2-\overline{z}_1-\overline{z}_2)\,\partial_{{}_{\scriptstyle z_2}}w+
(z_1-\overline{z}_2)\,\partial_{{}_{\scriptstyle \overline{z}_1}}w+
({}-z_1+2\overline{z}_1+2\overline{z}_2)\,\partial_{{}_{\scriptstyle \overline{z}_2}}w=0,
\hfill
$
\\[-0.1ex]
\mbox{}\hfill (2.39)
\\[0.75ex]
$
(2z_1-\overline{z}_1)\,\partial_{{}_{\scriptstyle z_1}}w+
({}-z_1+2z_2+\overline{z}_2)\,\partial_{{}_{\scriptstyle z_2}}w+
({}-z_1+3\overline{z}_1+\overline{z}_2)\,\partial_{{}_{\scriptstyle \overline{z}_1}}w+
(z_2-3\overline{z}_1+\overline{z}_2)\,\partial_{{}_{\scriptstyle \overline{z}_2}}w=0
\hfill
$
\\[2ex]
has the eigenvalue 
\vspace{0.25ex}
$\lambda_1^1=1$ with elementary divisor $(\lambda^1-1)^4$ 
corresponding to the  eigenvector $\nu^{0}=({}-1,1,{}-1,0)$ and to the generalized eigenvectors 
\vspace{0.5ex}
$\nu^{1}=(1,0,{}-1,{}-1),\ \nu^{2}=(1,{}-1,3,0),$ $\nu^{3}=({}-3,0,9,9).$  
The $\!{\mathbb R}\!$-differentiable functions (see the functional system (2.37))
\\[2ex]
\mbox{}\hfill                                                                 
$
\Psi_{11}^1\colon (z_1,z_2)\to
\dfrac{z_1-\overline{z}_1-\overline{z}_2}{{}-z_1+z_2-\overline{z}_1}$
\ \ for all $(z_1,z_2)\in \Omega,
\hfill
$
\\[2.75ex]
\mbox{}\hfill
$
\Psi_{21}^1\colon (z_1,z_2)\to
\dfrac{({}-z_1+z_2-\overline{z}_1)(z_1-z_2+3\overline{z}_1)-(z_1-\overline{z}_1-\overline{z}_2)^2}
{({}-z_1+z_2-\overline{z}_1)^{2}}$
\  \ for all $(z_1,z_2)\in \Omega,
\hfill
$
\\[2.75ex]
\mbox{}\hfill
$
\Psi_{31}^1\colon (z_1,z_2)\to
\dfrac{1}{({}-z_1+z_2-\overline{z}_1)^{3}}\,
\bigl(({}-3z_1+9\overline{z}_1+9\overline{z}_2)({}-z_1+z_2-\overline{z}_1)^2\, -
\hfill
$
\\[2.5ex]
\mbox{}\hfill
$
-\,3({}-z_1+z_2-\overline{z}_1)(z_1-\overline{z}_1-\overline{z}_2)(z_1-z_2+3\overline{z}_1) +
2(z_1-\overline{z}_1-\overline{z}_2)^3\,\bigr)$
\ for all $(z_1,z_2)\in \Omega,
\hfill
$
\\[2ex]
where a domain $\Omega\subset \{ (z_1,z_2)\colon z_1-z_2+\overline{z}_1\ne 0\}.$ 
\vspace{0.75ex}

The first integrals (by Theorem 2.11) of  the jacobian system (2.39) 
\\[1.75ex]
\mbox{}\hfill                                    
$
F_{1}\colon (z_1,z_2)\to \Psi_{21}^{1}(z_1,z_2)$
\ for all $(z_1,z_2)\in \Omega
\hfill                                    
$
\\[0.5ex]
and
\\[1.5ex]
\mbox{}\hfill                                    
$
F_{2}\colon (z_1,z_2)\to ({}-z_{1}+z_{2}-\overline{z}_{1})^{2}
\exp\bigl( {}-2\Psi_{11}^{1}(z_1,z_2)-\Psi_{31}^{1}(z_1,z_2)\bigr)$
\ for all $(z_1,z_2)\in \Omega
\hfill
$
\\[2ex]
are a basis of first integrals on the domain $\Omega$ of system (2.39).
\vspace{7ex}



}
\end{document}